\documentclass[3p,times,english,sort&compress]{elsarticle}
\usepackage{graphicx} 
\usepackage{amsthm}
\usepackage{amssymb}
\usepackage{amsmath}
\usepackage{latexsym}
\usepackage{bm}
\usepackage{color} 
\usepackage[final]{changes}
\usepackage{float}
\usepackage{caption}
\usepackage{subcaption}
\usepackage{placeins}

\everymath{
    \fontdimen14\textfont2=1.0ex
    \fontdimen17\textfont2=0.3ex
}

\usepackage[unicode=true,
bookmarks=true,bookmarksnumbered=false,bookmarksopen=false,
breaklinks=false,pdfborder={0 0 1},backref=false]
{hyperref}

\hypersetup{
    unicode=true,
    colorlinks,
    linkcolor={dark-red},
    citecolor={dark-blue},
    urlcolor={medium-blue},
}
\makeatletter

\pdfpageheight\paperheight
\pdfpagewidth\paperwidth

\begin{document}

\title{{\textit{e}{MOOD}}: Entropy-Based Multidimensional Optimal Order Detection for High-Order Shock-Capturing Finite-Volume Schemes}

\author[CU]{Micha{\l} Umi\'nski}
\author[CU]{Panagiotis Tsoutsanis\corref{cor1}}
\ead{panagiotis.tsoutsanis@cranfield.ac.uk}
\author[CU]{Ahmed A Sheikh Al-Shabab}
\author[AWE]{Alan S. Dawes}

\cortext[cor1]{Corresponding author}

\address[CU]{Faculty of Engineering and Applied Sciences, Cranfield University, Cranfield MK43 0AL, United Kingdom}
\address[AWE]{Computational Physics Group, AWE Nuclear Security Technologies, Aldermaston, Berkshire, United Kingdom}

\begin{abstract}
In this work, we present eMOOD, an entropy-guided evolution of the Multidimensional Optimal Order Detection (MOOD) framework for high-order finite-volume methods on unstructured meshes \cite{MOOD_OriginalPaper1, MOOD_OriginalPaper3, TsoutsanisMood}. The proposed approach introduces a novel troubled-cell indicator based on numerical entropy production, replacing the classical relaxed discrete maximum principle traditionally used in finite-volume MOOD schemes. By exploiting entropy generation as a physically motivated measure of non-smoothness, the new detector selectively targets shocks, contact discontinuities, and under-resolved flow features using a substantially more compact stencil.
Compared to conventional MOOD detection strategies, the proposed method reduces the number of flagged cells by up to 75\%, leading to improved computational efficiency, reduced memory-access overhead, and enhanced suitability for massively parallel simulations. The entropy-based detector is employed both within the classical MOOD cascade of progressively reduced-order schemes and in a novel MOOD-CWENO hybrid formulation, where inexpensive linear reconstruction is retained in smooth regions while nonlinear high-order reconstruction is activated only where required.
The methodology is implemented in the open-source parallel unstructured solver UCNS3D and assessed on a diverse set of challenging compressible-flow benchmarks, including the Shu-Osher shock--entropy-wave interaction, double Mach reflection, supersonic forward facing step, shock--helium-bubble interaction, and triple material point problems. These test cases span strong shock interactions, shock-induced instabilities, multi-material interfaces, and complex vortex-dominated flows, providing a rigorous assessment of robustness, shock-capturing capability, and accuracy on both structured and fully unstructured meshes.
The results demonstrate that eMOOD preserves the accuracy and robustness of reference CWENO and MOOD formulations while reducing computational cost, achieving speed-ups of up to 20\% without loss of solution quality. Owing to its compact detection stencil, reduced memory footprint, and improved localisation properties, the proposed framework provides a more scalable and computationally efficient approach for massively parallel high-order compressible-flow simulations. UK Ministry of Defence \copyright\  Crown owned
copyright 2026/AWE.
\end{abstract}

\maketitle

\section{Introduction}
Computational Fluid Dynamics (CFD) simulations are essential tools in engineering applications involving fluid flow, heat transfer, particulate transport, and other related physical phenomena. One of the major challenges in CFD is the accurate resolution of shock waves and discontinuities. High-order numerical methods are a common approach to this problem because they provide rapid error reduction with increasing resolution and are capable of accurately capturing complex nonlinear flow features. However, Godunov's theorem \cite{GodunovPhD,GodunovMethod_OriginalPaper} states that there are no second- or higher-order schemes, which are monotonic. The idea was further extended by Harten \cite{Harten1983} who showed that any linear scheme of second order or higher cannot be Total Variation Diminishing (TVD). In practice this means that high-order linear schemes produce spurious oscillations near discontinuities.

One of the first successful approaches to mitigate these oscillations was the Essentially Non-Oscillatory (ENO) method introduced by Harten et al. \cite{ENO_OriginalPaper}, which performs polynomial reconstructions over multiple stencils and selects the smoothest candidate. ENO was later largely superseded by the Weighted Essentially Non-Oscillatory (WENO) method \cite{WENO_OriginalPaper}, which blends candidate polynomials instead of selecting a single stencil.

The original derivation of WENO by Liu and Osher relied on uniform Cartesian meshes. Jiang and Shu later reformulated the smoothness indicators as integrals (WENO-JS) \cite{WENO-JS_OriginalPaper}, enabling the extension of the WENO schemes to arbitrary unstructured meshes. This generalisation was pursued, among others, by Tsoutsanis et al. \cite{Titarev2010,Tsoutsanis2011} within the massively parallel unstructured CFD framework UCNS3D \cite{UCNS3D_Paper}, as well as through alternative reconstruction bases such as Hermite polynomials \cite{HWENO_OriginalPaper1,HWENO_OriginalPaper2}.

Since then, WENO-based schemes have become one of the dominant approaches for high-order finite-volume (FV) methods and have undergone numerous developments and modifications, as reviewed by Bozorgpour et al. \cite{WENO_Review2025}. Two major areas of ongoing optimisation are the accuracy of the reconstruction near discontinuities and the high computational cost associated with multiple polynomial reconstructions and smoothness-indicator evaluations.

The first issue has been addressed by methods such as Mapped WENO \cite{Henrick2005} and WENO-Z \cite{BORGES20083191,Castro2011}, which modify the nonlinear blending weights to improve accuracy near shocks and discontinuities. Both issues are partially addressed by Central WENO (CWENO) schemes \cite{CWENO_OriginalPaper,Dumbser_2017,Tsoutsanis_Dumbser_CWENO}, which replace large directional stencils with smaller lower-order reconstructions and compensate through a modified central stencil. As a result, CWENO schemes reduce the computational cost while also decreasing the overall stencil size \cite{Tsoutsanis_Adebayo_Skote_2021}, while they can make switching from Discontinuous Galerkin to Finite-Volume more compact \cite{Maltsev_Skote_Tsoutsanis_dg_fv}. This idea was further extended in the Targeted ENO (TENO) framework \cite{TENO_OriginalPaper,TENO_Review2023}, which employs aggressive ENO-like stencil selection and sharper nonlinear weighting.

A large body of work has also focused on simplifying the computation of WENO smoothness indicators. Many approaches exploit properties specific to structured Cartesian meshes \cite{USVevek2021,HUANG2018498,KUMAR20181059}. More generally applicable approaches include the work of Zhang et al. \cite{Zhang2007}, who reduced the smoothness indicator to only two derivatives, as well as methods replacing derivative evaluations with simple difference-based measures \cite{Baeza2019,Pandey2022,Huang2023}. Furthermore, it was observed \cite{Teng2011,Wang2025} that the expensive nonlinear WENO reconstruction can often be avoided by employing only the central linear polynomial in sufficiently smooth regions and activating the full WENO procedure only when required.

In contrast to WENO, Clain, Diot, and Loub\`ere proposed the Multidimensional Optimal Order Detection (MOOD) method \cite{Diot2013,MOOD_OriginalPaper1,MOOD_OriginalPaper2,MOOD_OriginalPaper3}, which takes a fundamentally different approach to controlling oscillations in high-order schemes. Instead of reconstructing multiple candidate polynomials, MOOD first advances the solution using a high-order linear scheme, detects numerical oscillations a posteriori, and recomputes only troubled cells using progressively more robust lower-order schemes. Compared to WENO methods, MOOD eliminates the need for multiple reconstructions and smoothness-indicator evaluations while explicitly identifying numerical failures. In practice, only a small fraction of cells ($<5\%$) typically require re-computation, making the method computationally efficient.

MOOD has since been successfully applied to adaptive mesh refinement (AMR) \cite{Semplice2018}, unstructured meshes \cite{Farmakis2020,Tsoutsanis_ECCOMAS2022,Tsoutsanis2023}, and moving meshes \cite{Boscheri2015,Boscheri2022}, as well as to a wide range of governing equations including the Euler equations, ideal magnetohydrodynamics \cite{Gaburro2021}, shallow-water equations \cite{Ioriatti_PhD,Ioriatti2019}, and Arbitrary Lagrangian-Eulerian (ALE) remapping procedures \cite{Blanchard201683,Braeunig2018,Qing2023}. The method has also been incorporated into implicit schemes \cite{Michel-Dansac2022}, despite the increased complexity associated with troubled-cell detection in that context. Its scalability and robustness are further demonstrated by its implementation within the ExaHyPE framework \cite{ExaHyPE}, a massively parallel engine for solving hyperbolic PDEs on modern supercomputing architectures.

Beyond finite-volume methods, the MOOD paradigm has been extended to numerous other numerical frameworks. Applications include direct-ALE node motion \cite{Cirrottola2021}, THINC interface reconstruction \cite{Li2024}, finite-difference schemes \cite{Fernndez-Fidalgo2019, Bogey2004}, spectral methods \cite{Fernndez-Fidalgo2018,VelascoRomero2023}, residual distribution methods \cite{Bacigaluppi2023}, Gaussian-process reconstructions \cite{Bourgeois2022}, smoothed particle hydrodynamics \cite{Nogueira2016,Eirs2021}, and lattice Boltzmann methods \cite{Kozhanova2025}. However, its most widespread use has been as a limiter-activation strategy in Discontinuous Galerkin methods \cite{Vilar2019,Lee2020,Popov2021,Boscheri2017,Gaburro2021,Jiang2023,Nguyen2026,Ioriatti2019,Ioriatti2020}.

The original MOOD formulation \cite{MOOD_OriginalPaper1} employed a Discrete Maximum Principle (DMP) to detect oscillatory solutions in linear advection problems. When extended to the Euler equations \cite{MOOD_OriginalPaper3}, the criterion was relaxed heuristically to avoid incorrectly flagging smooth extrema. Since then, most finite-volume MOOD implementations have continued to rely on variants of the same detection strategy. In contrast, several non-FV implementations introduced alternative troubled-cell indicators, including smooth-extrema preserving criteria \cite{Bacigaluppi2023,VelascoRomero2023}, shock detectors \cite{Jiang2023}, and WENO-inspired smoothness indicators \cite{Nguyen2026}.

In this work, we present \textit{eMOOD}, an entropy-guided extension of the MOOD framework based on numerical entropy production. Unlike classical DMP-based MOOD detectors, the proposed criterion directly targets shocks, contact discontinuities, and under-resolved flow structures through a compact and physically motivated indicator. The proposed detector requires a significantly smaller stencil, improves localisation of troubled regions, and avoids the need for special treatment of smooth extrema. Since entropy production is strongest near discontinuities and shock waves \cite{Thornber2008}, the approach naturally focuses the MOOD recomputation process on regions where nonlinear stabilisation is genuinely required. Although entropy-based order adaptation was previously explored by Zappa and Semplice \cite{SharkEntropyPaper}, their work was restricted to one-dimensional structured grids, leaving its applicability to multidimensional unstructured meshes largely unexplored.

The proposed methodology is implemented in the open-source massively parallel UCNS3D CFD framework \cite{UCNS3D_Paper} and assessed on unstructured two-dimensional meshes using both the classical MOOD cascade and a novel MOOD-CWENO hybrid strategy inspired by the recent work of Tong and Tsoutsanis \cite{Tong2026}. The hybrid formulation preserves inexpensive linear reconstruction in smooth regions while activating the full nonlinear CWENO reconstruction only where necessary.

The remainder of this paper is organised as follows. Section \ref{Numerical_Methods} introduces the governing numerical methods together with the proposed entropy-based troubled-cell indicator and the MOOD-CWENO hybrid formulation. Sections \ref{Applications} and \ref{General_Performance} present the benchmark test cases and analyse the performance, robustness, and computational efficiency of the proposed methods. Finally, conclusions and directions for future work are discussed in Section \ref{Conclusions}.

\section{Numerical Methods}\label{Numerical_Methods}
\subsection{Explicit Finite Volume Method}

Consider the unsteady nonlinear hyperbolic system of conservation laws on a 2D (or 3D) domain $\Omega$, written in its conservative form:
\begin{equation}
\frac{\partial{\mathbf{U}}}{\partial t}+\nabla  \cdot (\mathbf{F}({\mathbf{U}}))=0, \label{eq:Conservative_equation}
\end{equation}
where ${\bm{U}}={\bm{U}}(\bm{x},t)$ is the vector of conserved variables, $\bm{x}$ denotes the coordinates of a point in the domain, $t$ is time, and $\bm{F}({\bm{U}})
$ is a nonlinear flux tensor.
The physical domain $\Omega$ is divided into disjoint polygonal (polyhedral in 3D) cells $V_i$ indexed by $i$ - $\Omega = \bigcap\limits_iV_i$. Let $S_{ij} = V_i\cap V_j$ be the edge (face) between cells $i$ and $j$ and $N_i = \{j:V_i\cap V_j\ne\varnothing\}$ be the set of cells which share an edge (respectively face) with $i$.

The global solution is discretised and approximated by storing in each cell a vector of cell averages:
\begin{equation}
\bm{U}_i^t=\frac{1}{\left|V_{i}\right|}\int_{V_{i}}\bm{U}(\bm{x},t)\ dV
\label{eq:FV_average}
\end{equation}
Integrating (\ref{eq:Conservative_equation}) over a cell and applying Stoke's theorem one gets
\begin{equation}
 \frac{d\bm{U}_i}{dt}= -\frac{1}{|V_{i}|}\sum\limits _{j\in N_i}\int\limits_{S_{ij}}\bm{F}(\mathbf{U(\bm{x},t}))\cdot\mathbf{n}_{ij}\ d\bm{x}
\label{eq:FV_formulation}
\end{equation}
where $\mathbf{n}_{ij}$ is the unit vector normal to the edge (face) $S_{ij}$. This is further discretised using Gaussian quadrature with number of points $N_{qp}$ depending on desired order of accuracy:
\begin{equation}
 \frac{d\bm{U}_{i}}{dt}= -\frac{1}{|V_i|}\sum\limits _{j\in N_{i}}\sum\limits _{\alpha=1}\limits^{N_{qp}}\omega_{\alpha}\bm{F}_{\bm{n}_{ij}}\left({\bm{U}}^{n,L}_{ij}(\bm{x}_{ij,\alpha},t),{\bm{U}}^{n,R}_{ij}(\bm{x}_{ij,\alpha},t)\right)|S_{ij}|
\label{eq:FV_formulation_discretised}
\end{equation}
where $x_{i,j,\alpha}$ are Gaussian quadrature points on the edge (respectively face) $S_ij$ and $\omega_\alpha$ are corresponding weights. $\bm{F}_{bm{n}_{ij}}$ is an approximate numerical flux in the direction normal to the edge (face), which depends on the approximate field values $\bm{U}^{n,L}_{ij}$ and $\bm{U}^{n,R}_{ij}$ at the quadrature point in the cells $i$ and $j$ respectively.

There are multiple ways of constructing numerical flux. In this work the Godunov's \cite{GodunovMethod_OriginalPaper} method with approximate HLL (Harten-Lax-van Leer) \cite{HLL_OriginalPaper} and HLLC (Harten-Lax-van Leer-Contact) of Toro \cite{HLLC_OriginalPaper} Riemann solvers is employed.
The HLL flux assumes that the Riemann problem produces two shock waves travelling with speeds $s_L<s_R$ and approximates the flux as:
\begin{equation}
    \bm{F}^{HLL}_{\bm{n}_{ij}}(\bm{U}^L, \bm{U}^R) = 
    \begin{cases}
    \bm{F}(\bm{U}^L) & if\ s_L \le 0     \\
    \bm{F}^\star     & if\ s_L > 0 > S_R \\
    \bm{F}(\bm{U}^R) & if\ s_R \ge 0     \\
    \end{cases}
    \label{eq:HLL_flux}
\end{equation}
where $\bm{F}^\star$ is found by considering the Rankine–Hugoniot conditions across both waves and solving a linear equation. The HLLC flux adds a contact discontinuity between the shock waves and follows an equivalent derivation with four cases and two intermediate fluxes $\bm{F}^{*L}$, $\bm{F}^{*R}$. For details see Toro's textbook \cite{ToroTextbook}.

\subsubsection{High-order in time}\label{Runge-Kutta_description}
Desired order of accuracy in time is achieved by integrating (\ref{eq:FV_formulation_discretised}) using an explicit Runge-Kutta (RK) method of appropriate order with the time step $\Delta t$ is selected according to the Courant-Friedrichs-Lewy condition 
\begin{equation}
    \Delta t = \frac{\text{CFL}}{2*dimension+1}\min_i\left(\frac{h_i}{s_i}\right),
\end{equation}
where $h_i$ is the characteristic length of the cell $i$, $s_i \approx |\bm{u}_i|+c_i$ is an estimate of the maximum wave propagation speed in the cell. The CFL number is adjusted depending on the mesh quality, complexity of the test problem and robustness of the method. In the presented test problems, unless specified otherwise, the value $\text{CFL}=0.75$ and the third-order Strong Stability Preserving (SSP) RK method \cite{SSP_Runge-Kutta} are used:
\begin{equation}\label{eq:Runge-Kutta}
\begin{split}
    \mathbf{U}_i^{I} &= \mathbf{U}_i^t - \frac{\Delta t}{|V_i|} \bm{F}(\bm{U}^t)\\
    \mathbf{U}_i^{II} &= \frac{3}{4}\mathbf{U}_i^t + \frac{1}{4}(\mathbf{U}_i^{I} - \frac{\Delta t}{|V_i|} \bm{F}(\mathbf{U}^{I})) \\
    \mathbf{U}_i^{t+\Delta t} &= \frac{1}{3}\mathbf{U}_i^t + \frac{2}{3}(\mathbf{U}_i^{II} - \frac{\Delta t}{|V_i|} \bm{F}(\mathbf{U}^{II}))
\end{split}
\end{equation}

\subsubsection{High-order in space}
Higher order in space is achieved by constructing an approximate polynomial distribution of $\textbf{U}$ in the cell $i$:
\begin{equation}
    \mathbf{U}_i(\mathbf x) \approx \mathbf{p}_i(\mathbf{x}) := \sum_{k=0}^K \mathbf{U}_i^{k}\phi_{i}^{k}(\mathbf x).
    \label{eq:polynomial_reconstruction}
\end{equation}
\noindent where $\phi_i(\mathbf{x})$ are polynomial basis functions and $\mathbf{U}_i^{k}$ are coefficients and interpolating $\mathbf{U}^{n,L}_{ij,L}(\mathbf{x}_{ij,\alpha},t) \approx p_i(x_{ij,\alpha})$ and $\mathbf{U}^{n,R}_{ij,L}(\mathbf{x}_{ij,\alpha},t) \approx p_j(x_{ij,\alpha})$ to be used in Equation (\ref{eq:FV_formulation_discretised}). 

The simplest approach for finding the approximate distribution $\bm{p}_i(\bm{x})$ would involve constructing a stencil (set) of at least $K$ direct and indirect neighbours around the cell $i$ and obtaining the polynomial weights by minimising the error $\left|V_j\right|\bm{U}_j - \int\limits_{V_j} \bm{p}_i(\bm{x})dV$ using the least squares method.
Using only one stencil results in a linear scheme and by Godunov's theorem \cite{GodunovPhD} such scheme cannot be Total Variation Diminishing (TVD), which in practice results in spurious oscillations in the solution.

\subsection{CWENO reconstruction}\label{CWENO-description}

The WENO method originally proposed by Liu et al. \cite{WENO_OriginalPaper} is a solution to this problem. In addition to the central stencil ($k=1$) a number of directional stencils ($k=2...s_t$), which include only neighbours on one side of the cell, are used. The final polynomial is a linear combination of the polynomials fitted to each stencil $p_i(x) = \sum\limits_{k=1}^{s_t}w^kp_i^k(x)$ where the weight $w^k$ depends on a smoothness indicator of the polynomial $p_i^k(x)$:
\begin{equation}\label{eq:SmoothnessIndicator}
\begin{split}
\mathcal{SI}_k = \sum\limits_{a,b:\ 1\leq a+b} &\int_{V_i} \left( \frac{\partial^{a+b}}{\partial x^a \partial y^b}p_i^k(x, y) \right)^2dV \\
\tau = \left(\frac{\sum\limits_{k=2}^{s_t}|\mathcal{SI}_k-\mathcal{SI}_1|}{{s_t}-1}\right)^4 \qquad \tilde{w}^k &= \lambda_k \left(1+
\frac{\tau}{\epsilon+\mathcal{SI}_k}\right) \qquad
w^k = \frac{\tilde{w}^k}{\sum\limits_{l=1}^{s_t}\tilde{w}^l}\\ 
\end{split}
\end{equation}
with a small constant $\epsilon = 10^{-14}$ added to avoid division by zero.

The original WENO scheme requires computing polynomial fits and smoothness indicators on a number of large stencils, which is computationally expensive. Compact WENO (CWENO) has been proposed \cite{CWENO_OriginalPaper} to address this issue by replacing the directional stencils with smaller ones and using first-order reconstruction polynomials on them. The central polynomial is modified
\begin{equation}\label{eq:CWENO_polynomail}
\begin{split}
p_{opt}(\mathbf x) &= \sum\limits _{s=1}^{s_t} \lambda_{s}p_{s}(\mathbf x) \\
p_{1}(\mathbf x) &= \frac{1}{\lambda _{1}} \left(p_{opt}(\mathbf x)-\sum\limits _{s=2}^{s_t} \lambda_{s}p_{s}(\mathbf x)\right).
\end{split}
\end{equation}
and the subsequent algorithm is identical as in WENO.
The method has been shown to be meaningfully faster while offering as good or better accuracy and robustness of the original WENO method \cite{Tsoutsanis_Adebayo_Skote_2021,Cravero2016}, so it is used in this work.
Compared to Cartesian meshes, on unstructured meshes there is more then one natural choice stencil size and cell selection and the decision has an effect of accuracy and robustness of the method. In this work, we use stencils of $2K$ calls for improved robustness and follow the stencil selection algorithm of Tsoutsanis et al. \cite{Tsoutsanis_stencil_selection}.

\subsection{MOOD}

\begin{figure}[!h]
    \vspace{-2ex}
    \centering
    \includegraphics{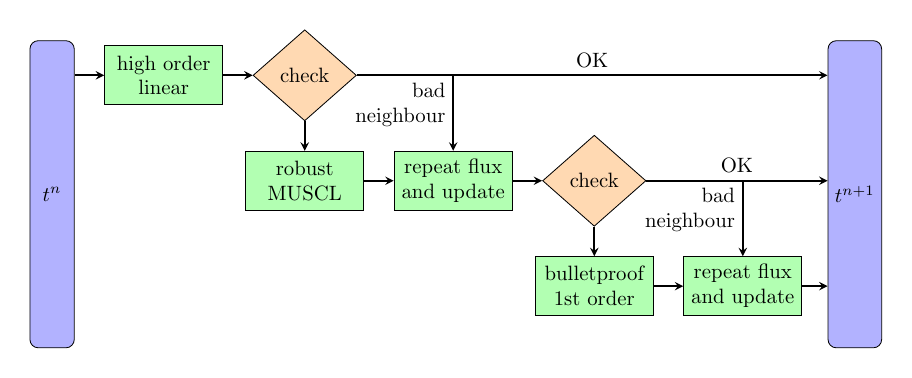}
    \vspace{-5ex}
    \caption{\label{fig:MOOD-diagram} Symbolic flow chart of the MOOD method.}
    \vspace{-1ex}
\end{figure}

A different solution to the limitations of linear high-order reconstruction is the MOOD method originally proposed by Clain, Diot, and Loub\'ere \cite{MOOD_OriginalPaper1, MOOD_OriginalPaper2}. It attempts to perform time step using linear scheme and checks each cell for errors and artefacts, if none are found the solution in the cell is accepted as final, otherwise the reconstruction in the cell is repeated using the $2^{nd}$ order MUSCL \cite{MUSCL_OriginalPaper} scheme or any variant of it \cite{Tsoutsanis_2018_bounds_limiter}, which is more robust, the fluxes are recomputed and the time step in the cell and its neighbours is repeated. The potential solution is checked again. If it fails, the step is repeated for the third time using the most robust first-order scheme. A simplified flow chart can be seen in the Figure \ref{fig:MOOD-diagram}.
The leading challenge in parallel implementation of MOOD is the need to repeat all the cell average updates contaminated by flux from inadmissible cells. This typically requires additional communication and synchronisation with very unpredictable and scattered pattern. As a results it is favourable to have as few cells flagged as possible and have them clustered in groups reducing the number of neighbours.

\subsection{MOOD-(C)WENO Hybrid}

\begin{figure}[!h]
    \vspace{-2ex}
    \centering
    \includegraphics{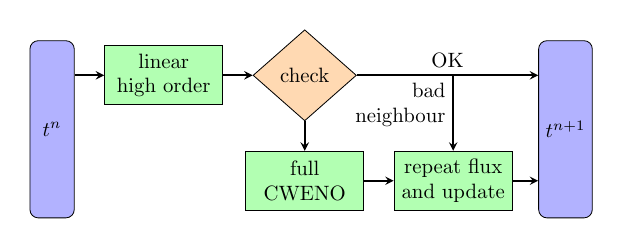}
    \vspace{-3ex}
    \caption{\label{fig:MOOD_hybrid-diagram} Symbolic flow chart of the hybrid MOOD-CWENO method. Compare with Figure \ref{fig:MOOD-diagram}}
    \vspace{-1ex}
\end{figure}

The MOOD method relays on $2^{nd}$ and $1^{st}$ order methods in the latter steps, which in principle can lead to reduced order of accuracy and difficulty with capturing higher-order phenomena. As a result, we investigate using in the second step the full (C)WENO method of the same order as the first step, as visualised in the Figure \ref{fig:MOOD_hybrid-diagram}.
The linear reconstruction from the first step is reused as the central stencil in CWENO, so in the second step it is only necessary to find polynomial reconstructions on the smaller directional stencils and compute all the weights. Similar to traditional MOOD, because linear scheme is sufficient in majority of cells this method should be computationally more efficient than CWENO, but by construction it has the same robustness and order of accuracy as traditional CWENO implementation. Compared to MOOD the method should maintain higher order of accuracy, and have lower parallelisation overhead, because the second round of re-updates and associated communication is eliminated.

\subsection{Switch Criteria}
In both traditional MOOD and the proposed MOOD-CWENO hybrid the switch for each element from a high-order linear solution to more robust algorithm is governed by a troubled cell indicator. To maintain the high order of convergence and to minimise the computational cost, the aim for the indicator is to trigger in as few cells as possible while still giving accurate solution. Secondarily, it should be quick to compute and preferably depend only on cell-local data or a limited number of cell neighbours, to simplify memory access patterns and communication overhead.  

\subsubsection{PAD \normalfont{\&} NAD indicator} \label{classic_MOOD_trigger}
For reference we use the criteria proposed in the original papers introducing MOOD paradigm \citep{MOOD_OriginalPaper3,MOOD_OriginalPaper1,Dumbser2014,Vilar2019}, namely the Physical Admissible Detector (PAD) and Numerical Admissible Detector (NAD) \citep{MOOD_OriginalPaper3}.
The PAD flags the cell as troubled if any component of the potential solution is a NaN (Not a Number), and when the density or pressure become negative. The NAD criteria checks that no new extrema are created according to the discrete maximum principle (DMP). The relaxed variant with added margin $\delta$ is used meaning that the cell is flagged if either density or pressure component of the potential $\bm{U}^{t+\Delta t,*}$ solution does not satisfy
\begin{equation}\label{eq:mood_pad}
\begin{split}
    \bm{U}^{t+\Delta t,*}_i \in &\left[\min_{j\in N_i} \bm{U}^{t}_j - \delta\ ,\ \max_{j\in N_i} \bm{U}^{t}_j + \delta\right]\\
    \delta =& \max\left(10^{-4},\alpha\cdot(\displaystyle{\max_{j\in N_i}}\bm{U}^{t}_j-\displaystyle{\min_{j\in N_i}}\bm{U}^{t}_j\right).
\end{split}
\end{equation}
The choice of the relaxation parameter $\alpha$ regulates the degree of oscillations allowed in the solution. In this work, the value of $\alpha=10^{-3}$ as suggested in \citep{Fidalgo2018} is used in the first pass and $\alpha=10^{-1}$ in the second pass. The condition $\delta \geq 10^{-4}$ is introduced to prevent the indicator from spuriously triggering in smooth regions where 
$\bm{U}(\bm{x}) \approx \min U^t_j \approx \max U^t_j$.

\subsubsection{Entropy Production Indicator}\label{entropy_MOOD_trigger}
As an alternative to PAD and NAD conditions, we propose to use the entropy production as a troubled cell indicator. Entropy production can be calculated as $\Delta s = |S(\bm{U}^{t+\Delta t,*}(x)) - S(\bm{U}^t(x))|$. Where $\bm{U}^{t+\Delta t,*}(x)$ is the potential solution in the cell and $\bm{U}^t(x)$ is the solution in the previous time step. We want to mark the cell as ``troubled'' if $\Delta s > \beta$. However, entropy production scales with mesh resolution and time step, to compensate for it we rescale 
$\Delta s_{rscl} = \Delta s / \bar{h}$ where $\bar{h}$ is the average characteristic cell size defined as $\bar{h} = \frac{1}{N}\sum\limits_{i=1}^{N}\sqrt{|V_i|}$. The final condition for flagging is:
\begin{equation}
    \left|S(\bm{U}^{t+\Delta t,*}_i) - S(\bm{U}^t_i)\right| \frac{1}{\bar{h}} > \beta
    \label{eq:entropy_MOOD_condition}
\end{equation}
The threshold $\beta$ takes two values smaller $\beta_1$ for the first pass and larger $\beta_2$ for the second pass, and is selected for each test case based on characteristic velocity, density, and pressure.
For the test cases with single ideal gas the entropy is calculated as
\begin{equation}
    S_{ideal}(\rho, p) = -\rho\ log\left(\frac{p}{(\gamma - 1)\ \rho^\gamma}\right) 
    \label{eq:idealGasEntropy}
\end{equation}
and for tests with two gases mixed in ratio $vf : (1-vf)$ the entropy is the mixture of the two ideal entropies:
\begin{equation}
    S_{mix}(\rho_1, \rho_2, p, vf) = vf \cdot S_{ideal}(\rho_1,p) + (1- vf)\cdot S_{ideal}(\rho_2,p)
    \label{eq:mixtureEntropy}
\end{equation}

\subsubsection{Entropy Production Scaling}
For the entropy production described in \ref{entropy_MOOD_trigger} to be a useful switching criterion, it must work reliably on unstructured meshes with a large range of cell sizes. Additionally, it would be desirable if the thresholds $\beta_{1,2}$ were independent of the mesh resolution so that they can be set quickly on a coarse mesh and reused for final simulations on refined meshes. Because numerical entropy production depends on the time step and mesh resolution \cite{SharkEntropyPaper} it is plausible that entropy production may need to be rescaled.
To optimise for the first requirement an experiment was conducted on a heavily deformed mesh with an extreme range of cell sizes.
\begin{figure}[H]
    \centering
    \vspace{-2ex}
    \includegraphics[width=\linewidth]{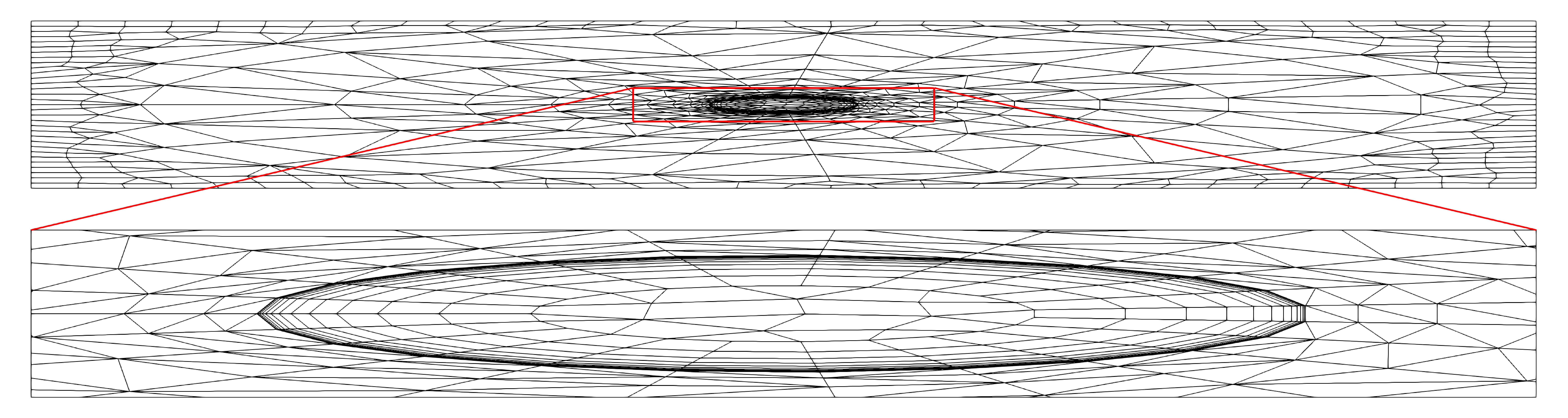}
    \vspace{-4ex}
    \caption{\label{fig:bad_mesh} Heavily deformed unstructured mesh used for testing the scaling in the entropy production MOOD trigger. Full mesh shown at the top, and the central region enlarger at the bottom.}
    \vspace{-2ex}
\end{figure}
The Shu-Osher test is used \ref{ShuOsher}, and five scaling strategies are considered: $\Delta s$, $\Delta s/h$, $\Delta s\cdot h$, $\Delta s/|V_i|$, $\Delta s\cdot |V_i|$, where $|V_i|$ is the cell's area and $h_i := \sqrt|V_i|$ is the characteristic length of the cell. For each strategy independently good values for thresholds $\beta_1$ and $\beta_2 = 5\beta_1$ were found and the simulation results were compared.
\begin{figure}[H]
    \centering
    \begin{subfigure}{\textwidth}
        \centering
        \includegraphics[width=0.75\linewidth]{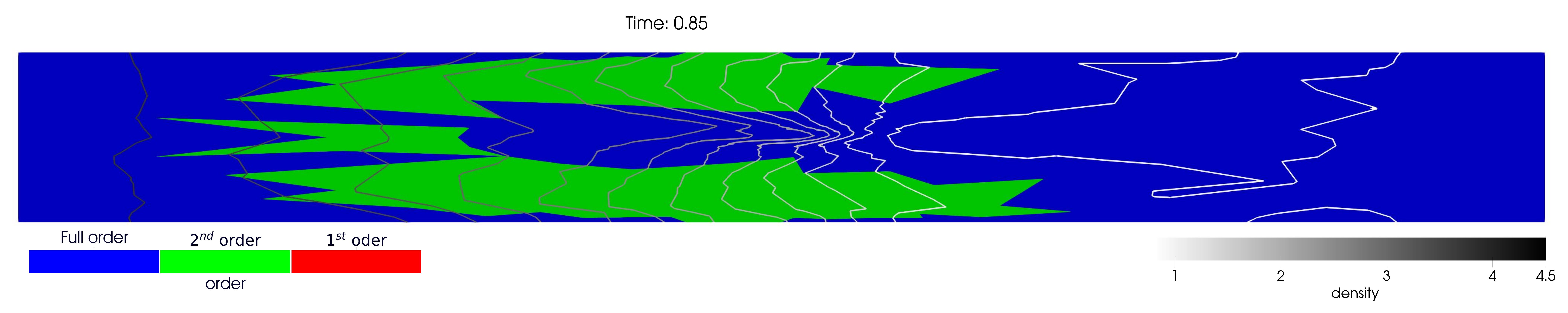}
        \vspace{-4ex}
        \caption{\label{fig:times_cell_area_badmesh_scaling}$\Delta s\cdot |V_i|$}
    \end{subfigure}
    \begin{subfigure}{\textwidth}
        \centering
        \includegraphics[width=0.75\linewidth]{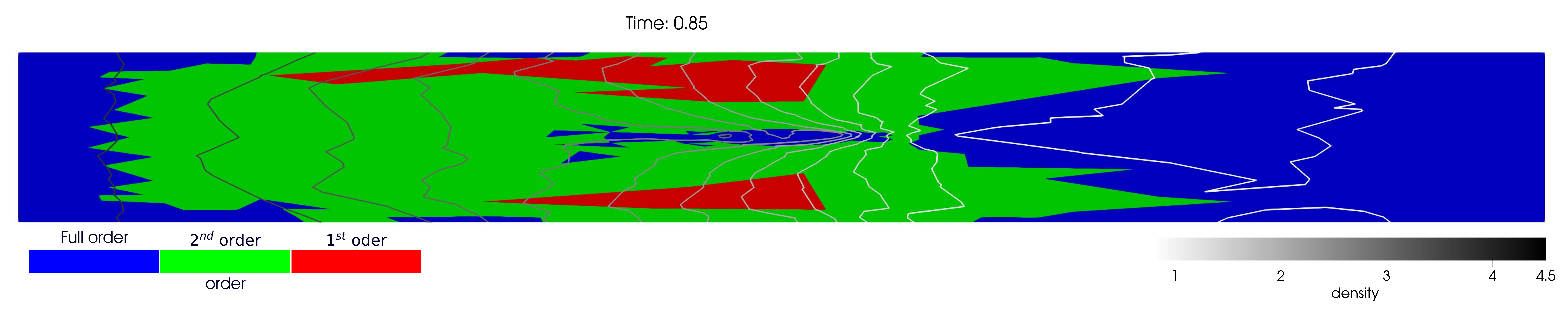}
        \vspace{-4ex}
        \caption{\label{fig:times_cell_size_badmesh_scaling}$\Delta s\cdot h_i$}
    \end{subfigure}
    \begin{subfigure}{\textwidth}
        \centering
        \includegraphics[width=0.75\linewidth]{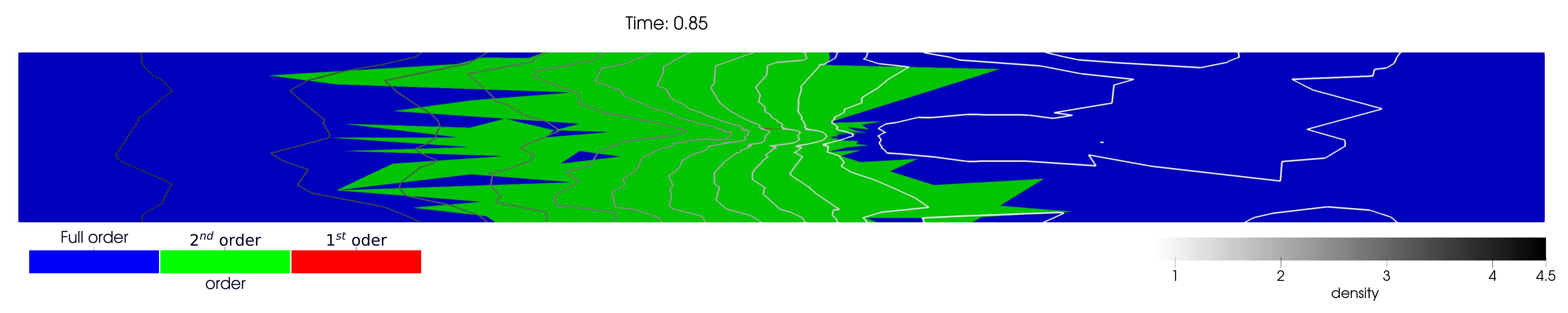}
        \vspace{-4ex}
        \caption{\label{fig:nothing_cell_badmesh_scaling}$\Delta s$}
    \end{subfigure}
    \begin{subfigure}{\textwidth}
        \centering
        \includegraphics[width=0.75\linewidth]{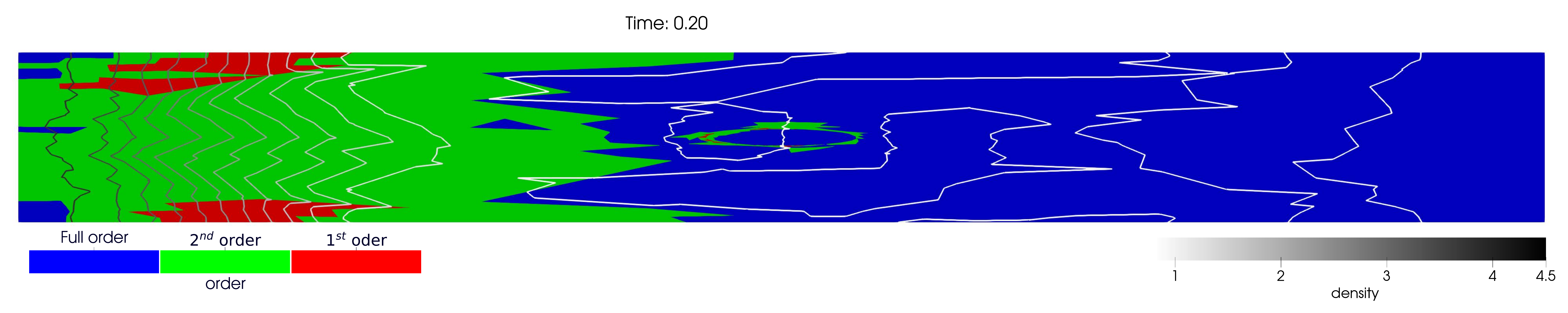}
        \vspace{-4ex}
        \caption{\label{fig:div_cell_size_badmesh_scaling}$\Delta s/h_i$}
    \end{subfigure}
    \begin{subfigure}{\textwidth}
        \centering
        \includegraphics[width=0.75\linewidth]{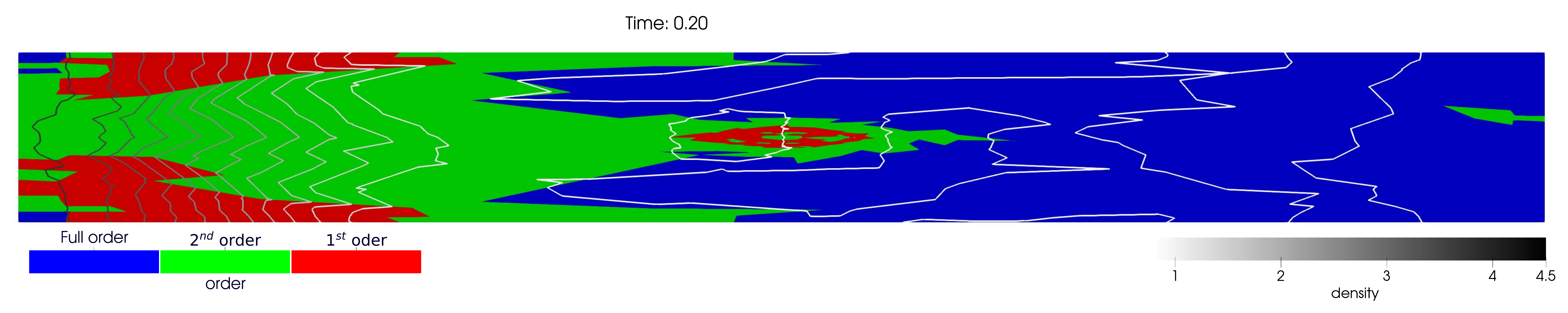}
        \vspace{-4ex}
        \caption{\label{fig:div_cell_area_badmesh_scaling}$\Delta s/|V_i|$}
    \end{subfigure}
    \caption{\label{fig:deformed_mesh_scaling} Comparison of different entropy production scalings in the MOOD trigger on a heavily deformed mesh (Figure \ref{fig:bad_mesh}). Green designates the cells which have been flagged in the first pass, red designates cells which have been flagged in both passes and hence have used $1^{st}$ order reconstruction.}
\end{figure}

From Figure \ref{fig:deformed_mesh_scaling} it is visible that multiplying the entropy production by cell area causes the trigger to never flag the smallest cells. On the other hand dividing the entropy production by cell area causes the trigger to flag small cells prematurely. Replacing the cell area with characteristic length in either scaling strategy reduces the magnitude of the flaws but does not eliminate them. The best scaling strategy out of five considered is \ref{fig:nothing_cell_badmesh_scaling} - no scaling with cell size.

To test the entropy production scaling with mesh refinement the Forward Facing Step problem was ran on a sequence of four unstructured meshes each doubling the resolution. Unstructured meshes do not have one unambiguously defined resolution, as a result average cell size $\overline{h} := \frac{1}{N}\sum\sqrt{|V_i|}$ and average cell area $\overline{A} := \frac{1}{N}\sum|V_i|$ are used as proxies.
They are selected because they scale linearly and quadraticly respectively with refinement when the whole mesh is refined. Selecting optimal resolution proxies for meshes where refinement has been applied only locally is left for future work.
Four scaling strategies were considered: $\Delta s \cdot \overline{A}$, $\Delta s \cdot \overline{h}$, $\Delta s$, and $\Delta s/ \overline{h}$.
The thresholds $\beta_{1,2}$ are adjusted to produce identical results on the coarse mesh and results on the finer meshes are compared.

\begin{figure}[H]
    \centering
    \begin{subfigure}{0.245\textwidth}
        \centering
        \includegraphics[width=\linewidth]{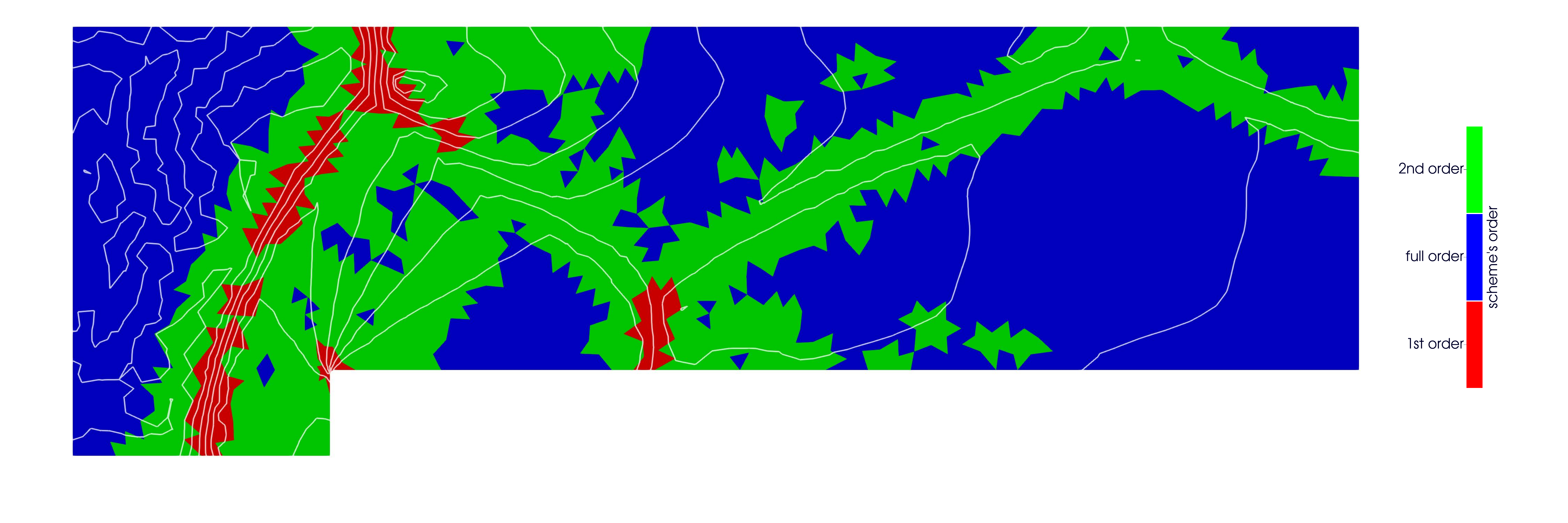}
        \includegraphics[width=\linewidth]{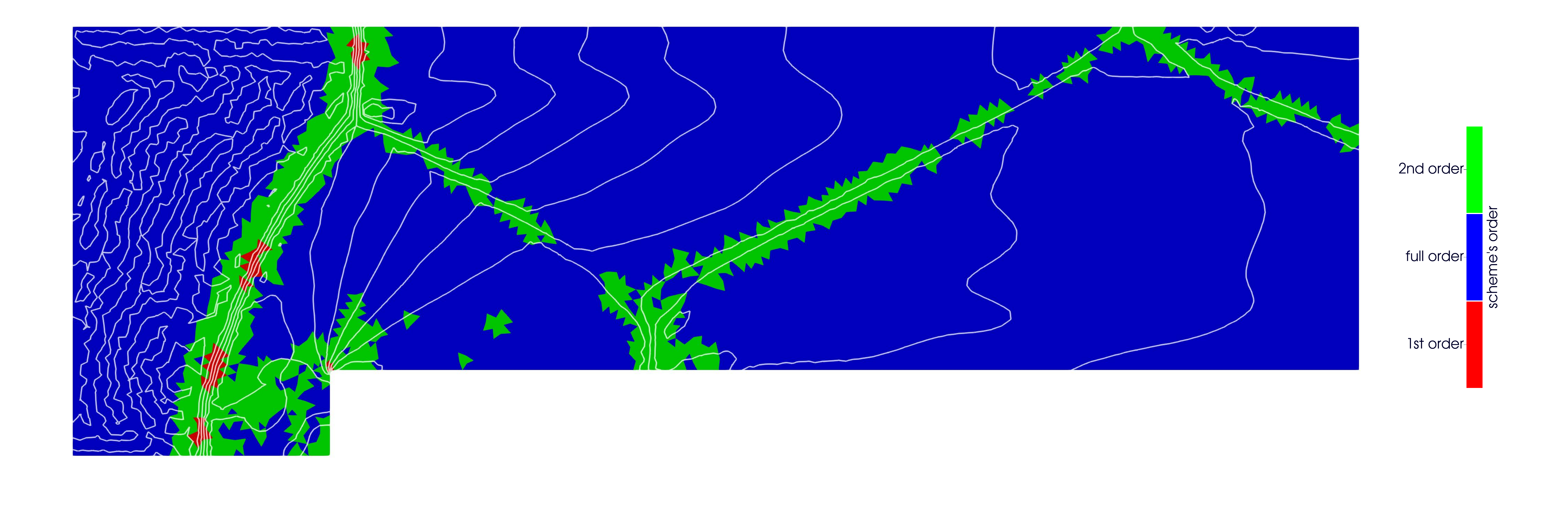}
        \includegraphics[width=\linewidth]{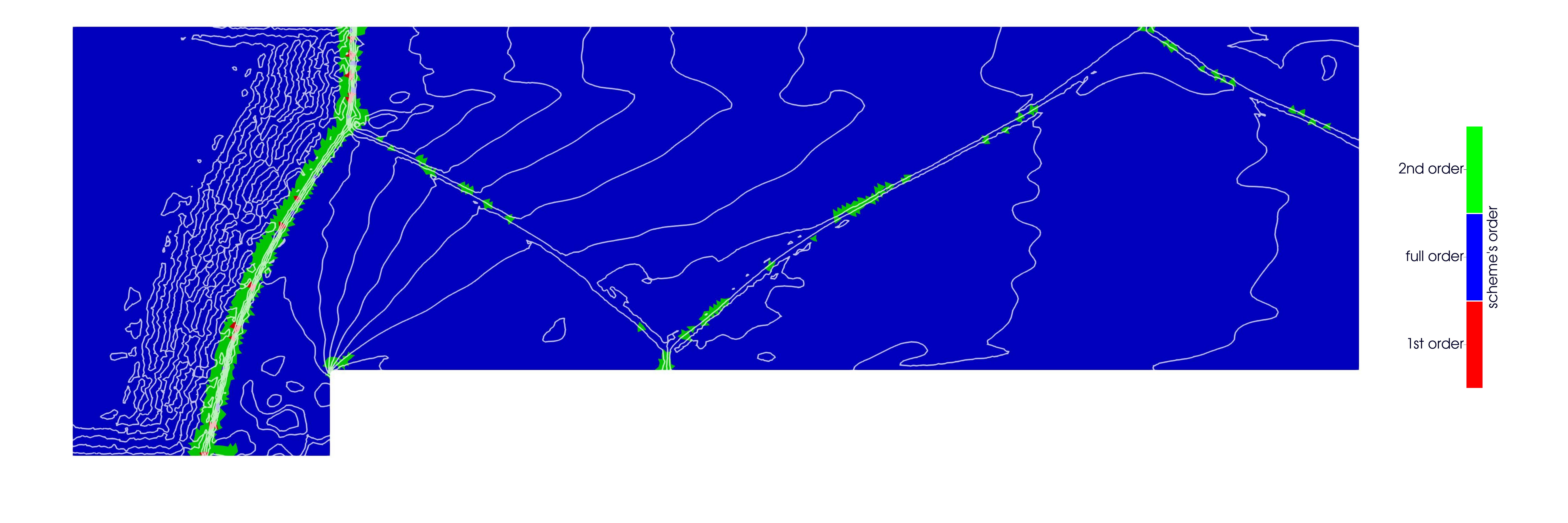}
        \includegraphics[width=\linewidth]{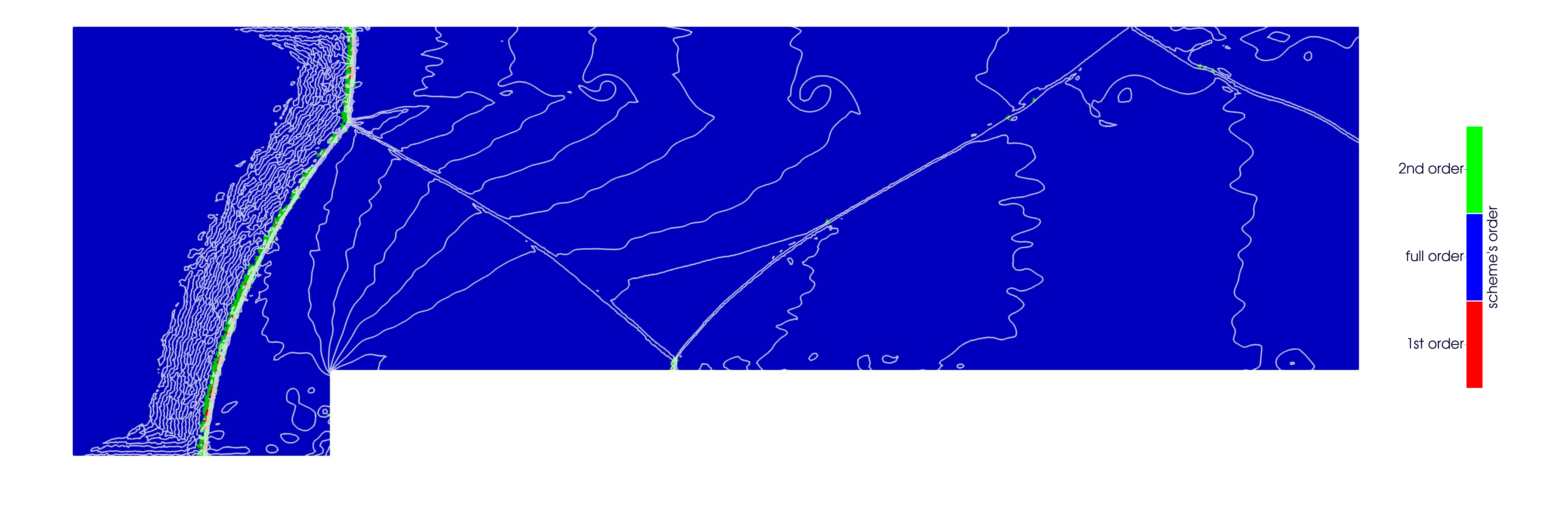}
        \vspace{-6ex}
        \caption{\label{fig:times_avg_cell_area_scaling}$\Delta s \cdot \overline{A}$}
    \end{subfigure}
    \begin{subfigure}{0.245\textwidth}
        \centering
        \includegraphics[width=\linewidth]{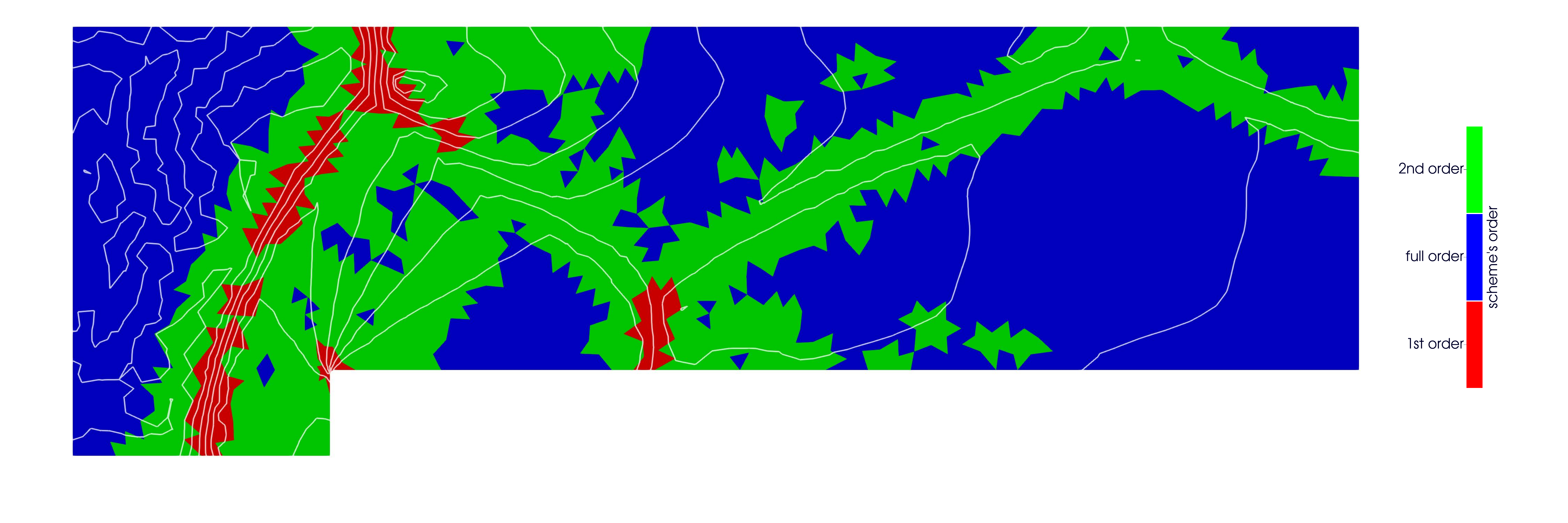}
        \includegraphics[width=\linewidth]{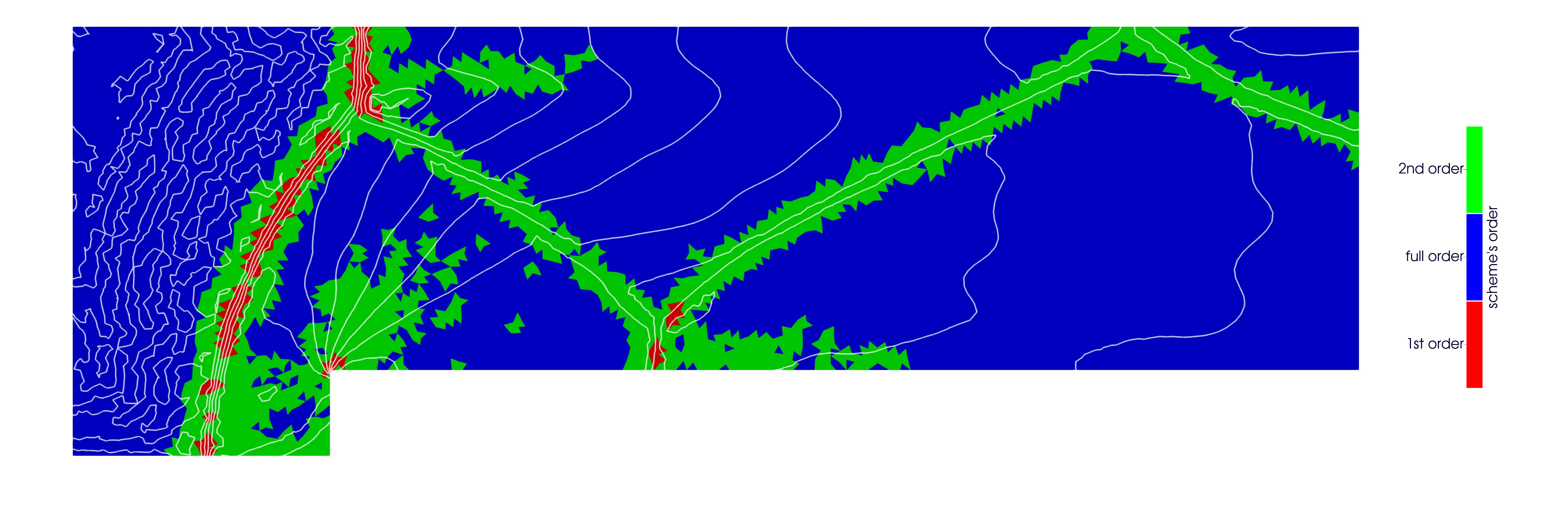}
        \includegraphics[width=\linewidth]{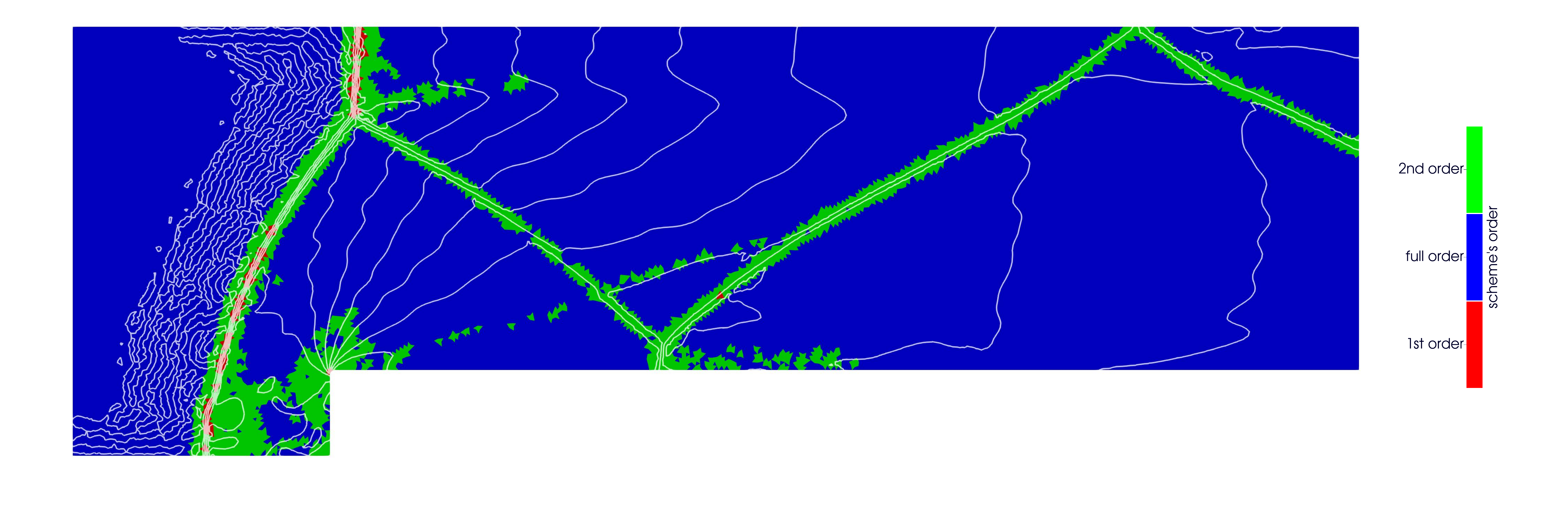}
        \includegraphics[width=\linewidth]{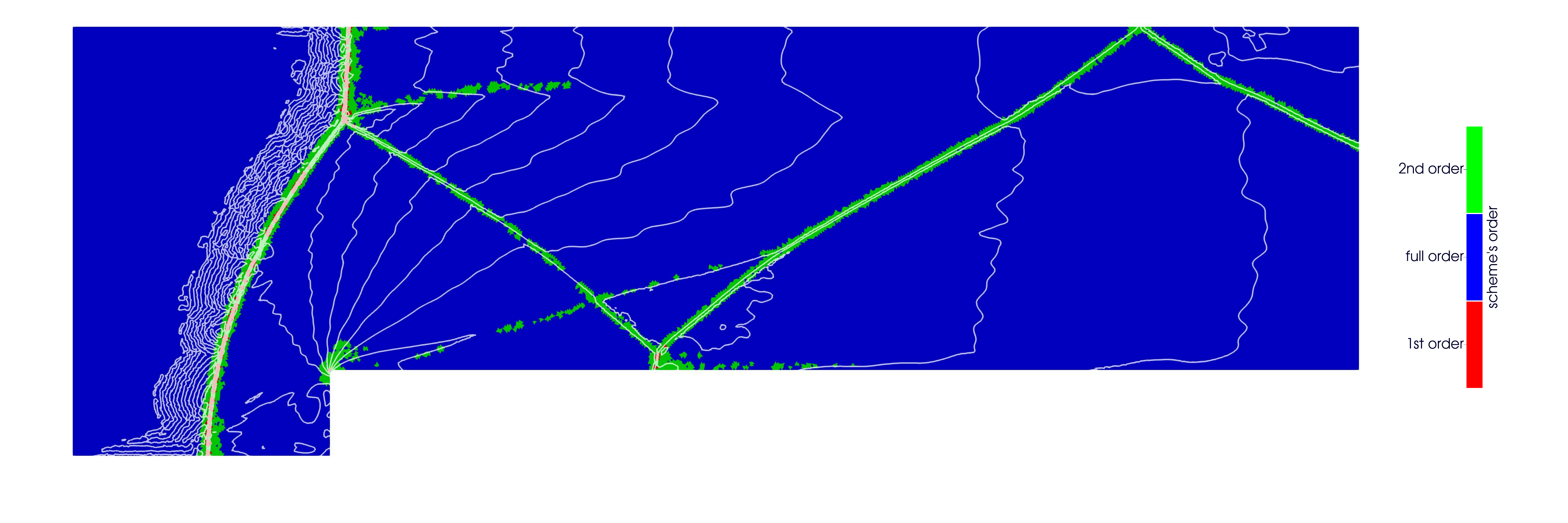}
        \vspace{-6ex}
        \caption{\label{fig:times_avg_cell_size_scaling}$\Delta s \cdot \overline{h}$}
    \end{subfigure}
    \begin{subfigure}{0.245\textwidth}
        \centering
        \includegraphics[width=\linewidth]{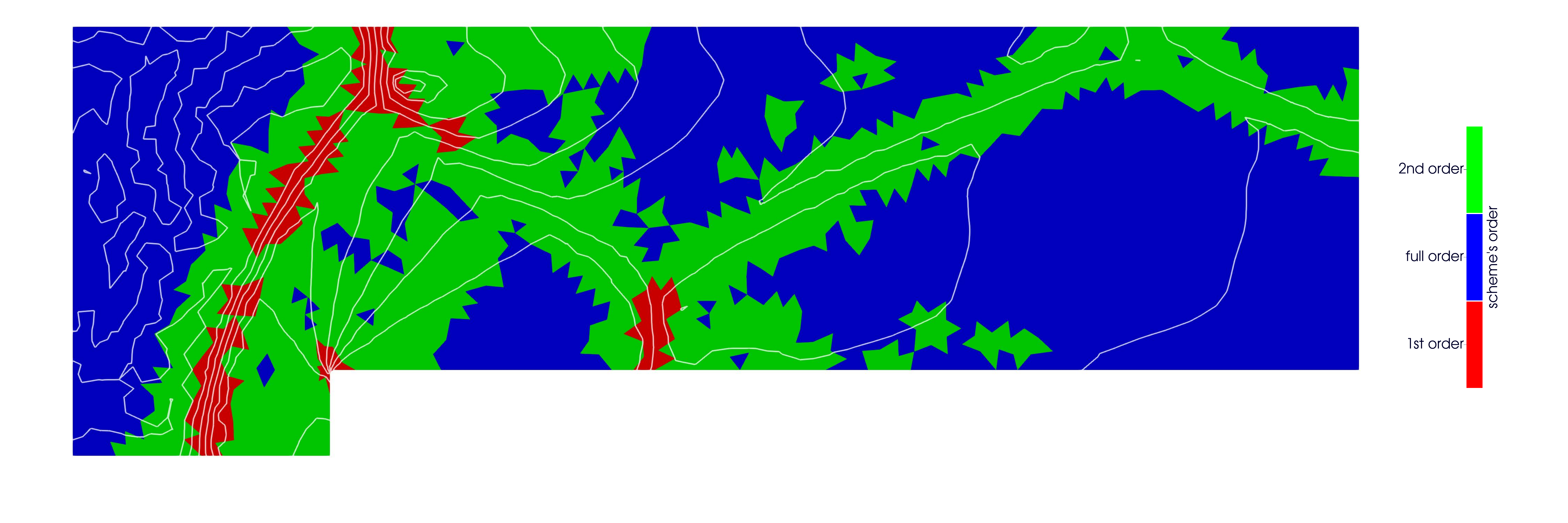}
        \includegraphics[width=\linewidth]{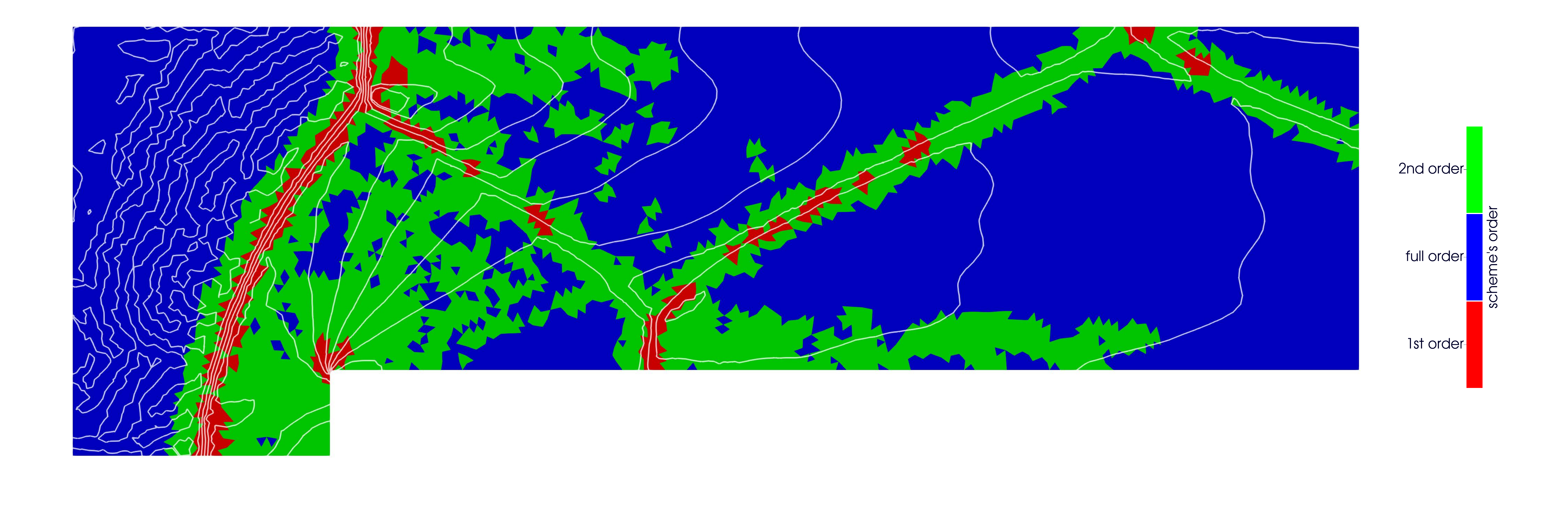}
        \includegraphics[width=\linewidth]{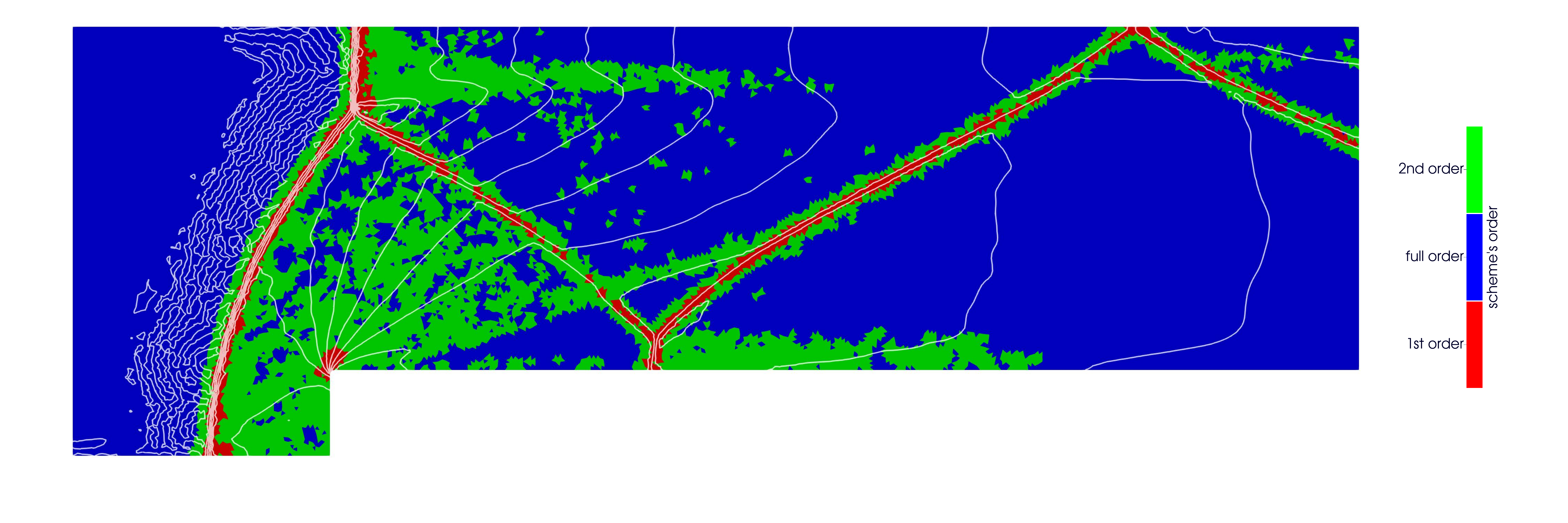}
        \includegraphics[width=\linewidth]{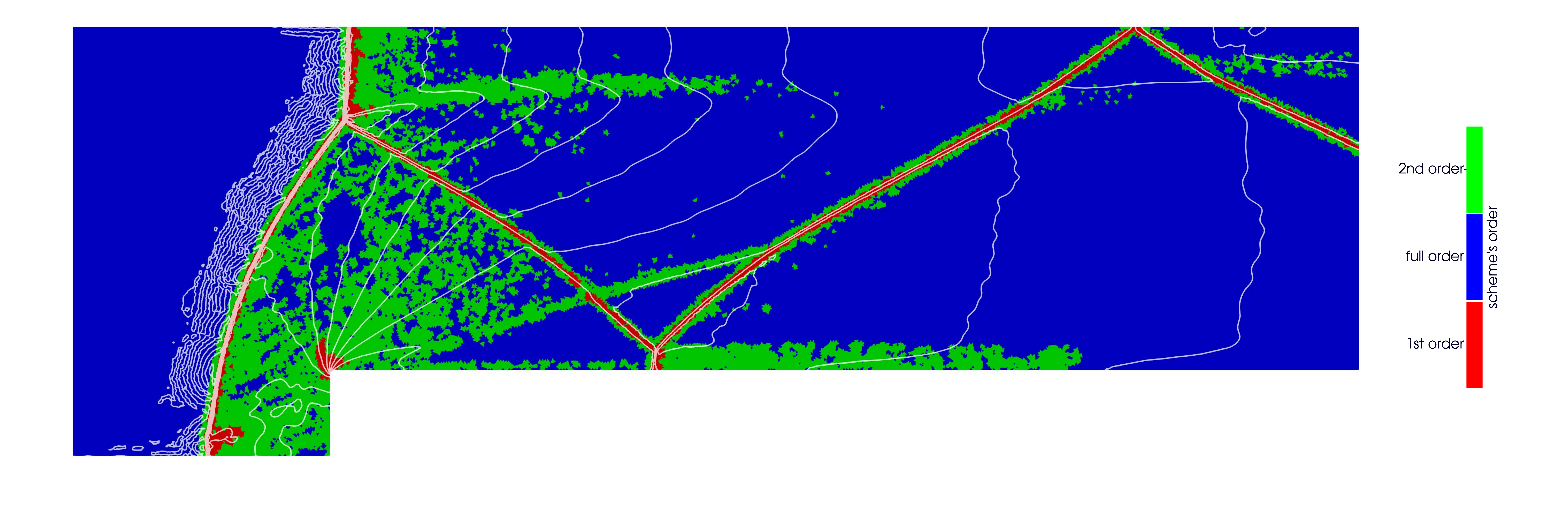}
        \vspace{-6ex}
        \caption{\label{fig:no_resolution_scaling}$\Delta s $}
    \end{subfigure}
    \begin{subfigure}{0.245\textwidth}
        \centering
        \includegraphics[width=\linewidth]{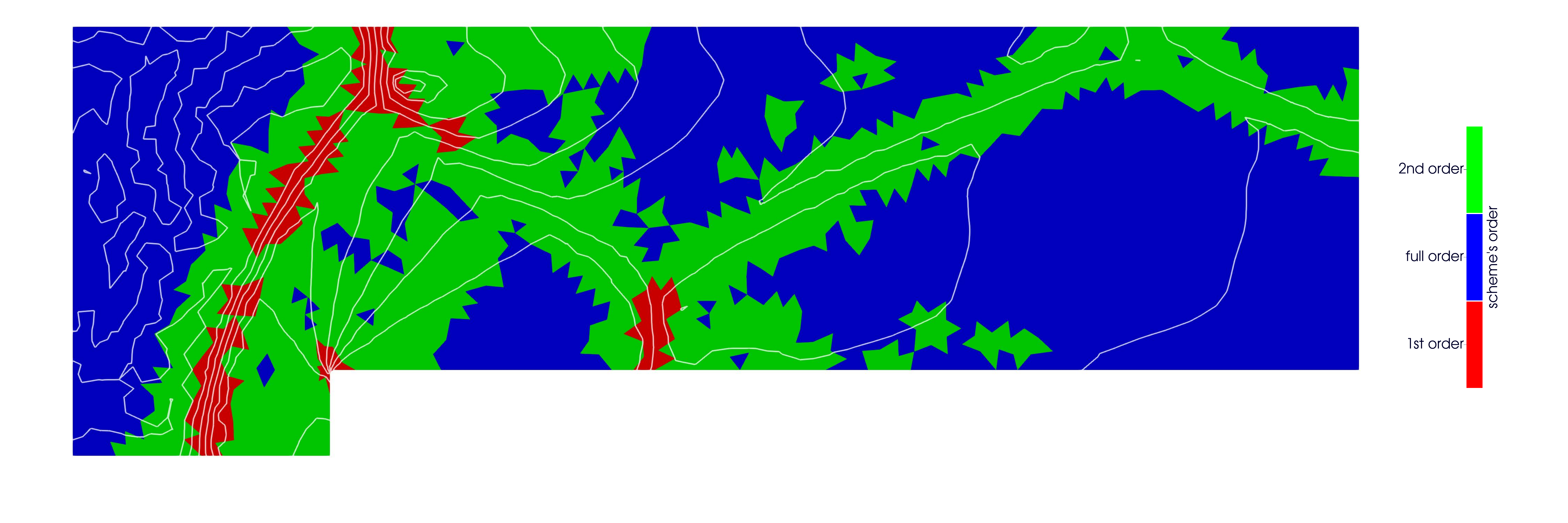}
        \includegraphics[width=\linewidth]{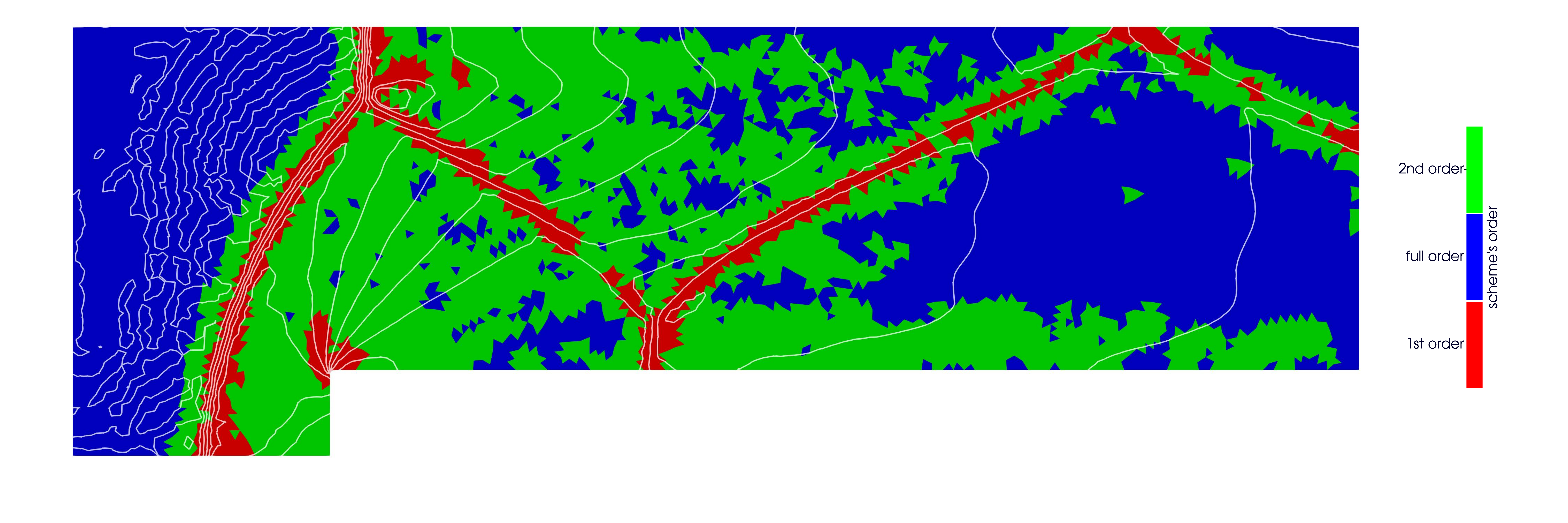}
        \includegraphics[width=\linewidth]{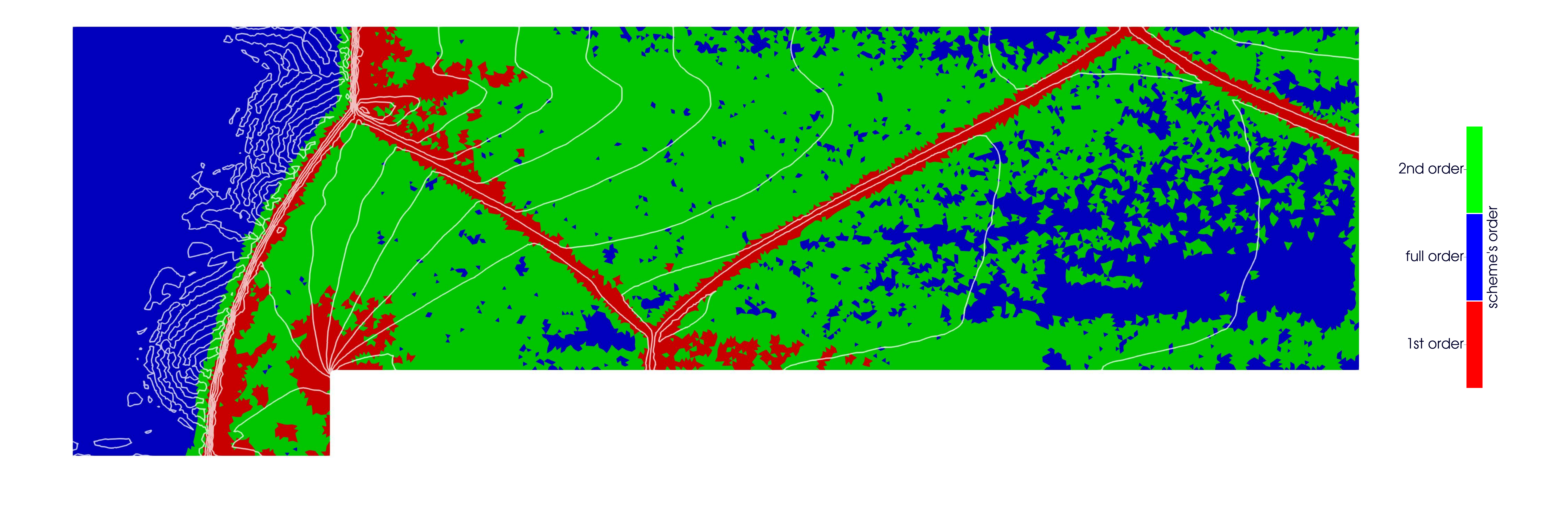}
        \includegraphics[width=\linewidth]{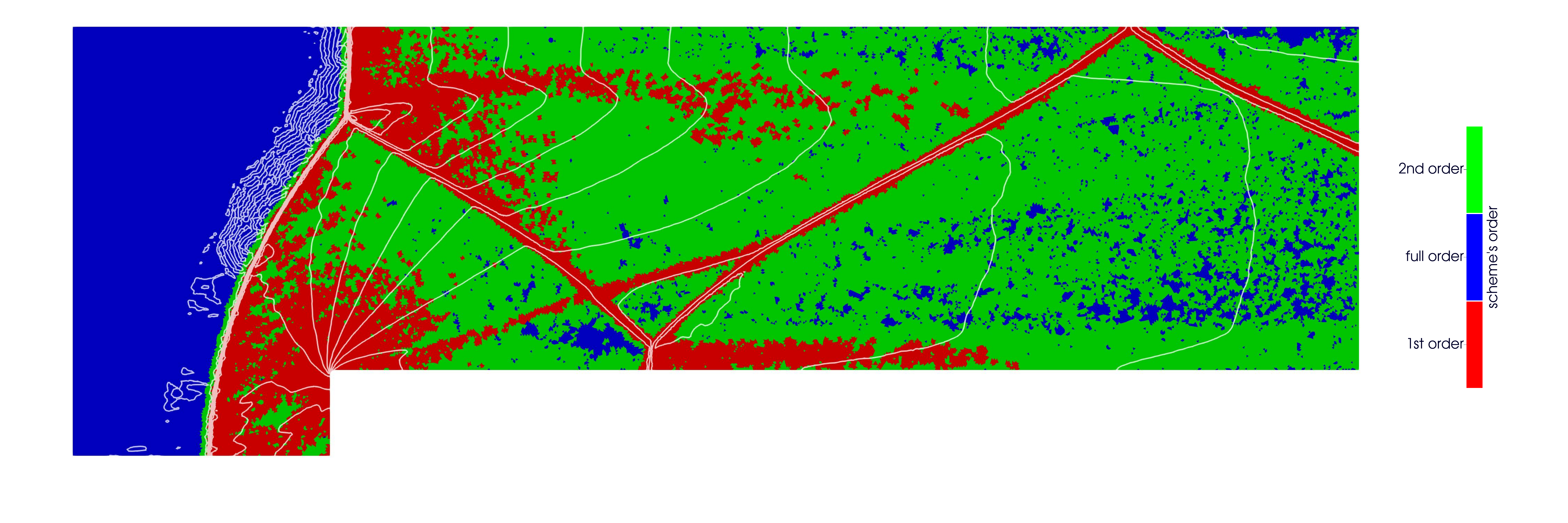}
        \vspace{-6ex}
        \caption{\label{fig:div_avg_cell_size_scaling}$\Delta s / \overline{h}$}
    \end{subfigure}
    \caption{\label{fig:resolution_scaling} Comparison of different strategies of scaling the entropy production trigger depending on the mesh resolution tested on Forward Facing Step problem on four mesh resolutions. Green designates the cells which have been flagged in the first pass, red designates cells which have been flagged in both passes and hence have used $1^{st}$ order reconstruction.}
\end{figure}

As seen in Figure \ref{fig:resolution_scaling} the scaling $\Delta s/ \overline{h}$ results in an increasingly larger fraction of the domain occupied by flagged cells, which is opposite to the desired effect. Not including any scaling results in the fraction of flagged cells remaining approximately constant, which could be useful in some scenarios, but is not desirable in presented test problems, in which refining the mesh should result in sharper resolution of the waves and interfaces therefore fewer cells should require flagging. This is achieved by the $\Delta s \cdot \overline{h}$ scaling (see Figure \ref{fig:times_avg_cell_size_scaling}). This is consistent with entropy production scaling with $ 1/dt $ reported in \cite{SharkEntropyPaper}, because $ dt \propto \overline{h} $ due to the CFL condition. $\Delta s \cdot \overline{A}$ scaling also reduces the fraction of flagged cells, but does so at too high a rate which causes a problem opposite to the first two triggers - it fails to flag the shock regions which should be recomputed. Hence $\Delta s \cdot \overline{h}$ scaling is recommanded and it is used in all subsequent tests which involve mesh refinement.

\section{Applications}\label{Applications}

The methods were tested in a large number of progressively more complex test problems. 
Firstly, a 1D flows run with a well-known solution for validation and measuring the accuracy,
then 2D tests with complex flow patterns which can be found in industrial applications, to test robustness, and finally 2D multi-component flows to validate compatibility with more complex sets of equations. 



\subsection{Shu-Osher Test}\label{ShuOsher}

The 1D test problem is the well-established benchmark proposed by Shu-Osher \citep{Shu198932} which involves the interaction of a shock wave with a smooth entropy wave. The domain is $[-4.5,4.5]$ and the initial profile consists of a shock wave at $x=-4$ with post-shock conditions $\left(\rho,u,p\right)=\left(3.857143,\ 2.629369,\ 10.333333\right)$ travelling right over an entropy wave defined in the rest of the domain as $\left(\rho,u,p\right)=\left(1+0.2sin(5x),0,1\right)$.
For compatibility with the presented code the test problem is extended to 2D by setting the domain to $[-4.5, 4.5]\times[0, 1]$, the \textit{y}-velocity to zero, and applying periodic boundary conditions along the long sides, so that the solution should have no variation in the \textit{y}-direction.
Supersonic inflow and outflow are applied on the left and right ends of the domain respectively.
Progressively finer 2D unstructured meshes with resolutions $h = 1/12.5,\ 1/25,\ 1/50$ and approximately \textit{\replaced{1800, 7300, 28200}{1815, 7287, 28195}} triangular and quadrilateral elements are utilised. The calculation is run until $t=1.8$. 

\begin{figure}[H]
  \vspace{-2ex}
  \centering
  \begin{subfigure}[t]{0.195\textwidth}
    \centering
    \includegraphics[width=\textwidth]{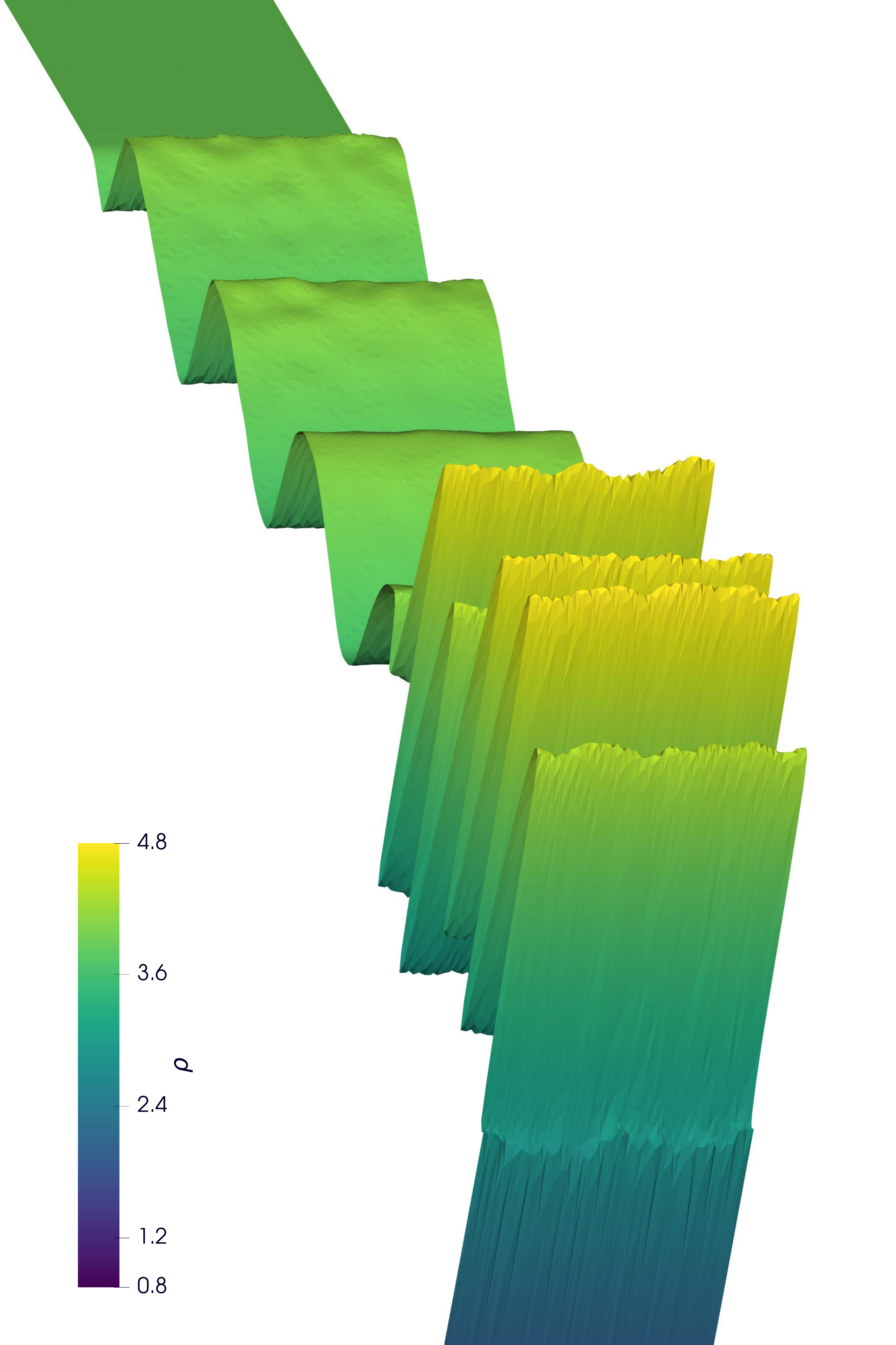}
    \caption{\label{fig:ShuOsher_CWENO} \centering CWENO5}
  \end{subfigure}
  \begin{subfigure}[t]{0.195\textwidth}
      \centering
      \includegraphics[width=\textwidth]{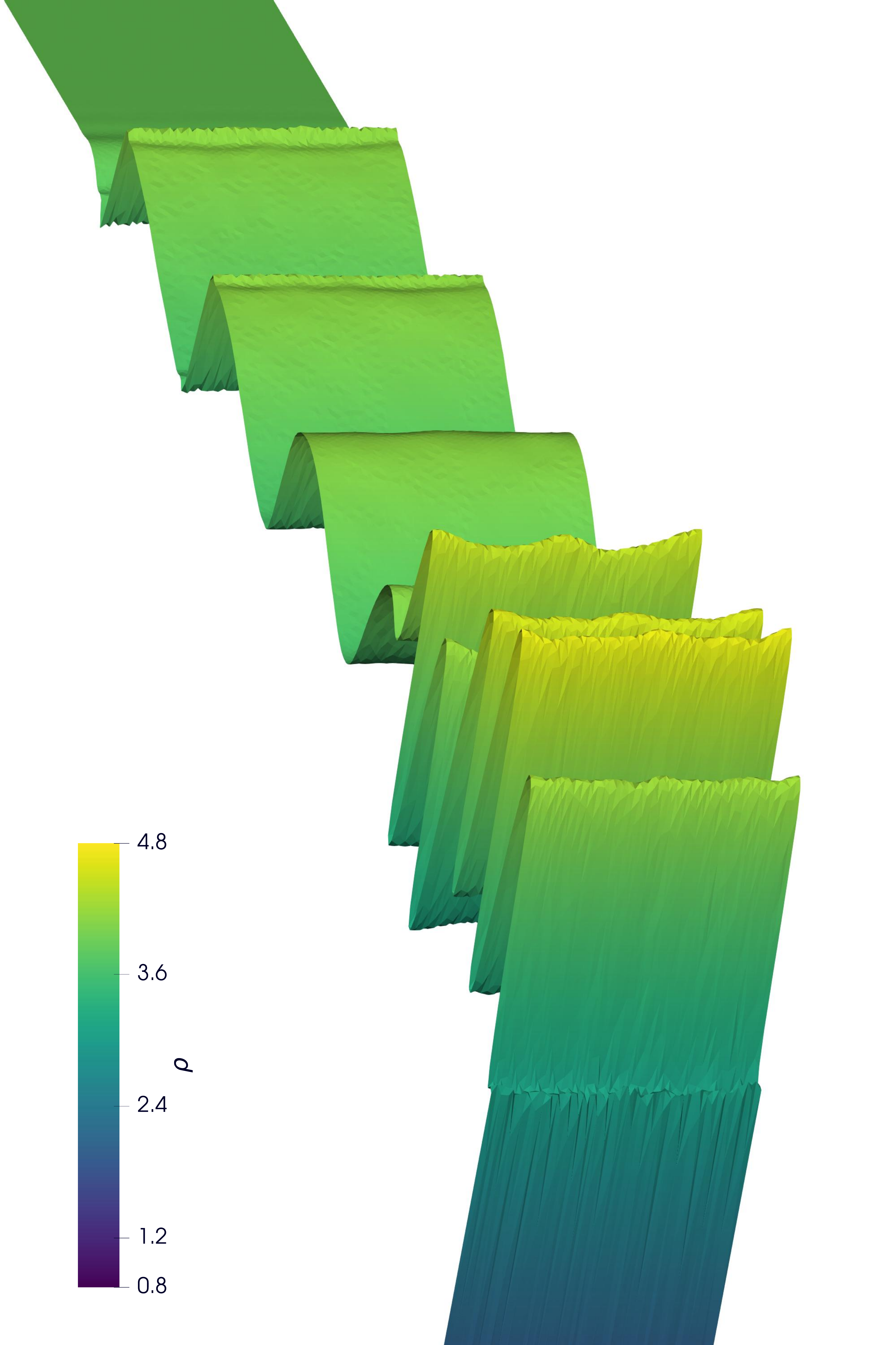}
      \caption{\label{fig:ShuOsher_ClassicMOOD} \centering Original MOOD}
  \end{subfigure}
  \begin{subfigure}[t]{0.195\textwidth}
      \centering
      \includegraphics[width=\textwidth]{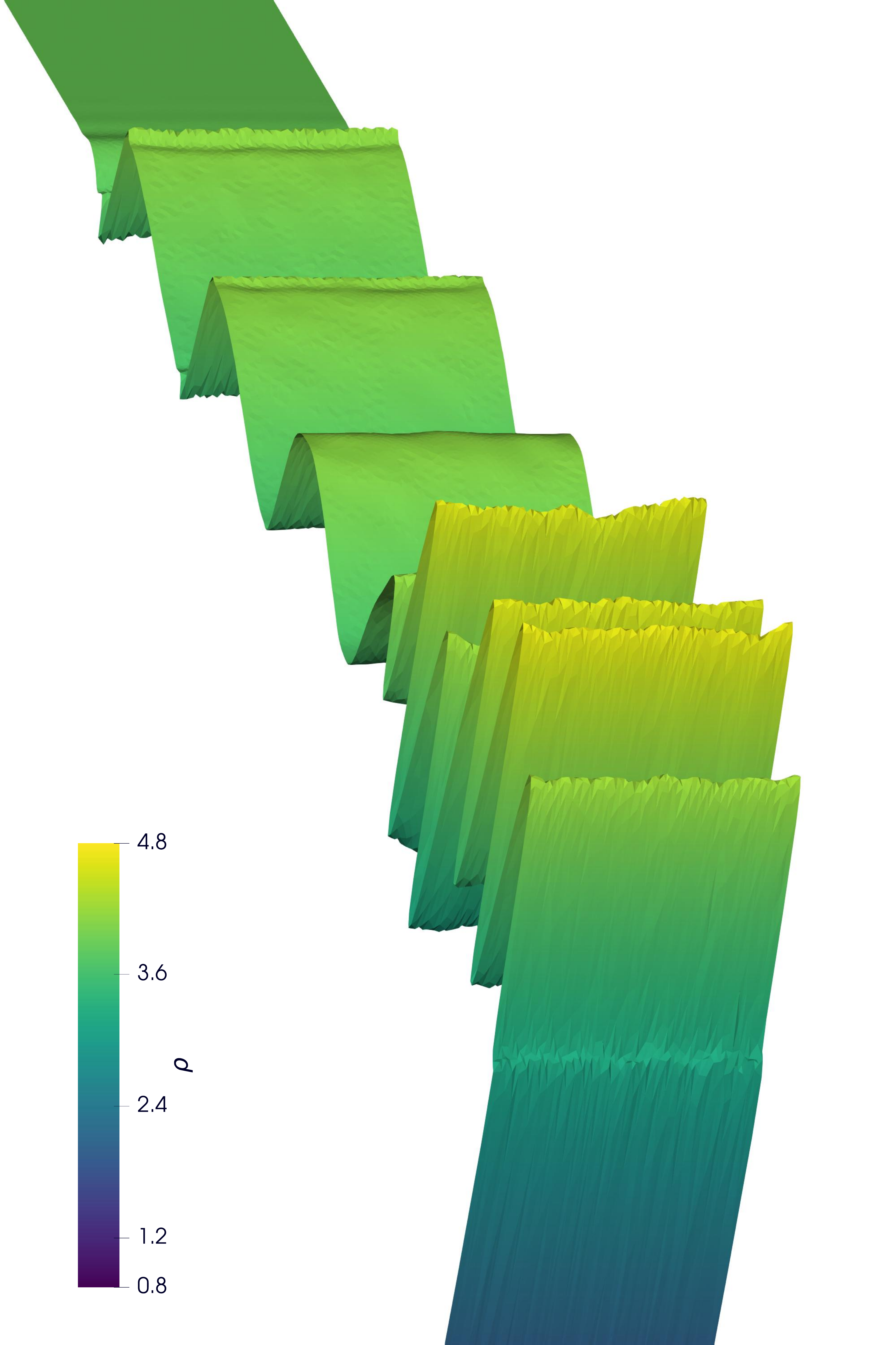}
      \caption{\label{fig:ShuOsher_EntropyMOOD} \centering entropy production MOOD}
  \end{subfigure}
  \begin{subfigure}[t]{0.195\textwidth}
      \centering
      \includegraphics[width=\textwidth]{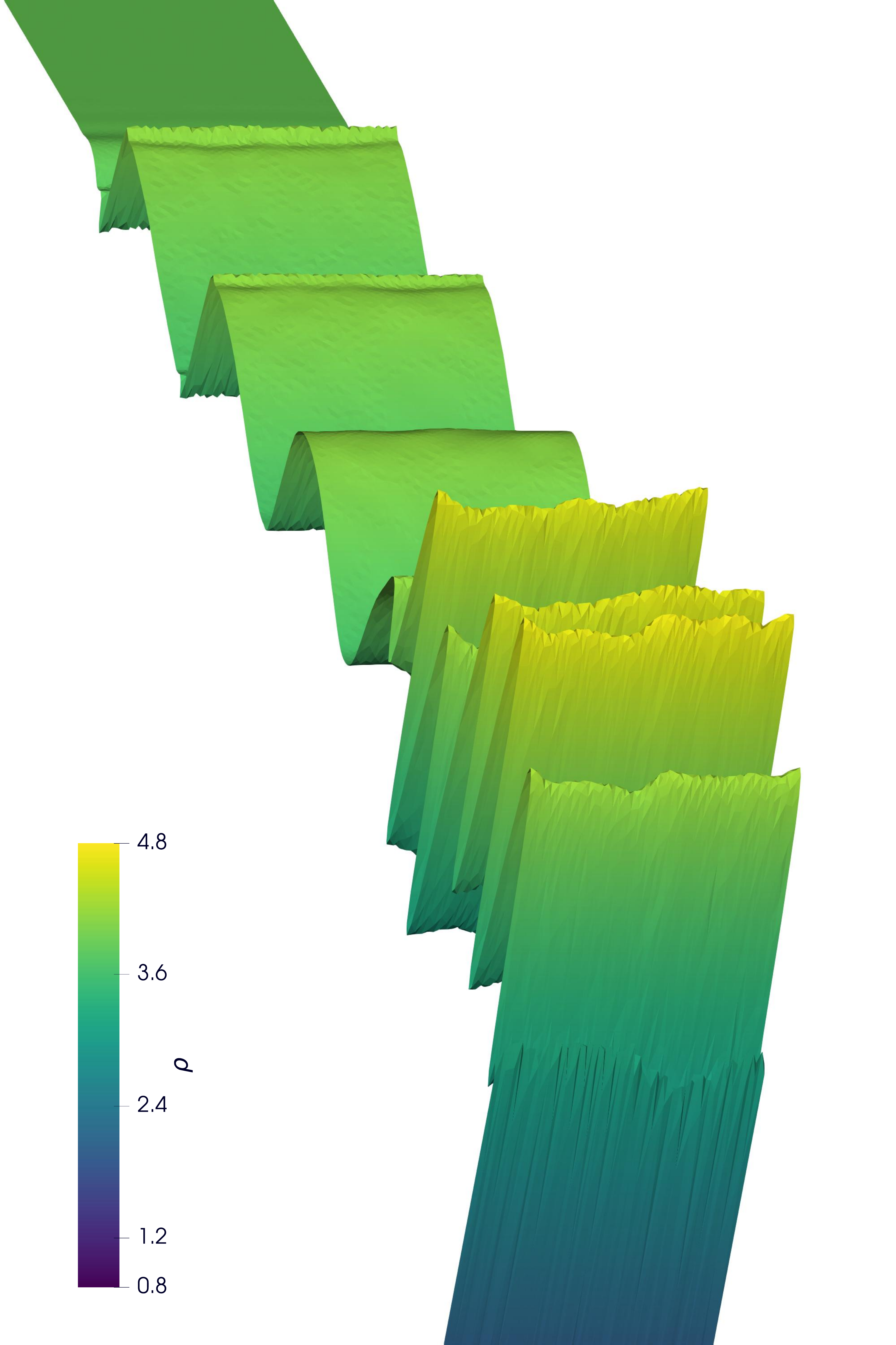}
      \caption{\label{fig:ShuOsher_CWENO-MOOD_classic} \centering MOOD-CWENO hybrid with original MOOD trigger}
  \end{subfigure}
  \begin{subfigure}[t]{0.195\textwidth}
      \centering
      \includegraphics[width=\textwidth]{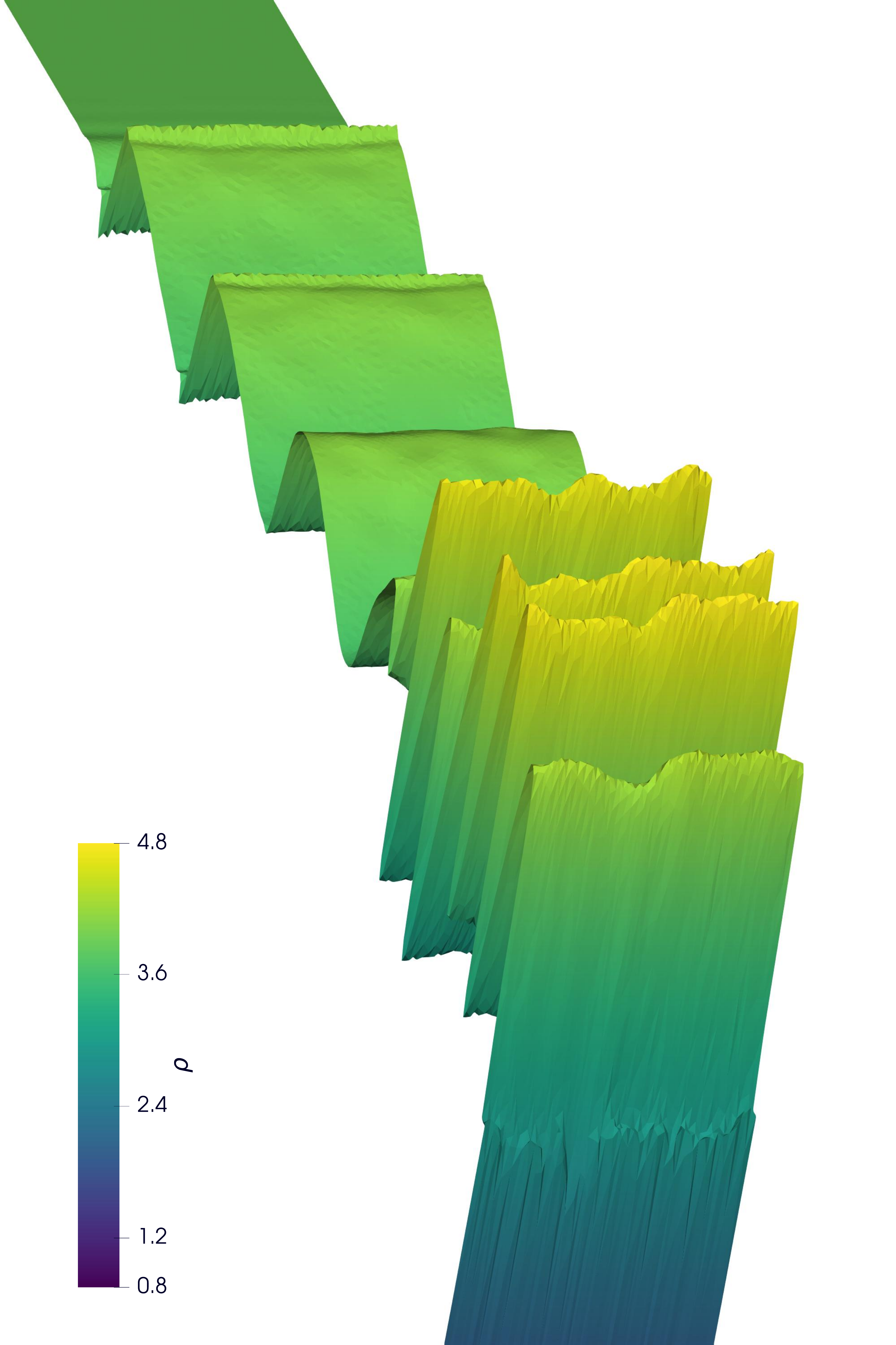}
      \caption{\label{fig:ShuOsher_CWENO-MOOD_entropy} \centering MOOD-CWENO hybrid with entropy production trigger}
  \end{subfigure}
  \vspace{-4ex}
\end{figure}

\begin{figure}[H]
  \ContinuedFloat
  \centering
  \includegraphics[width=0.667\textwidth]{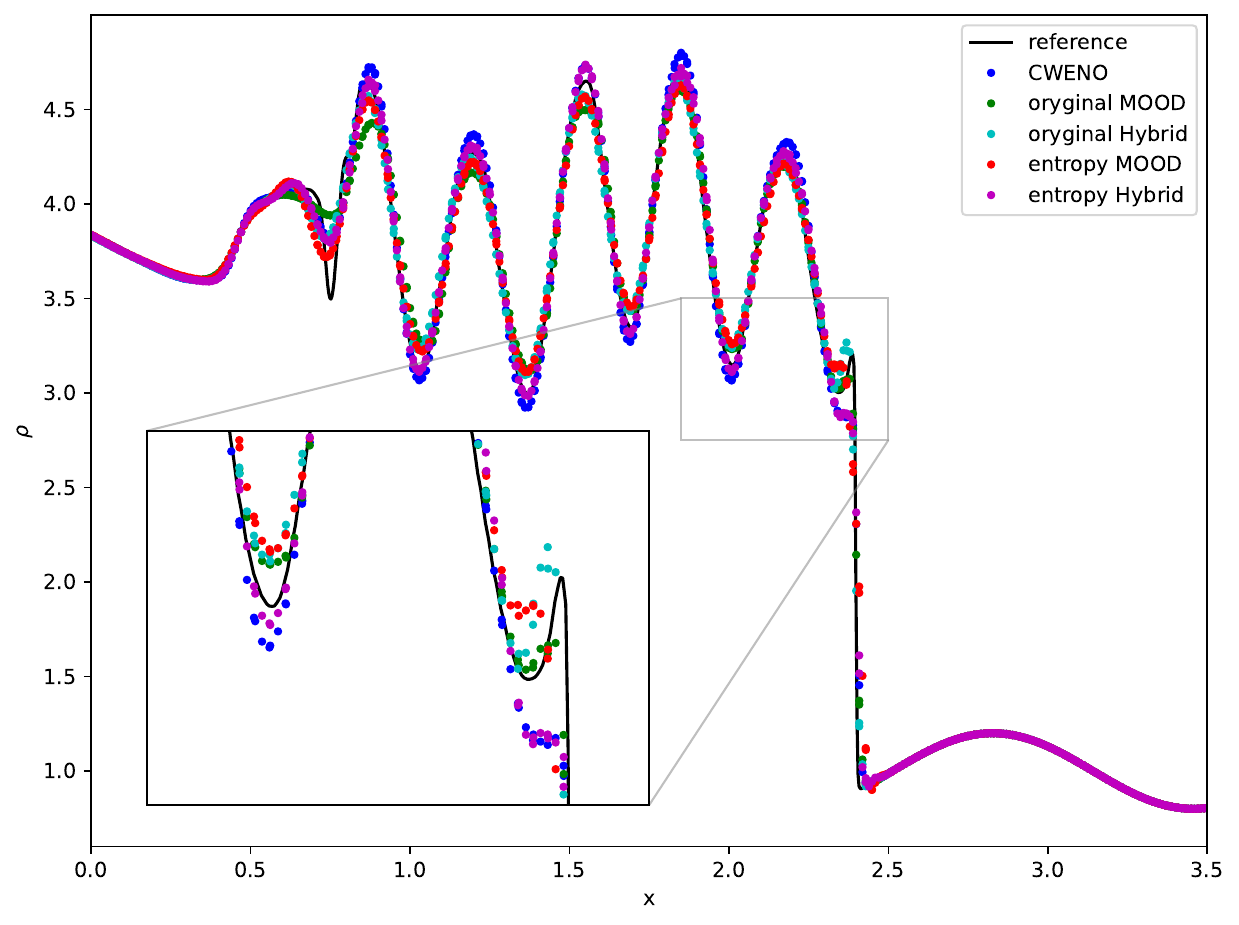}
  \vspace{-2ex}
  \caption{\label{fig:ShuOsher_densityPlots} Final density profiles of Shu-Osher test case computed using five compared methods on the finest mesh ($h = 1/50$). Top: full two-dimensional solution. Bottom: cross-section along the line in the centre of the domain. Close-up on the key interaction region behind the shock. Reference solution computed using very high resolution one-dimensional solver.}
\end{figure}

As seen in Figure \ref{fig:ShuOsher_densityPlots}\textcolor{blue}{a-e} all five methods produce very similar solutions, but the cross section reveals some differences in the heights of the shocked waves can be observed. The MOOD-CWENO hybrid with entropy trigger (\ref{fig:ShuOsher_CWENO-MOOD_entropy}) is the closest to the reference CWENO5 implementation (\ref{fig:ShuOsher_CWENO}).

\subsection{Double Mach Reflection}\label{DoubleMachReflection}

The popular double Mach reflection test case of Woodward and Colella \citep{Woodward1984115} is employed to assess the non-oscillatory properties of the hybrid schemes and their accuracy in resolving strong shocks. The problem consists of a shock wave moving diagonally with Mach number $M=10$ and hitting a ramp inclined by $\alpha=60^{\circ}$, reflecting and gradually generating secondary shock waves and triple points. Ahead of the shock the fluid is at rest with uniform density and pressure $(\rho,\ p)=(1.4,\ 1.0)$, and the post-shock conditions are assumed behind. Reflecting boundary conditions are applied at the bottom of the domain. At the top, exact solution of an isolated moving oblique shock wave with $M=10$ is prescribed. At the left and right side of the domain, supersonic inflow and outflow boundary conditions are applied, respectively.
The computational domain is given by $[0,\ 4]\times[0,\ 1]$ and the initial shock position is $x=1/6$. For additional details regarding the setup of the problem the reader is referred to the original work of Woodward and Colella \citep{Woodward1984115}. An unstructured mesh of approximately $2910$, $11540$, $46204$, $189029$ triangular and quadrilateral cells was used. That corresponds to the edge sizes of $h=1/20$ to $h=1/160$ at the boundaries. 

\begin{figure}[H]
  \vspace{-2ex}
  \centering
  \begin{subfigure}{0.75\textwidth}
    \centering
    \includegraphics[width=\textwidth, trim = {0 8cm 0 8cm}, clip]{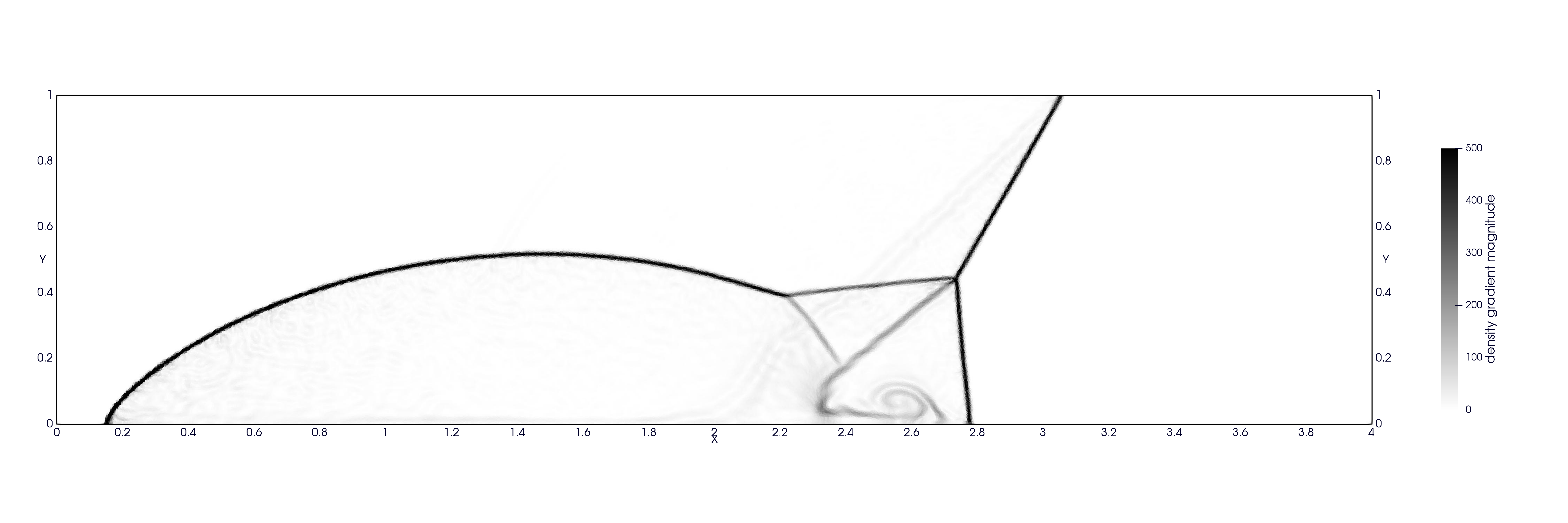}
    \caption{\label{fig:DMR_CWENO-superfine} reference solution - CWENO5 on finer mesh (189029 cells)}
  \end{subfigure}
  \begin{subfigure}{0.75\textwidth}
    \centering
    \includegraphics[width=\textwidth, trim = {0 8cm 0 6cm}, clip]{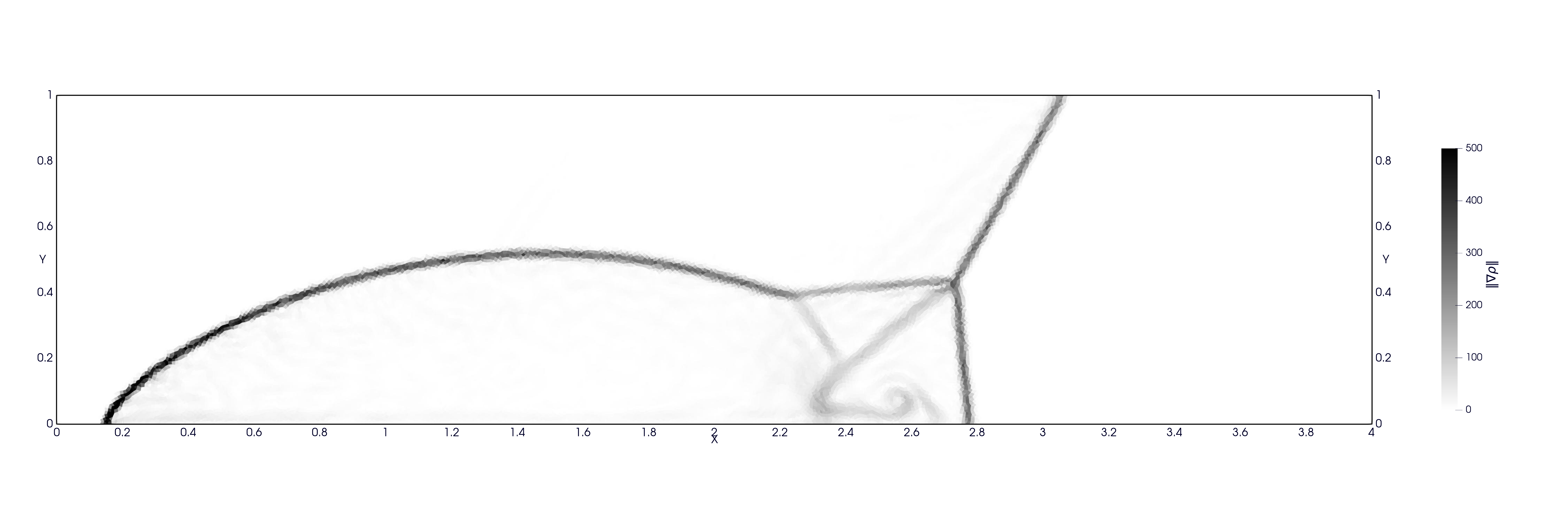}
    \caption{\label{fig:DMR_CWENO-fine} CWENO5}
  \end{subfigure}
  \begin{subfigure}{0.75\textwidth}
      \centering
      \includegraphics[width=\textwidth, trim = {0 2cm 0 4cm}, clip]{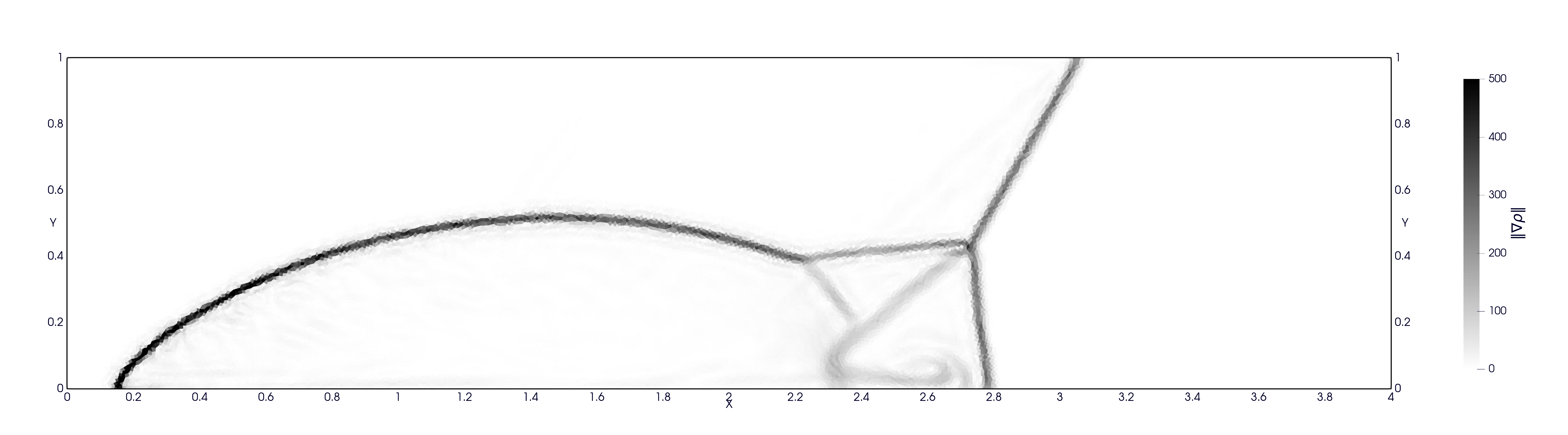}
      \caption{\label{fig:DMR_ClassicMOOD} Original MOOD}
  \end{subfigure}
  \begin{subfigure}{0.75\textwidth}
      \centering
      \includegraphics[width=\textwidth, trim = {0 2cm 0 4cm}, clip]{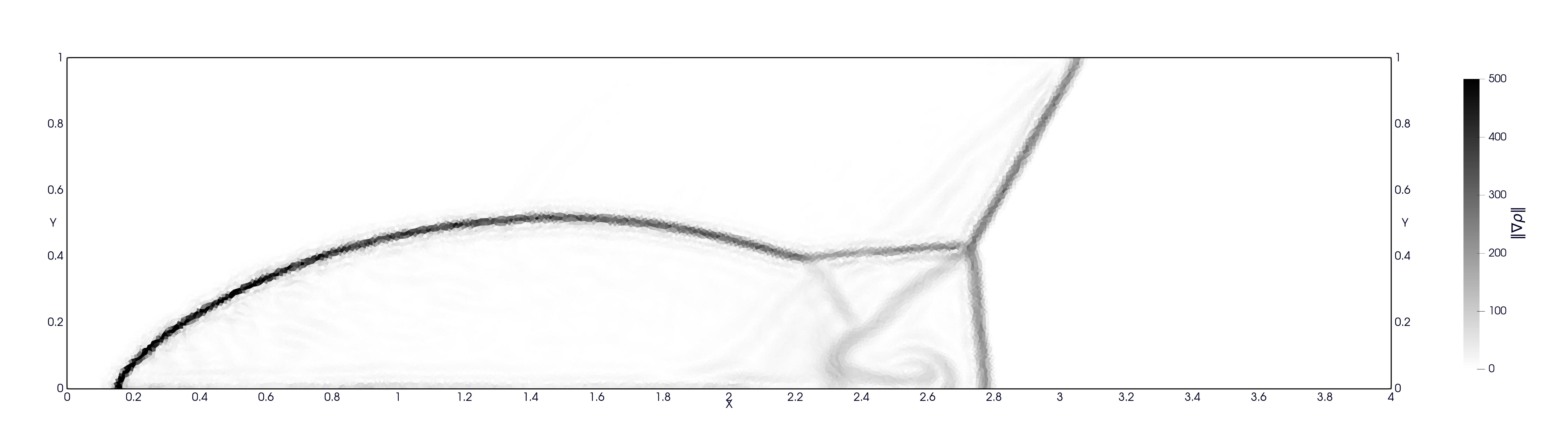}
      \caption{\label{fig:DMR_EntropyMOOD} Entropy MOOD}
  \end{subfigure}
  \begin{subfigure}{0.75\textwidth}
      \centering
      \includegraphics[width=\textwidth, trim = {0 2cm 0 4cm}, clip]{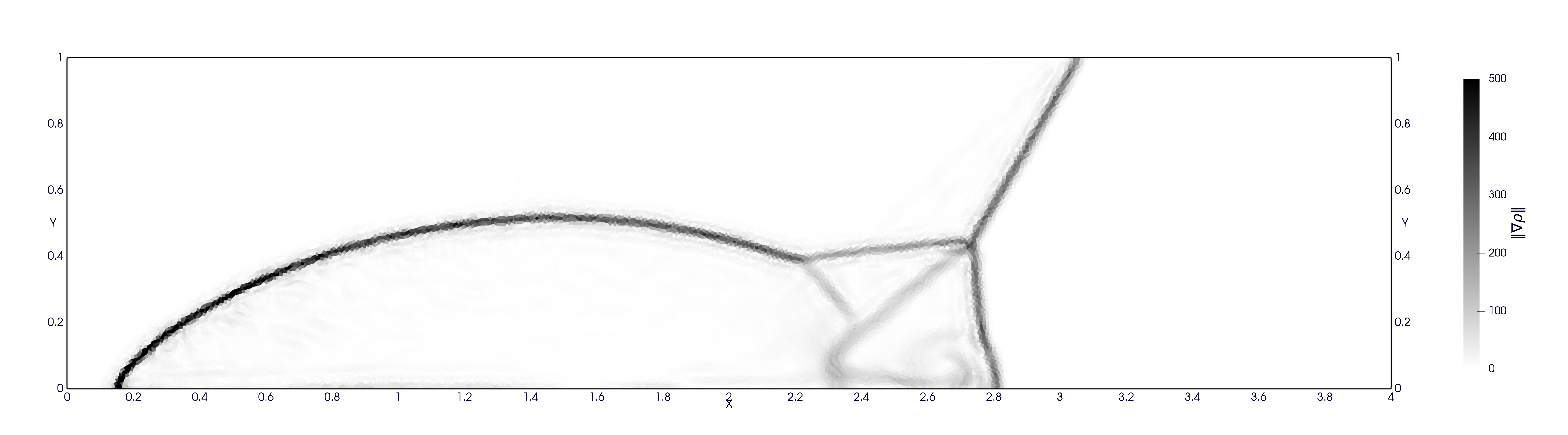}
      \caption{\label{fig:DMR_CWENO-MOOD_classic} MOOD-CWENO hybrid with original MOOD trigger}
  \end{subfigure}
  \begin{subfigure}{0.75\textwidth}
      \centering
      \includegraphics[width=\textwidth, trim = {0 2cm 0 4cm}, clip]{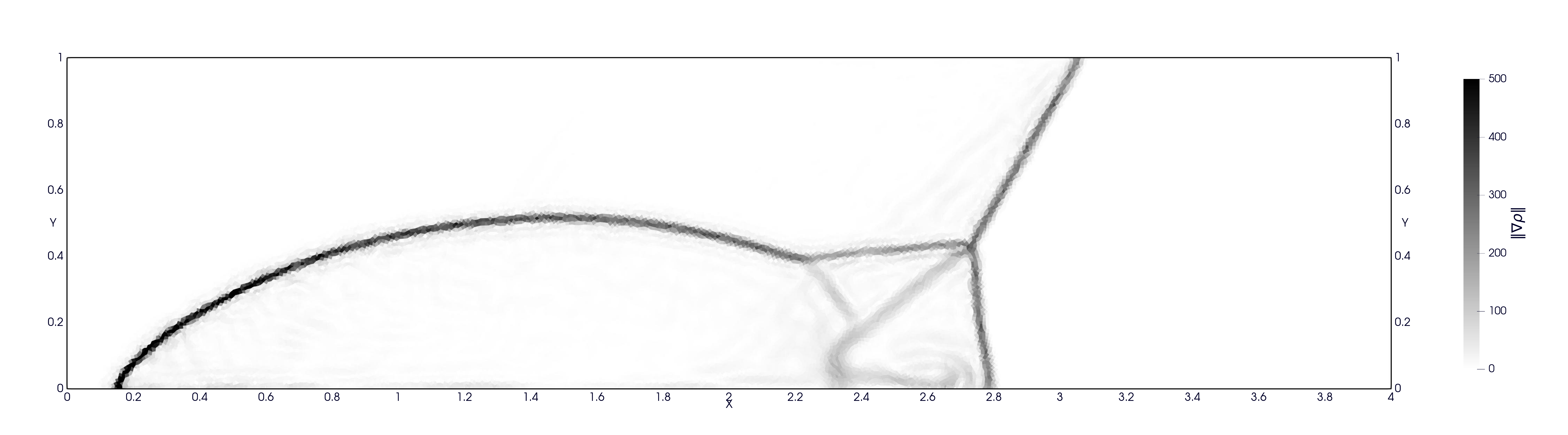}
      \caption{\label{fig:DMR_CWENO-MOOD_entropy} MOOD-CWENO hybrid with entropy trigger}
  \end{subfigure}
  \caption{\label{fig:DMR_DensityPlots} Computed magnitude of density gradient for the Double Mach Reflection test case at the final time $t=0.2$ with CWENO5 scheme, MOOD and MOOD-CWENO hybrid methods with two decision algorithms, all using HLLC numerical flux on an unstructured mesh with 46204 cells.}
\end{figure}

As seen in Figure \ref{fig:DMR_DensityPlots} the original MOOD (\ref{fig:DMR_ClassicMOOD}) and the MOOD-CWENO hybrid with entropy trigger (\ref{fig:DMR_CWENO-MOOD_entropy}) produce solutions that match very well with the reference (\ref{fig:DMR_CWENO-fine}). The solution using MOOD with the entropy based trigger (\ref{fig:DMR_EntropyMOOD}) also matches, but has slightly more diffused front shock. The MOOD-CWENO hybrid with the original trigger produces incorrect position and shape of the front shock at the point of contact with the wall.

\begin{figure}[H]
  \centering
  \begin{subfigure}{0.475\textwidth}
    \centering
    \includegraphics[width=\textwidth, trim = {0 2cm 0 2cm}, clip]{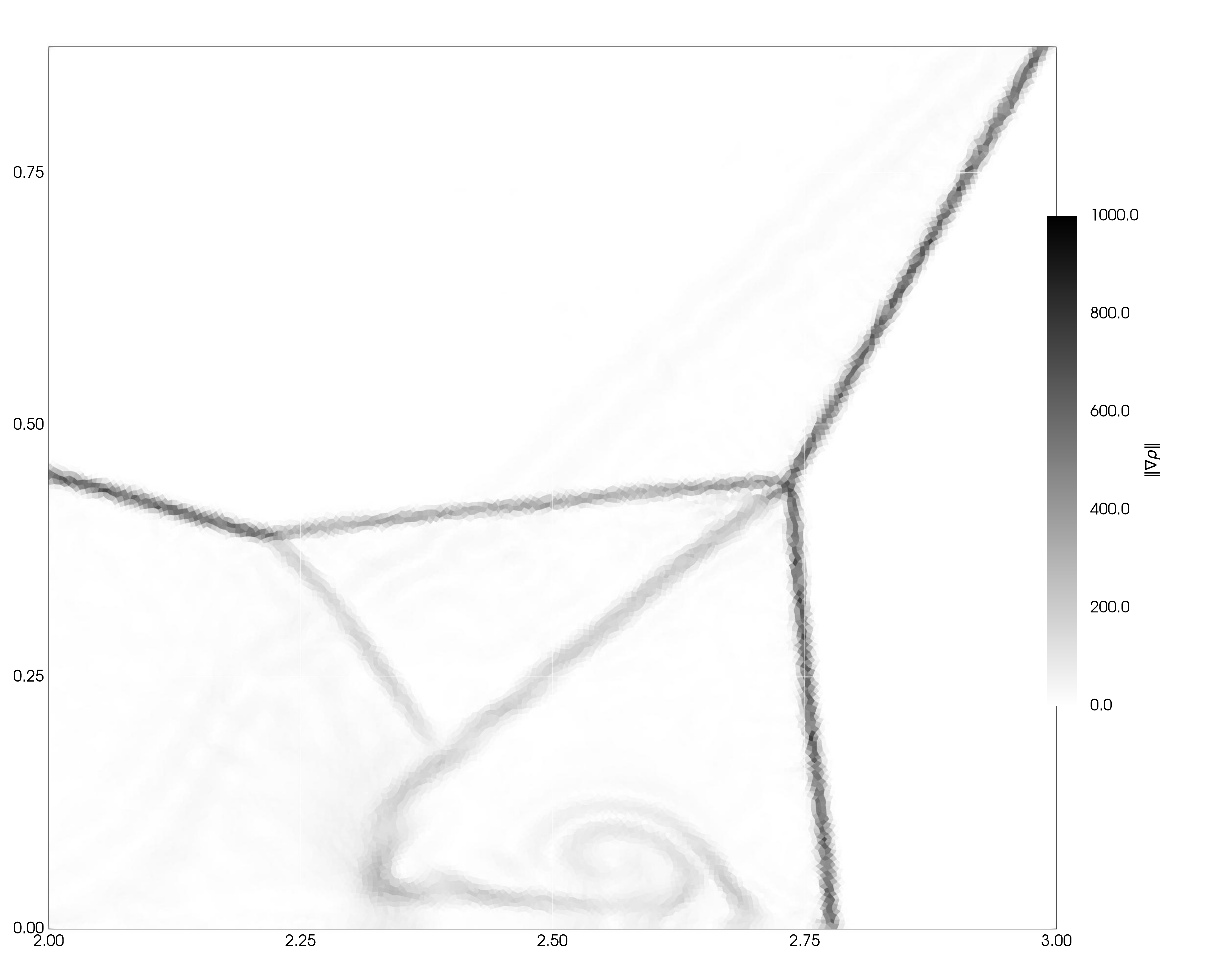}
    \caption{\label{fig:DMR_CWENO-superfine_closeup} CWENO5}
  \end{subfigure}
  \begin{subfigure}{0.475\textwidth}
    \centering
    \includegraphics[width=\textwidth, trim = {0 2cm 0 2cm}, clip]{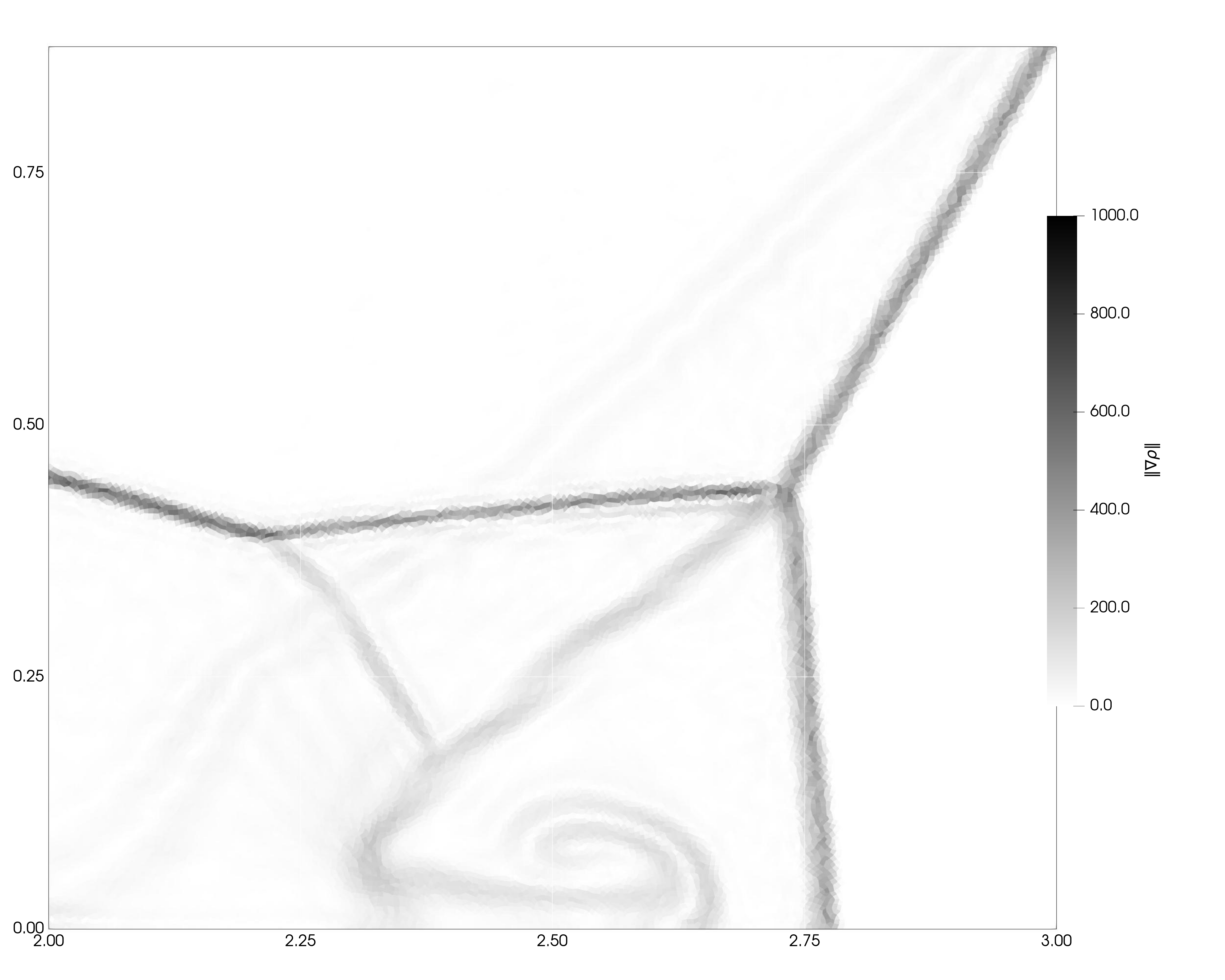}
    \caption{\label{fig:DMR_EntropyMOOD-superfine_closeup} MOOD with entropy production trigger}
  \end{subfigure}
  \caption{\label{fig:DMR_DensityPlotsCloseup} Computed magnitude of density gradient for the double Mach reflection test case at final time $t=0.2$ with CWENO5 scheme and MOOD scheme with the proposed entropy production trigger an unstructured mesh with 189029 cells.}
\end{figure}

In the Figure \ref{fig:DMR_DensityPlotsCloseup} one can confirm the close agreement between the CWENO reference and the proposed Entropy-MOOD method. MOOD produces minimally more diffused initial shock, but the secondary waves behind have a slightly higher gradient magnitude, meaning the jump is sharper.

\subsection{Forward Facing Step}\label{ForwardFacingStep}

The second test case of Woodward and Colella \citep{Woodward1984115} is the Forward Facing Step. It models a Mach 3 wind tunnel in which the airflow hits a sharp step at the bottom of the tunnel creating a shock wave that repeatedly reflects from the walls and itself creating a complex shape. The computational domain is $[0, 3]\times[0, 1]$. The step has height 0.2 and is positioned at 0.6. The gas in the tunnel is modelled as an ideal gas with $\gamma = 1.4$, initial velocity $\vec{u} = (3.0,\ 0.0)$, density $\rho = 1.0$, and pressure $p = \rho/\gamma$.
For contrast with the previous test, the simulations were run on purely triangular meshes in 4 resolutions with 3366, 12310, 50760, and 164598 cells.

\begin{figure}[H]
  \vspace{-2ex}
  \centering
  \begin{subfigure}{0.495\textwidth}
    \centering
    \includegraphics[width=\textwidth,trim={4cm 0 4cm 0},clip]{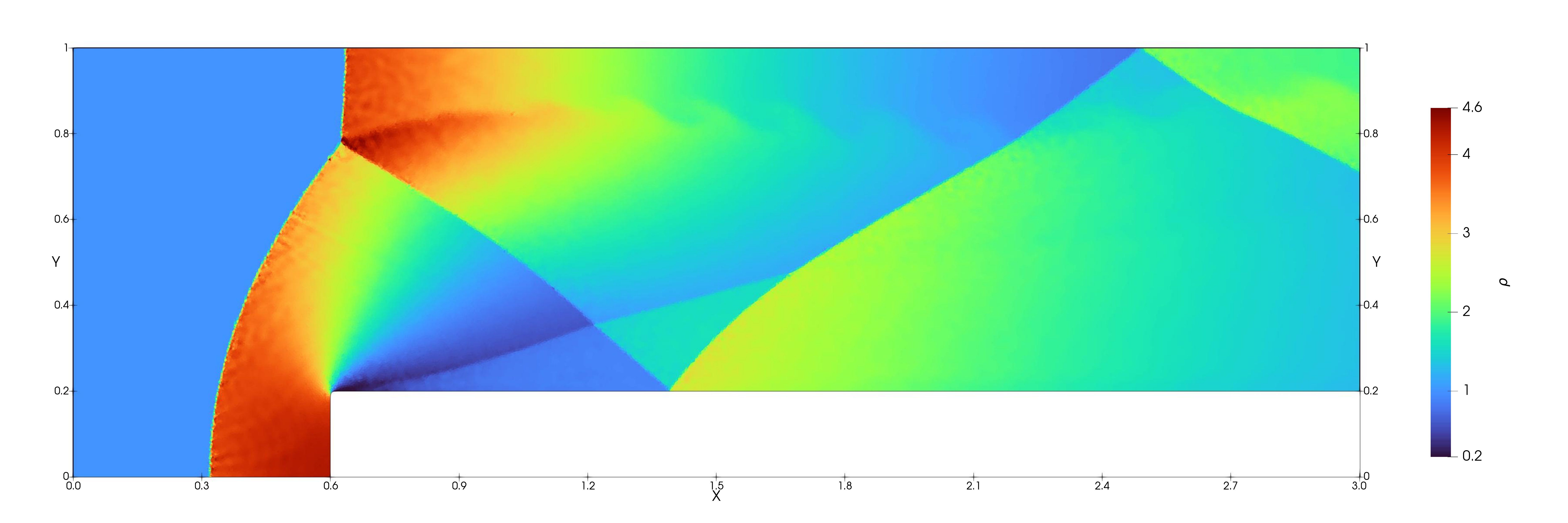}
    \vspace{-4ex}
    \caption{\label{fig:FFS_CWENO5-superfine} reference solution - CWENO5 on super fine mesh (164598 cells)}
  \end{subfigure}
  \begin{subfigure}{0.495\textwidth}
    \centering
    \includegraphics[width=\textwidth,trim={4cm 0 4cm 0},clip]{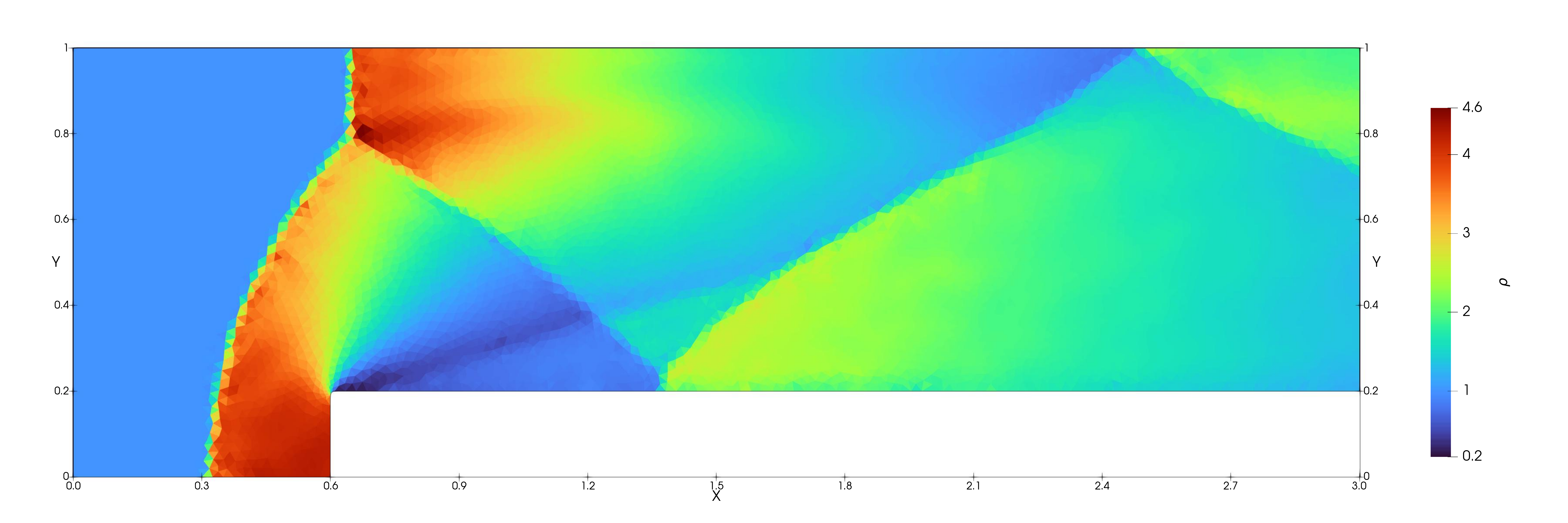}
    \vspace{-4ex}
    \caption{\label{fig:FFS_CWENO5-medium} CWENO5 on medium mesh}
  \end{subfigure}
  \begin{subfigure}{0.495\textwidth}
      \centering
      \includegraphics[width=\textwidth,trim={4cm 0 4cm 0},clip]{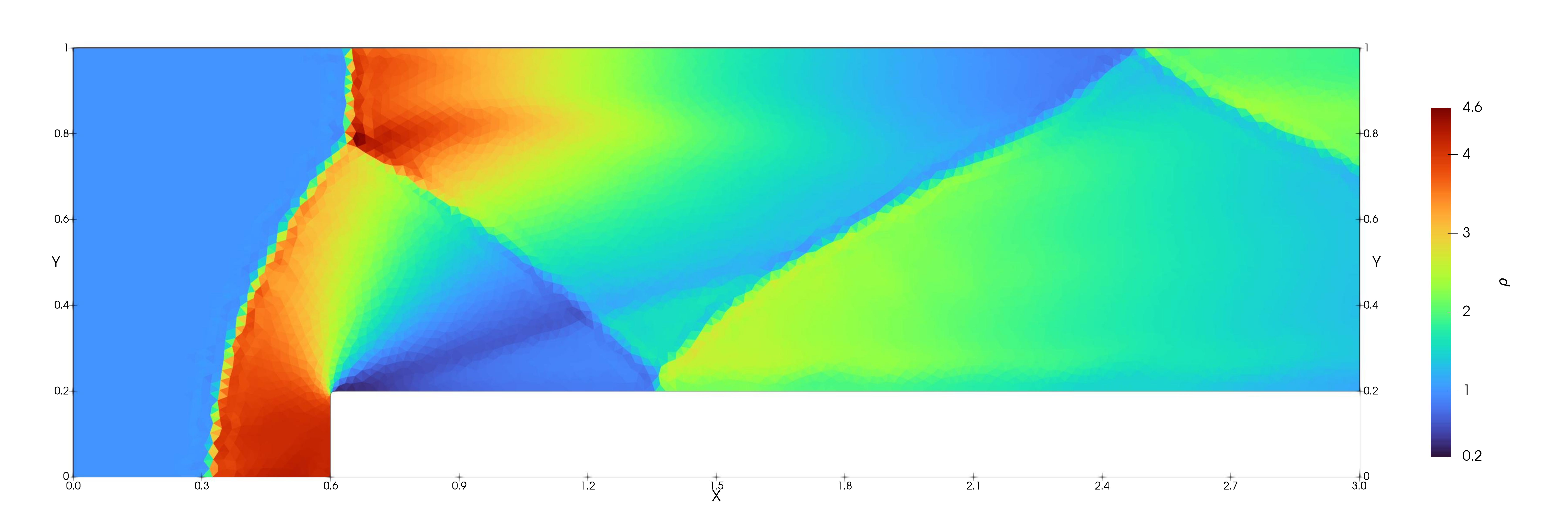}
      \vspace{-4ex}
      \caption{\label{fig:FFS_ClassicMOOD} Original MOOD}
  \end{subfigure}
  \begin{subfigure}{0.495\textwidth}
      \centering
      \includegraphics[width=\textwidth,trim={4cm 0 4cm 0},clip]{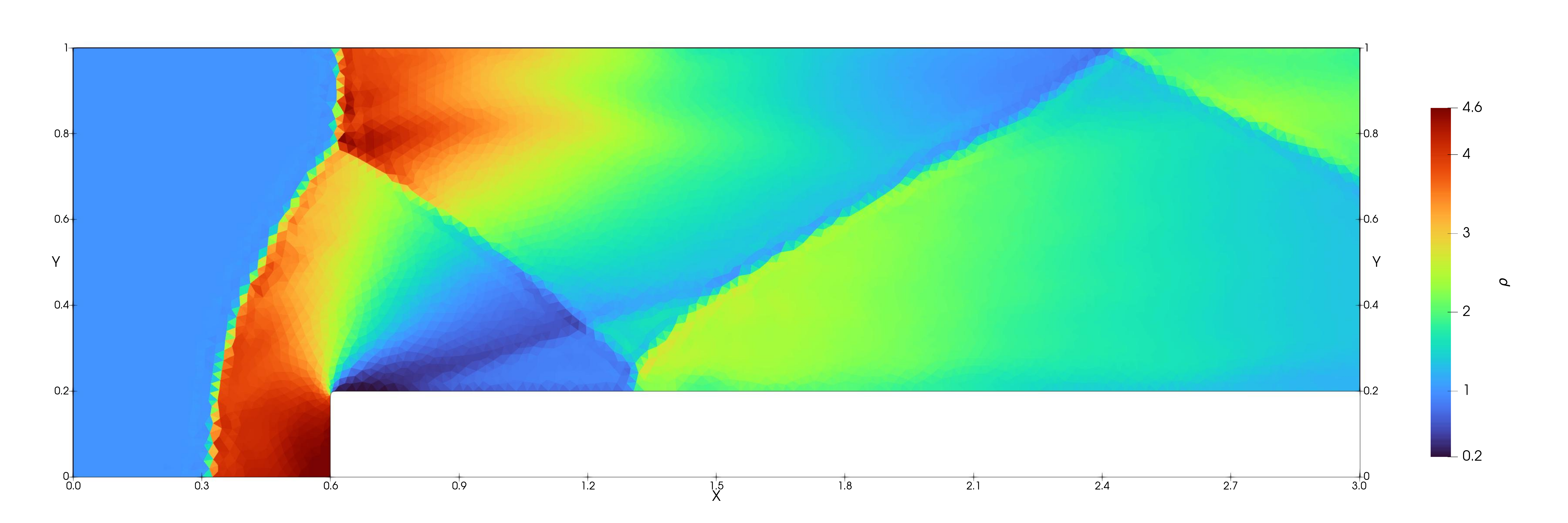}
      \vspace{-4ex}
      \caption{\label{fig:FFS_EntropyMOOD} Entropy MOOD}
  \end{subfigure}
  \begin{subfigure}{0.495\textwidth}
      \centering
      \includegraphics[width=\textwidth,trim={4cm 0 4cm 0},clip]{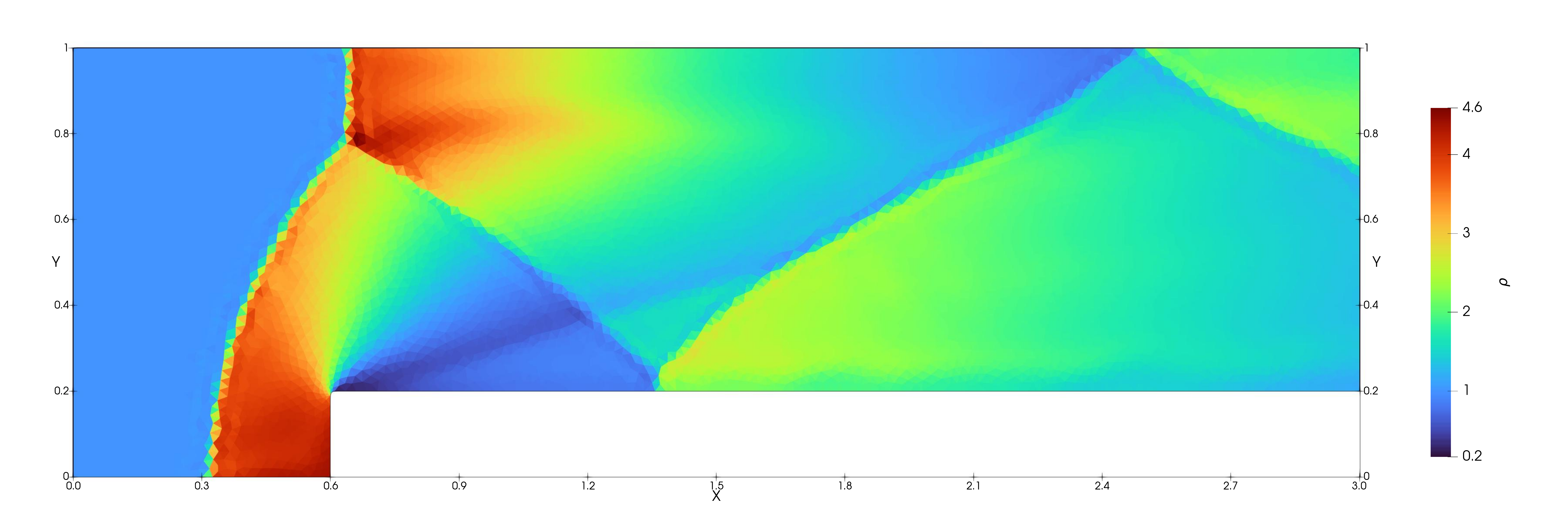}
      \vspace{-4ex}
      \caption{\label{fig:FFS_CWENO-MOOD_classic} CWENO-MOOD hybrid with original MOOD trigger}
  \end{subfigure}
  \begin{subfigure}{0.495\textwidth}
      \centering
      \includegraphics[width=\textwidth,trim={4cm 0 4cm 0},clip]{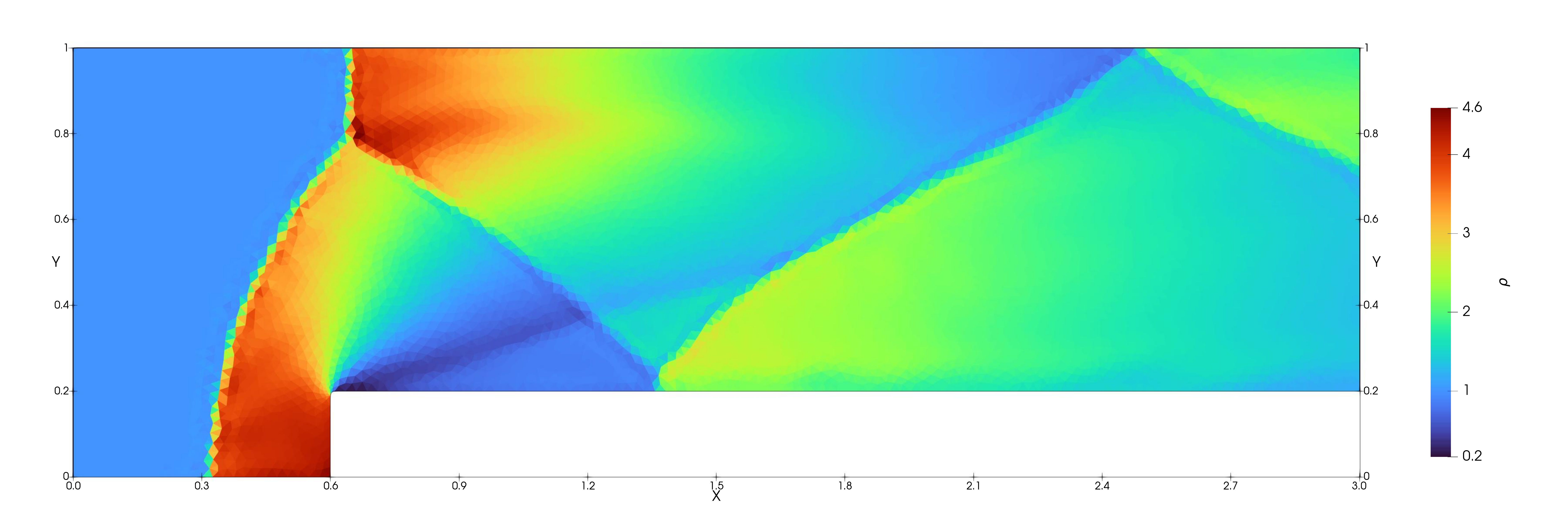}
      \vspace{-4ex}
      \caption{\label{fig:FFS_CWENO-MOOD_entropy} CWENO-MOOD hybrid with entropy trigger}
  \end{subfigure}
  \caption{\label{fig:FFS_DensityPlots} Density profiles for the Forward Facing Step problem at $t=4s$ computed using five methods on medium mesh (12310 triangular cells) and CWENO5 on a finer mesh for reference.}
\end{figure}

The results are presented in the Figure \ref{fig:FFS_DensityPlots}. All five methods produce practically indistinguishable flow profiles, which are a close match with the reference obtained on the mesh with four times higher resolution. 

\subsection{Helium Bubble - Shock Wave Interaction}\label{HeliumBubble}

The accuracy of the schemes when used with diffuse interface method is assessed using a popular numerical and experimental benchmark consisting of the interaction of a weak shock wave in air with a helium bubble. \cite{HeliumBubblePaperExperimental,HeliumBubblePaperSimulation}
The helium bubble is initially contaminated by air up to 28\% concentration by mass, and the initial conditions are defined as in Table \ref{tab:heliumbubbleIC}.
\vspace{-2ex}
\begin{table}[H]
  \centering
  \caption{\label{tab:heliumbubbleIC} Initial Conditions for the helium bubble shock wave test case}
  \vspace{-2ex}
  \begin{tabular}{|l|c|c|c|c|c|c|c|}
    \hline
    Location & $a_1\rho_1$ & $a_2\rho_2$ & $u$ & $v$ & $p$ & $a_1$ \\
    \hline 
    Post-shock &  0.0 & 1.204 & 0 & 0 & 101325 & 0  \\
    Pre-shock &  0.0 & 1.658 & -114.49 & 0 & 159060 & 0  \\
    Bubble    & 0.158 & 0.061 & 0 & 0 & 101325 & 0.95  \\
    \hline
  \end{tabular}
\end{table}
The computations are performed on a pair of 2D unstructured mesh composed of 9232 and 142976 triangular and quadrilateral elements. That corresponds to an edge lengths at the boundaries of the shock bubble interaction zone of $1/20$ and $1/80$ of the initial bubble diameter. Slip-wall boundary conditions are applied to top and bottom of the domain, and inflow and outflow are applied to right and left respectively. The final simulation time is $t_f=983\mu s$.

\begin{figure}[H]
  \vspace{-2ex}
  \centering
  \begin{subfigure}{0.475\textwidth}
    \centering
    \includegraphics[width=\textwidth, trim = {0 4cm 0 6cm}, clip]{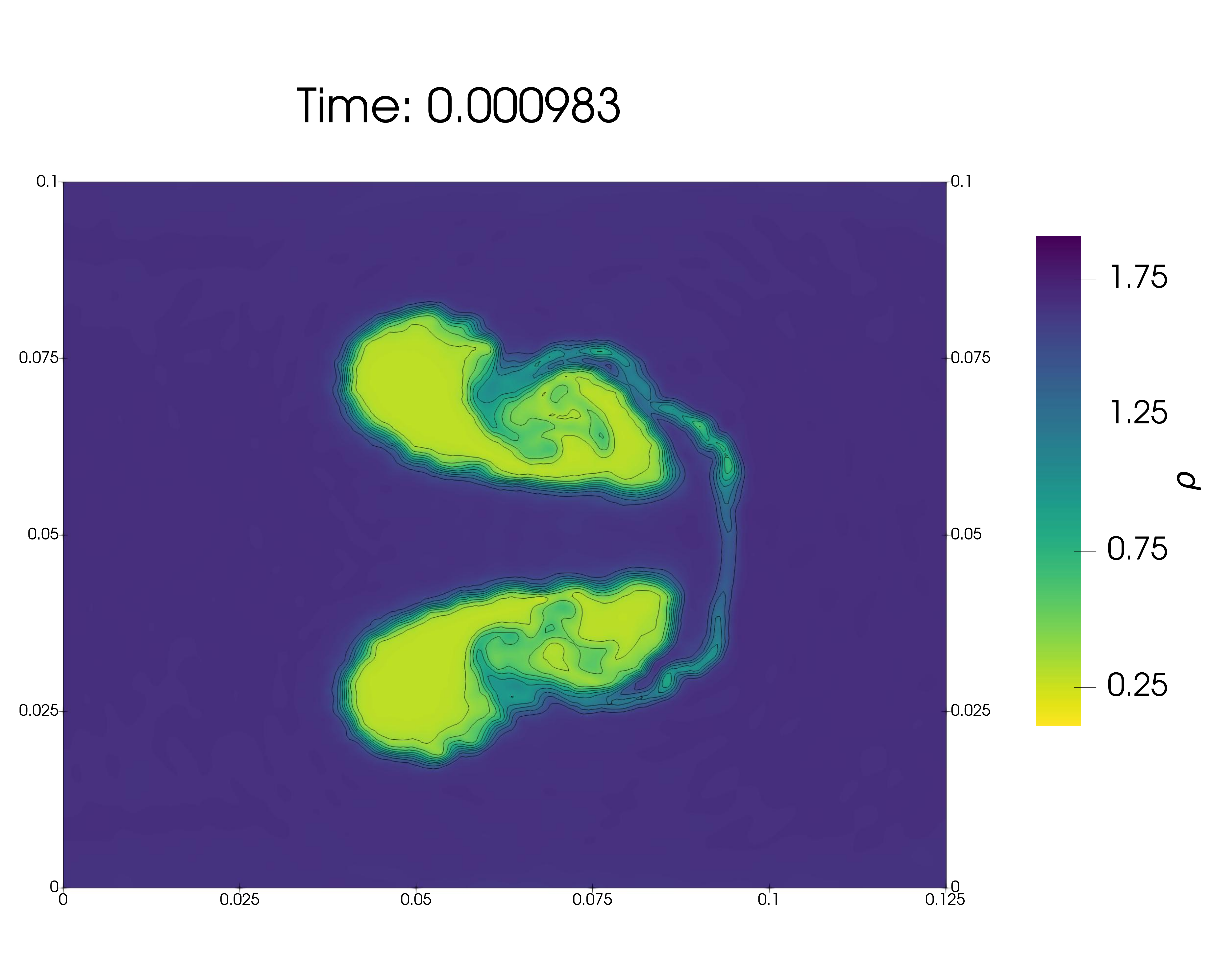}
    \caption{\label{fig:Bubble_CWENO5-superfine} reference solution - CWENO5 on the finest mesh (142976 cells)}
  \end{subfigure}
  \begin{subfigure}{0.475\textwidth}
    \centering
    \includegraphics[width=\textwidth, trim = {0 4cm 0 6cm}, clip]{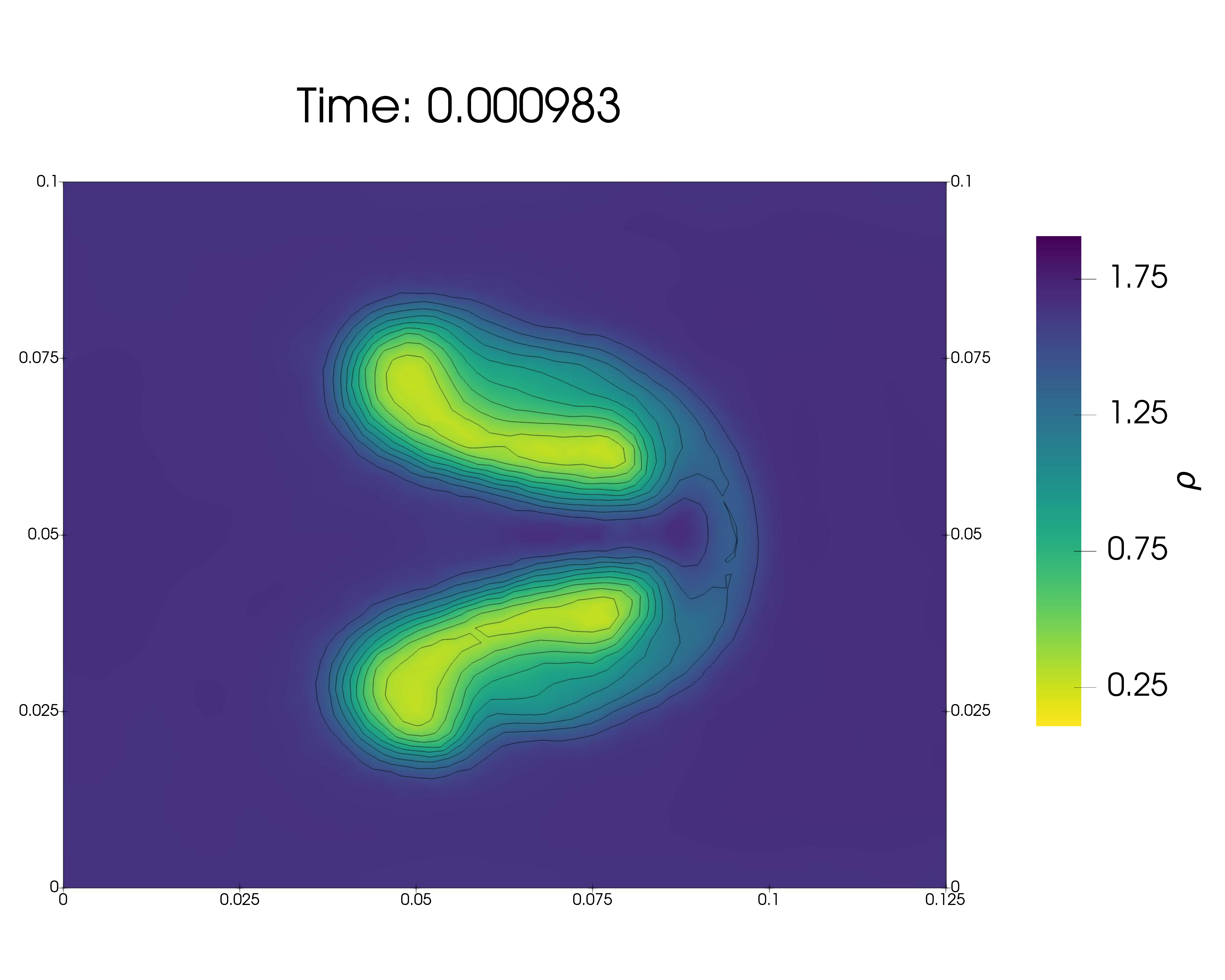}
    \caption{\label{fig:Bubble_CWENO5-medium} reference solution - CWENO5 on medium mesh}
  \end{subfigure}
  \begin{subfigure}{0.475\textwidth}
      \centering
      \includegraphics[width=\textwidth, trim = {0 4cm 0 2cm}, clip]{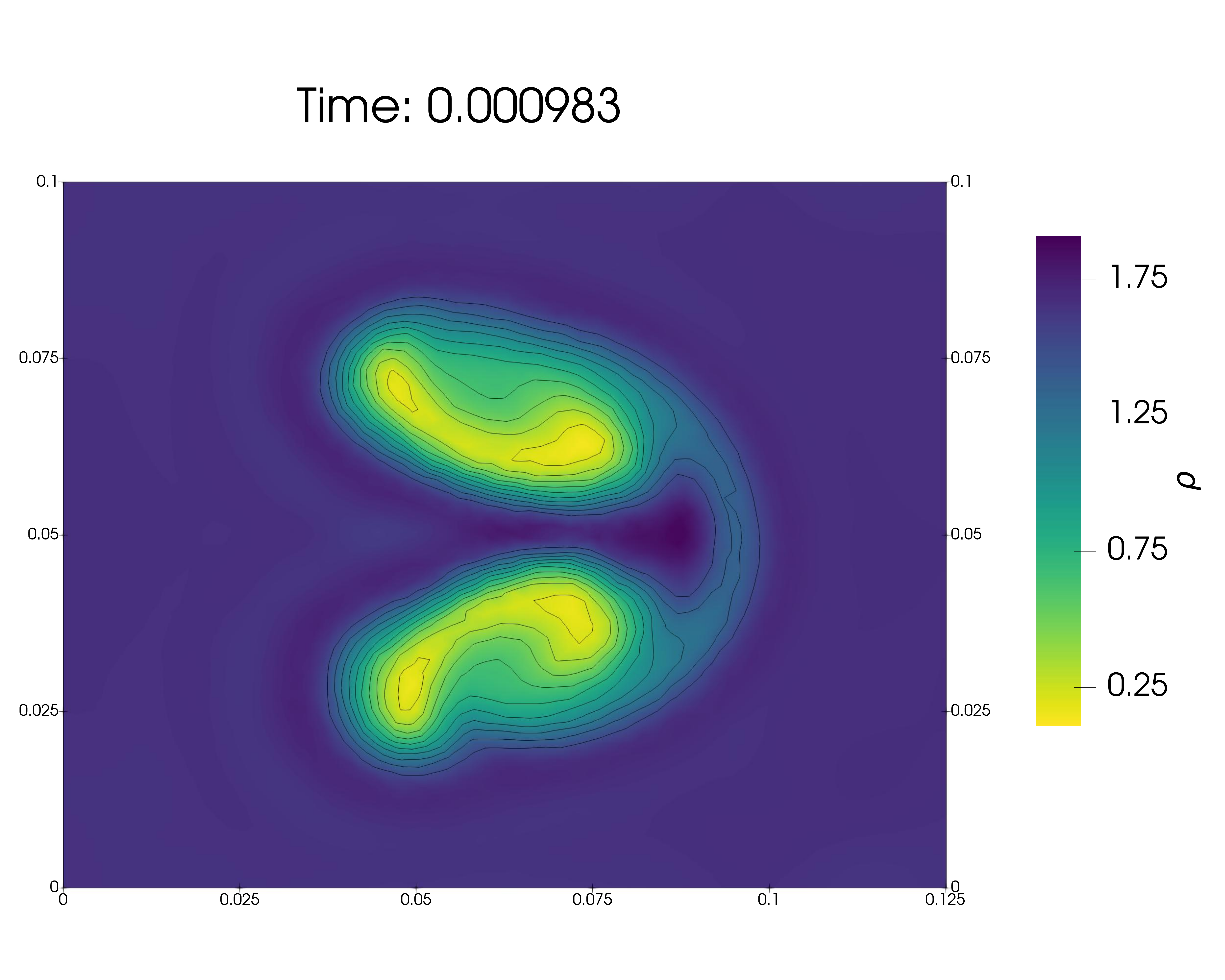}
      \caption{\label{fig:Bubble_ClassicMOOD} Original MOOD}
  \end{subfigure}
  \begin{subfigure}{0.475\textwidth}
      \centering
      \includegraphics[width=\textwidth, trim = {0 4cm 0 2cm}, clip]{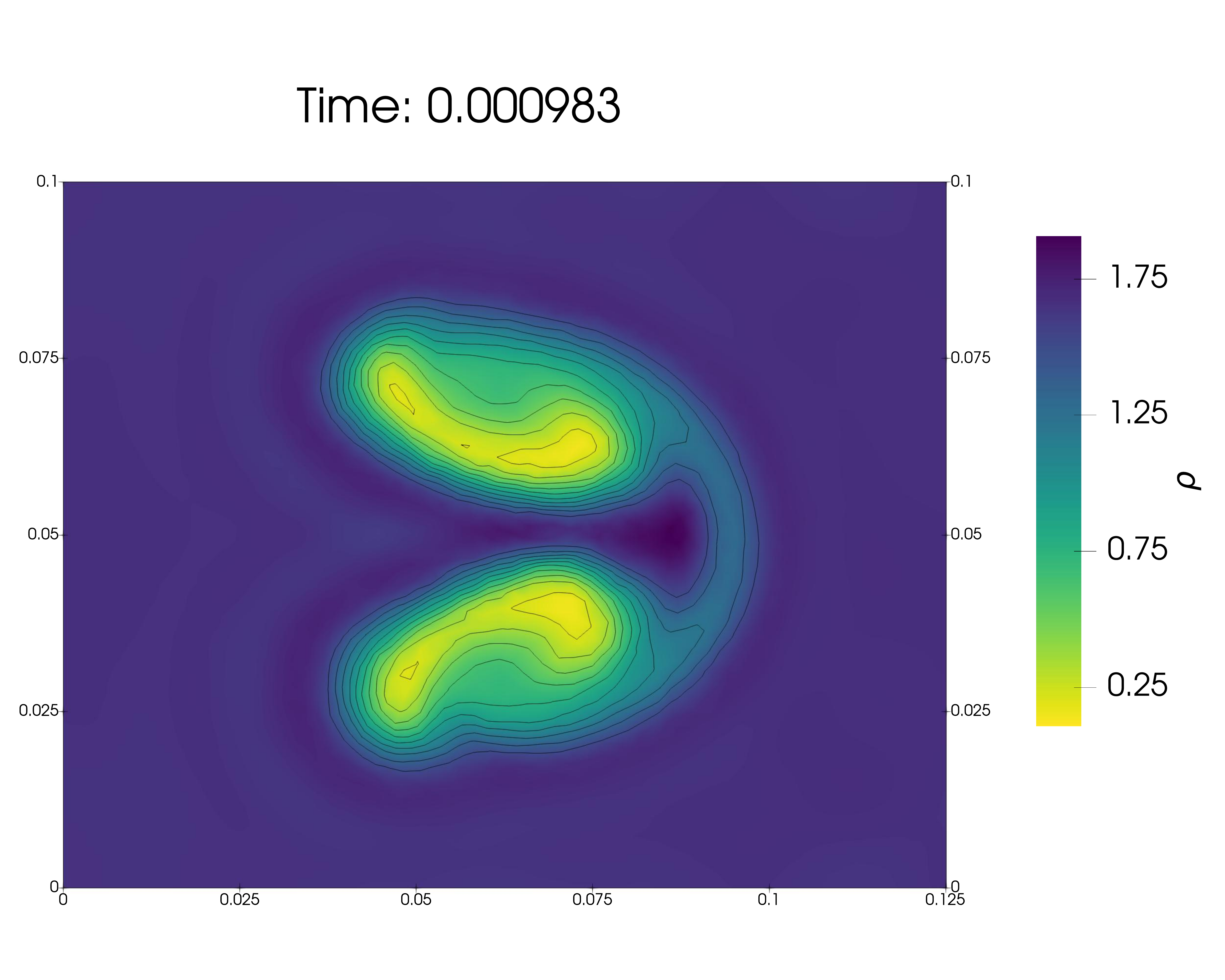}
      \caption{\label{fig:Bubble_EntropyMOOD} Entropy MOOD}
  \end{subfigure}
  \begin{subfigure}{0.475\textwidth}
      \centering
      \includegraphics[width=\textwidth, trim = {0 4cm 0 2cm}, clip]{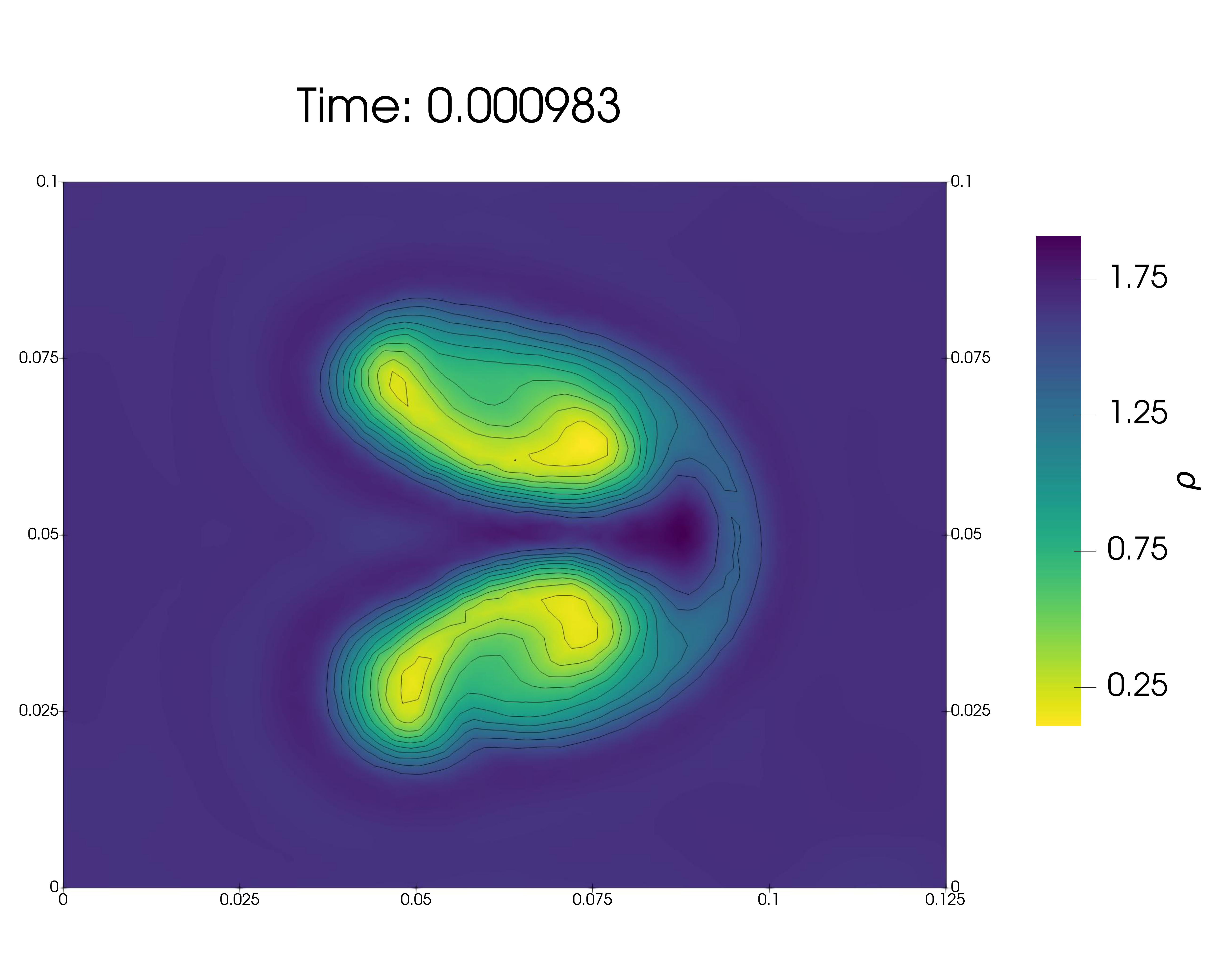}
      \caption{\label{fig:Bubble_CWENO-MOOD_classic} MOOD-CWENO hybrid with original MOOD trigger}
  \end{subfigure}
  \begin{subfigure}{0.475\textwidth}
      \centering
      \includegraphics[width=\textwidth, trim = {0 4cm 0 2cm}, clip]{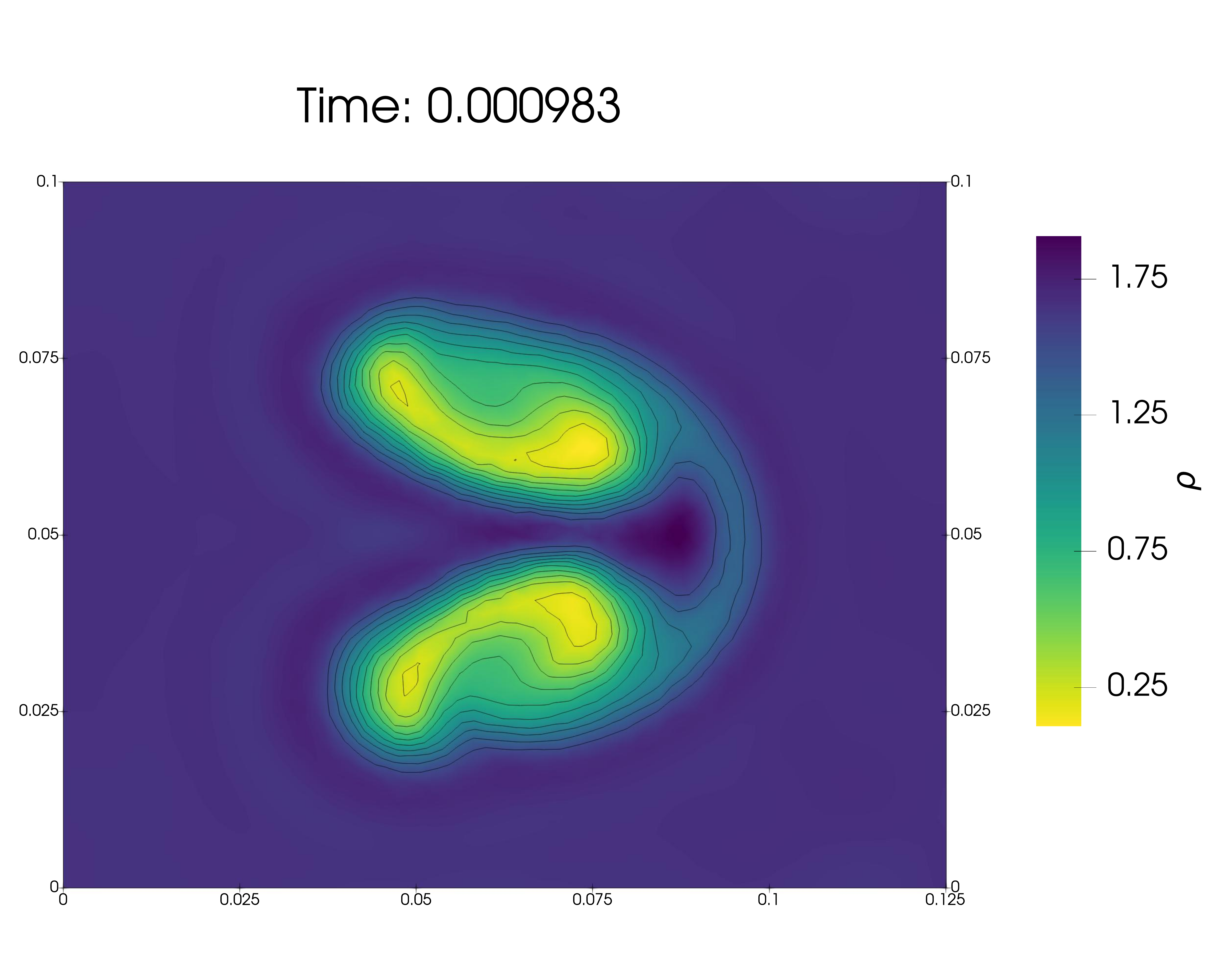}
      \caption{\label{fig:Bubble_CWENO-MOOD_entropy} MOOD-CWENO hybrid with entropy trigger}
  \end{subfigure}
  \caption{\label{fig:Bubble_DensityPlots} Computed density (colour) and volume fraction (contours spaced at $0.1$ intervals) for the shock - helium bubble interaction test problem at the final time $983\ \mu s$ simulated on an unstructured mesh with 9232 cells.}
\end{figure}

As seen in Figure \ref{fig:Bubble_DensityPlots} the four MOOD-based methods produce effectively indistinguishable results. Compared to the CWENO5 solution (\ref{fig:Bubble_CWENO5-medium}) a ring of higher density around the bubble becomes noticeable. The solution on the finer mesh (\ref{fig:Bubble_CWENO5-superfine}) visualises the increase of accuracy a high-order method can produce when the resolution is quadrupled.

\subsection{Triple Material Point}\label{TripleMaterialPoint}

\begin{figure}[H]
    \centering
    \vspace{-1ex}
    \caption{\label{fig:TriplePointInitialDiagram} Initial configuration of the Triple Point problem}
    \vspace{-2.5ex}
    \includegraphics[width=\textwidth, trim = 0cm 6cm 3cm 3.75cm, clip]{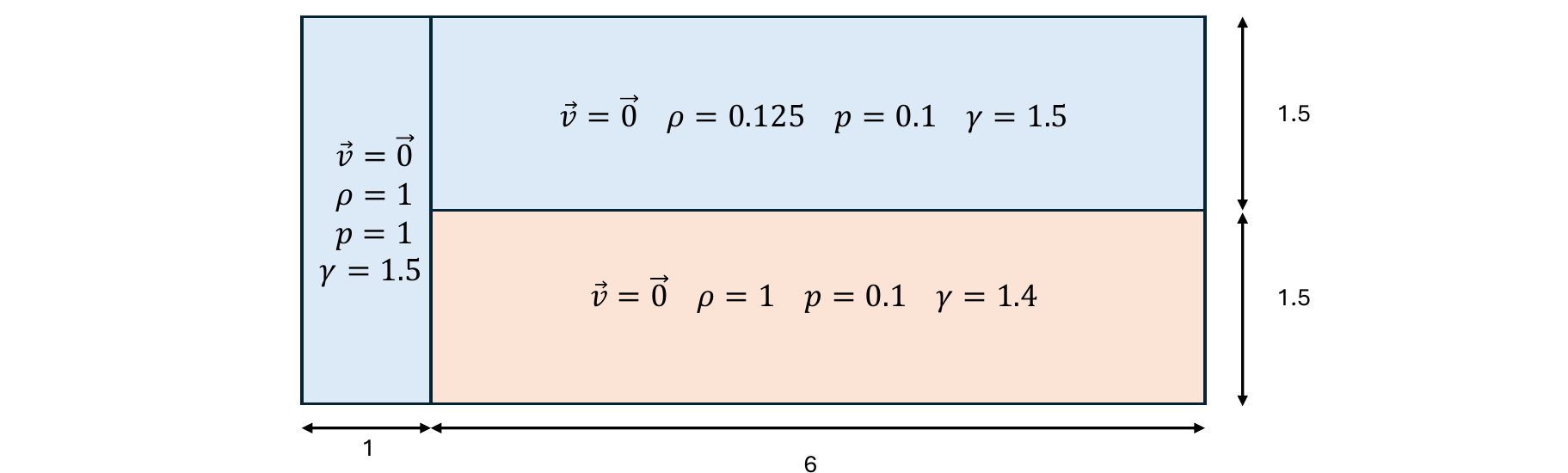}
    \vspace{-3ex}
\end{figure}

The second problem used to test the methods' robustness when simulation multi-material flows is a Triple-Point problem \cite{TriplePointPaper}. In this problem, three gases meet at a point and the interactions between them create an unstable vortex creating a complex flow. The domain is $[0,3]\times[0,7]$ subdivided into a region $1\times3$ and two regions $6\times1.5$ as seen in \textbf{Fig}.\ref{fig:TriplePointInitialDiagram}. All three regions are initially at rest and their initial conditions $(\rho, p, \gamma)$ are $(1, 1, 1.5)$ in the left region, $(0.125, 0.1 1.5)$ in the top right region, and $(1, 0.1, 1.4)$ in the bottom right. Their interaction is simulated until time $t=10s$.
Since the left and bottom regions share the same adiabatic index, the problem is simulated using a two-material diffuse-interface method, and the interface between those two regions is treated as a discontinuity within one material. Two unstructured meshes are used with approximately 3800 and 15000 triangular and quadrilateral elements.
The total simulated time is $10s$ but the intermediate results at $t = 4.5s$ are plotted for comparison with the existing literature.

\begin{figure}[H]
  \centering
  \vspace{-3ex}
  \begin{subfigure}{0.75\textwidth}
    \centering
    \includegraphics[width=\textwidth, trim = {0 2cm 0 1cm}, clip]{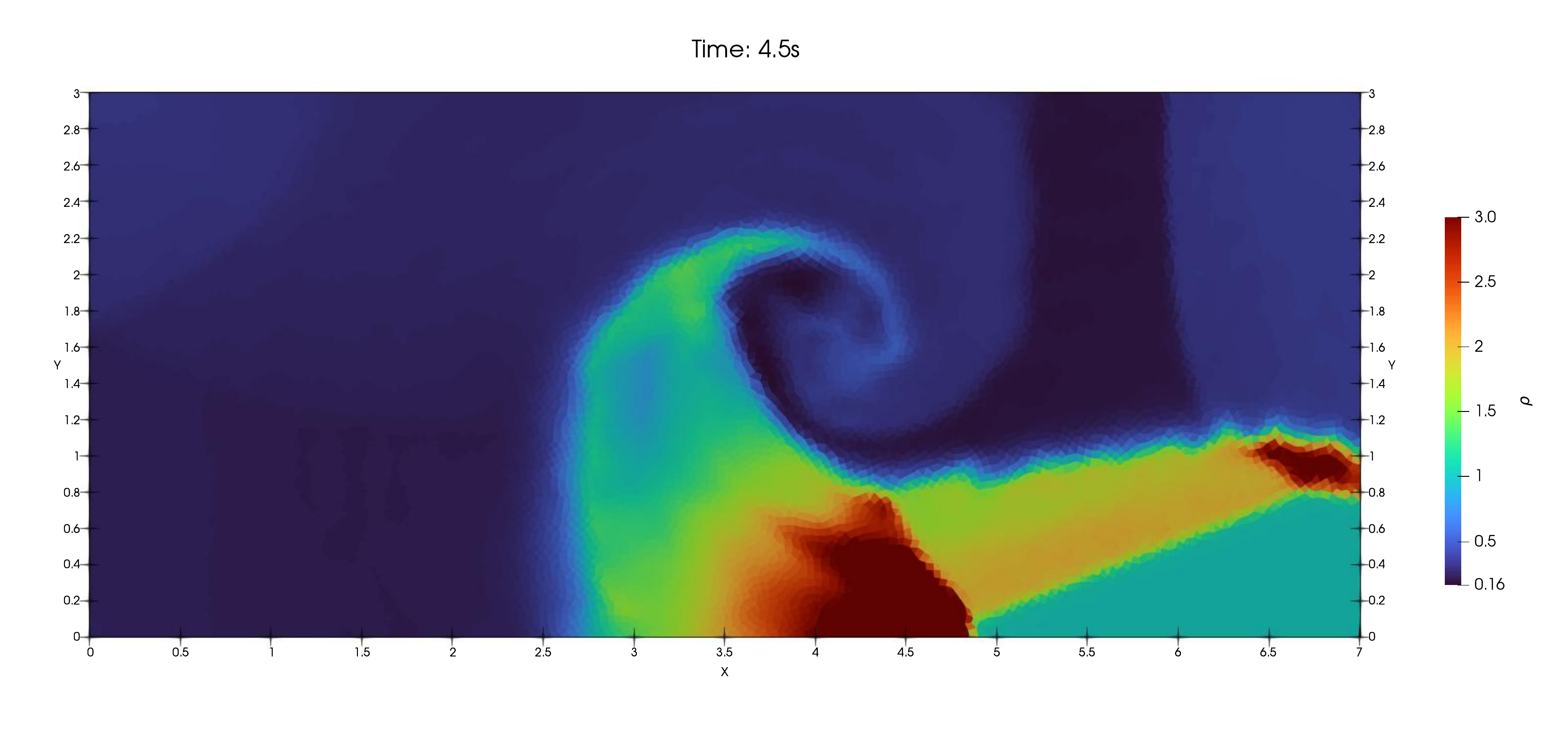}
    \vspace{-3ex}
    \caption{\label{fig:TripleMaterial_CWENO-fine} reference solution - CWENO5 scheme on a finer mesh (15000 cells)}
  \end{subfigure}
  \begin{subfigure}{0.475\textwidth}
      \centering
      \includegraphics[width=\textwidth, trim = {0 2cm 0 0cm}, clip]{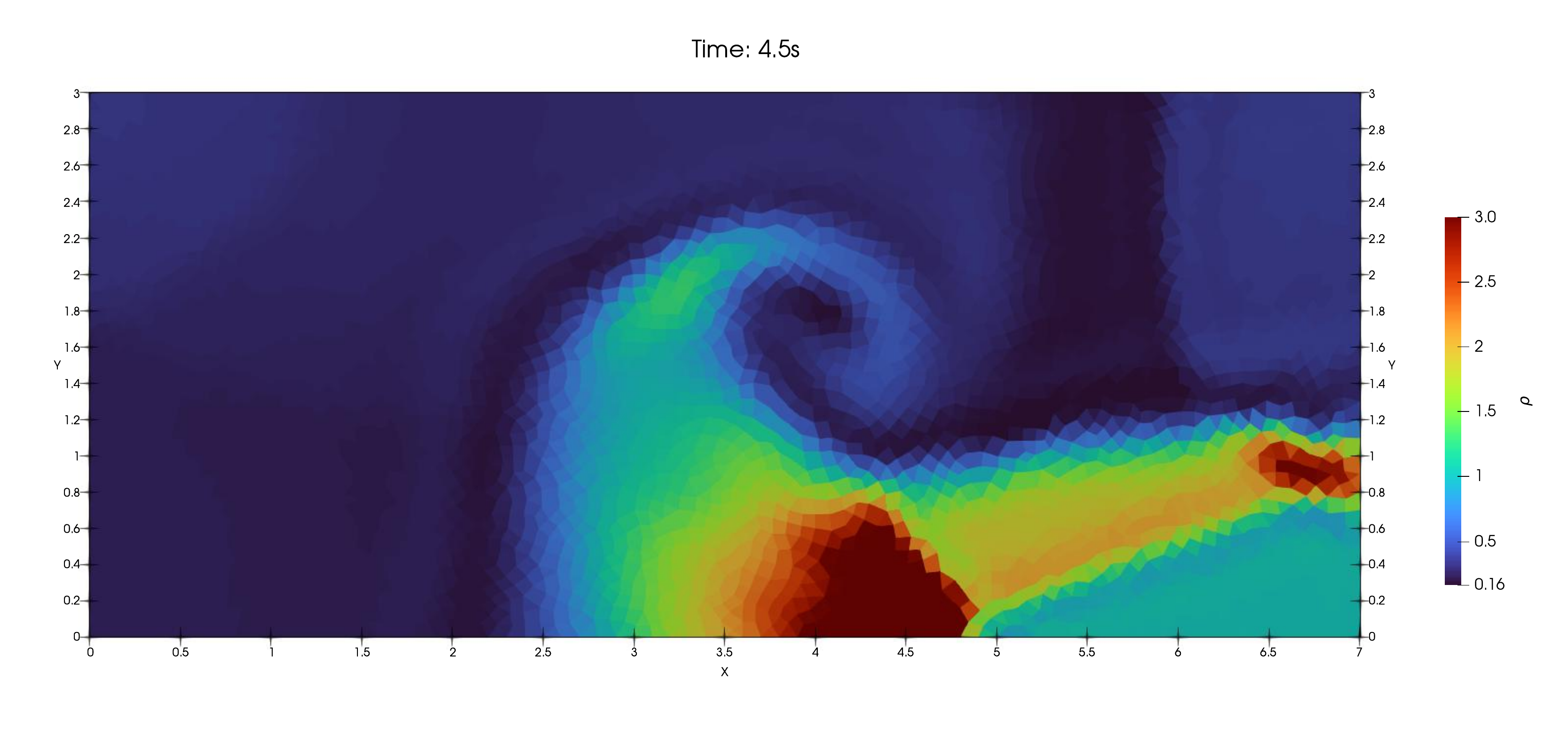}
      \vspace{-3.5ex}
      \caption{\label{fig:TripleMaterial_ClassicMOOD} Original MOOD}
  \end{subfigure}
  \begin{subfigure}{0.475\textwidth}
      \centering
      \includegraphics[width=\textwidth, trim = {0 2cm 0 0cm}, clip]{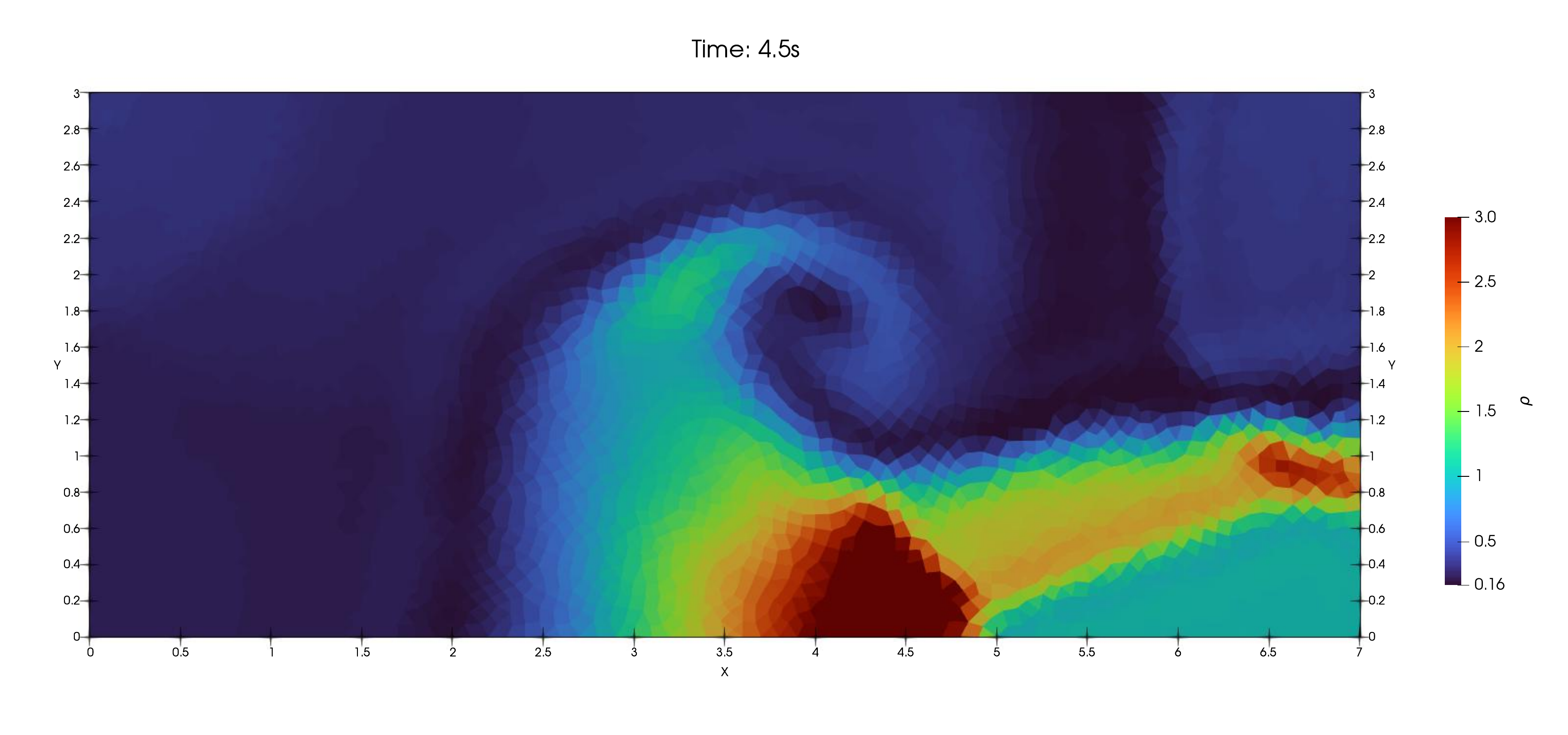}
      \vspace{-3.5ex}
      \caption{\label{fig:TripleMaterial_EntropyMOOD} Entropy MOOD}
  \end{subfigure}
  \begin{subfigure}{0.475\textwidth}
      \centering
      \includegraphics[width=\textwidth, trim = {0 2cm 0 0cm}, clip]{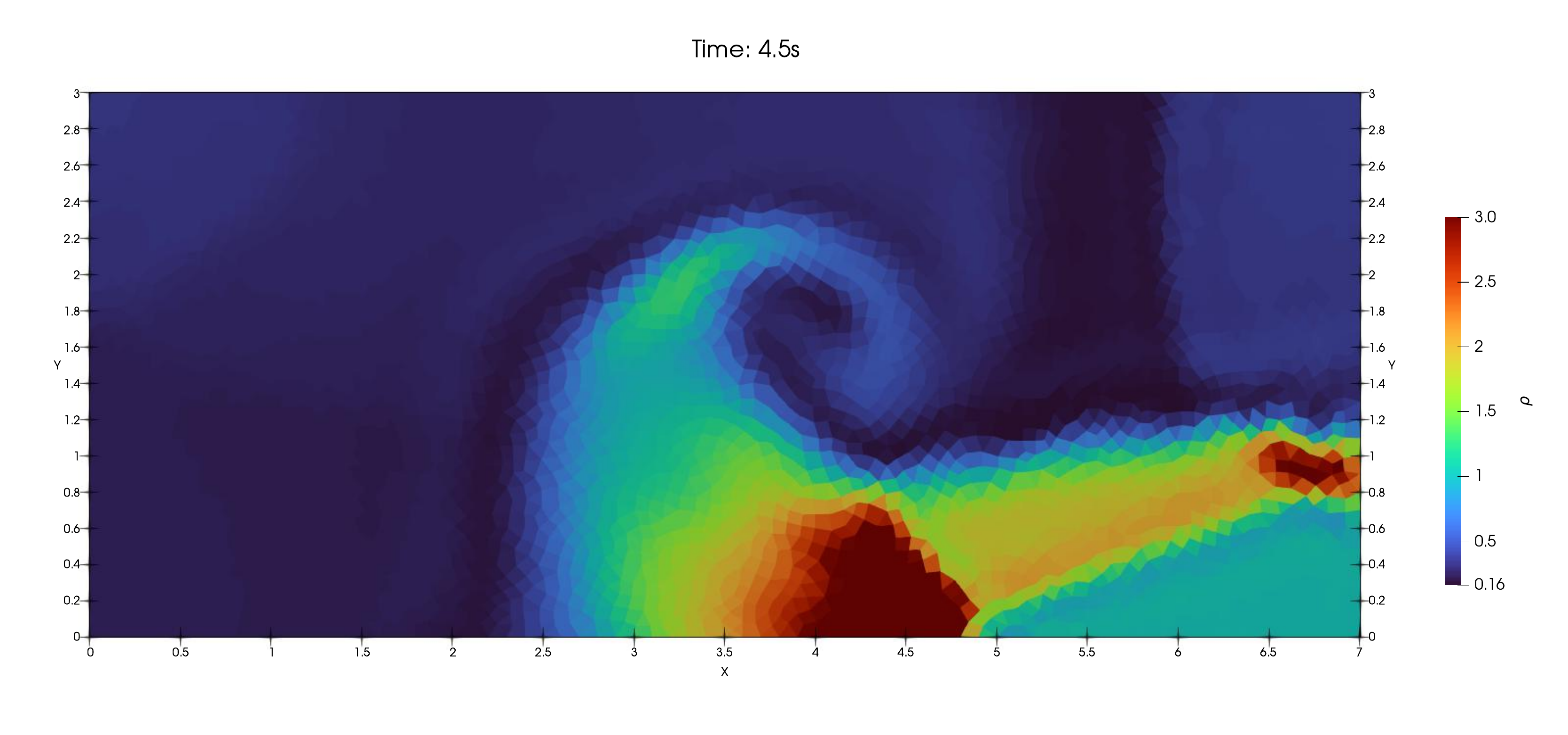}
      \vspace{-3.5ex}
      \caption{\label{fig:TripleMaterial_CWENO-MOOD_classic} MOOD-CWENO hybrid with original MOOD trigger}
  \end{subfigure}
  \begin{subfigure}{0.475\textwidth}
      \centering
      \includegraphics[width=\textwidth, trim = {0 2cm 0 0cm}, clip]{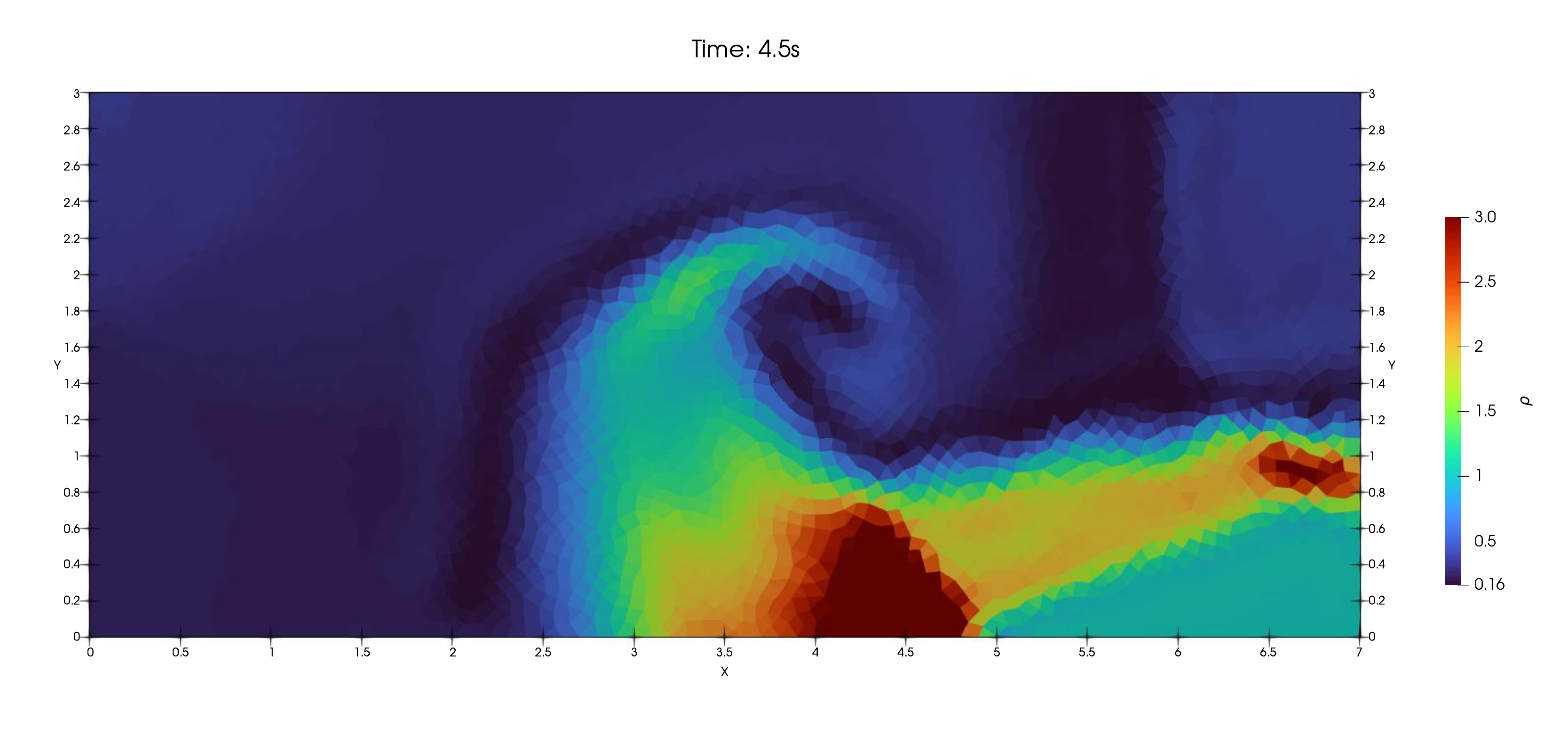}
      \vspace{-3.5ex}
      \caption{\label{fig:TripleMaterial_CWENO-MOOD_entropy} MOOD-CWENO hybrid with entropy trigger}
  \end{subfigure}
  \caption{\label{fig:TripleMaterial_DensityPlots} Computed density for the Triple Material Point test problem at the intermediate time $t = 4.5s$ simulated using the four MOOD methods on an unstructured mesh with 3800 cells and CWENO5 on a finer mesh for reference.}
\end{figure}

As in the previous test the results obtained with the four MOOD-based methods are very similar as seen in the Figure \ref{fig:TripleMaterial_DensityPlots}. Comparing MOOD and the MOOD-CWENO hybrid the later produces less diffused material interface on the left side of the vertex. There does not appear to be any difference between the classic MOOD trigger and the proposed entropy one. Comparing to the reference solution (\ref{fig:TripleMaterial_CWENO-fine}) a band of increased density around the material interface becomes noticeable, especially under the reflected wave on the right-hand side of the domain.

\FloatBarrier
\section{General Performance}\label{General_Performance}
Based on the results of the test cases presented in Section \ref{Applications} it is possible to make more general conclusions about the performance of the proposed numerical methods.

\subsection{Shock detection}\label{Conclusions-shock_detection}

\begin{figure}[H]
    \centering
    \begin{subfigure}{\textwidth}
        \centering
        \includegraphics[width=0.75\textwidth]{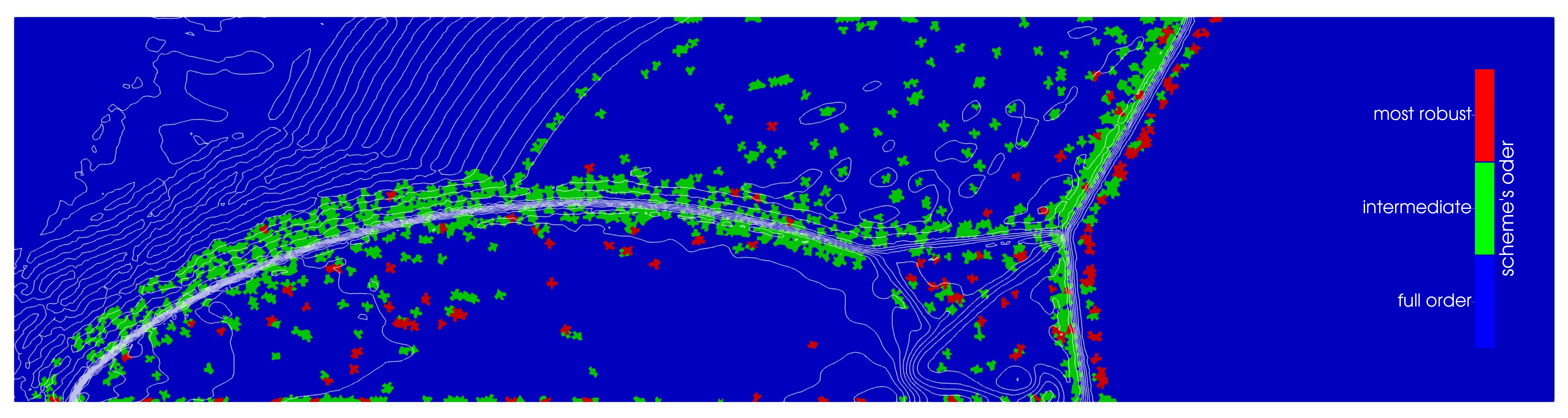}
        \vspace{-1.5ex}
        \caption{\label{fig:DMR_ClassicMOOD_order} Original MOOD}
        \vspace{1ex}
    \end{subfigure}
    \begin{subfigure}{\textwidth}
        \centering
        \includegraphics[width=0.75\textwidth]{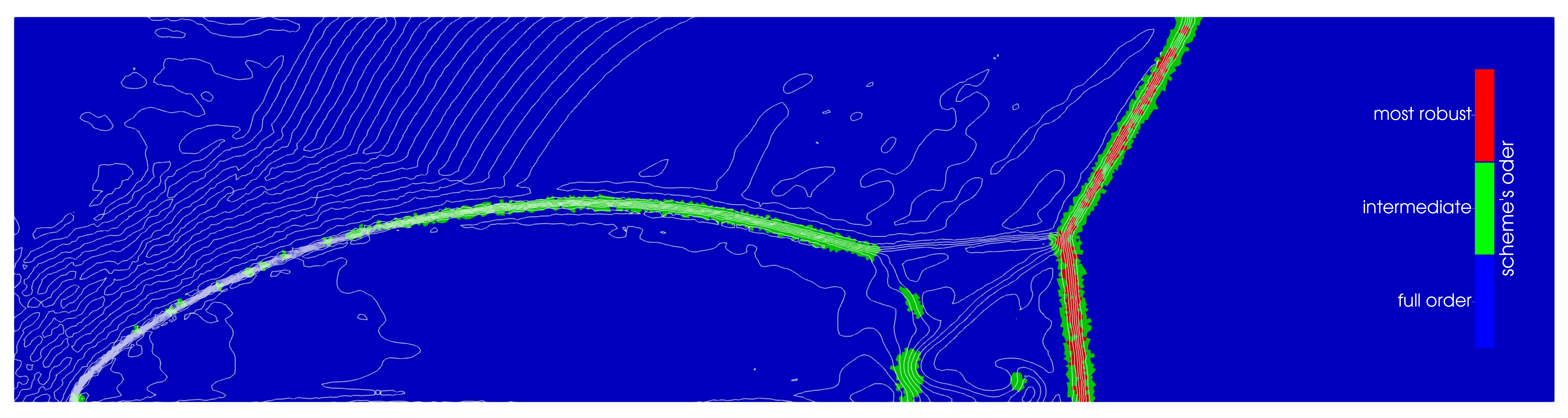}
        \vspace{-1.5ex}
        \caption{\label{fig:DMR_EntropyMOOD_order} Entropy MOOD}
        \vspace{1ex}
    \end{subfigure}
    \begin{subfigure}{\textwidth}
        \centering
        \includegraphics[width=0.75\textwidth]{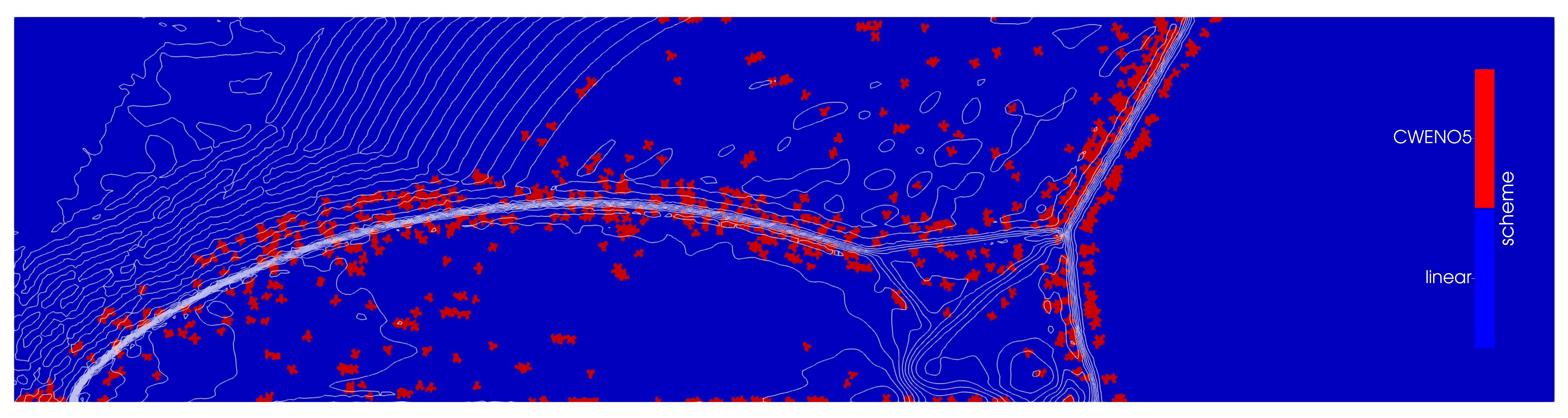}
        \vspace{-1.5ex}
        \caption{\label{fig:DMR_CWENO-MOOD_classic_order} MOOD-CWENO hybrid with original MOOD trigger}
        \vspace{1ex}
    \end{subfigure}
    \begin{subfigure}{\textwidth}
        \centering
        \includegraphics[width=0.75\textwidth]{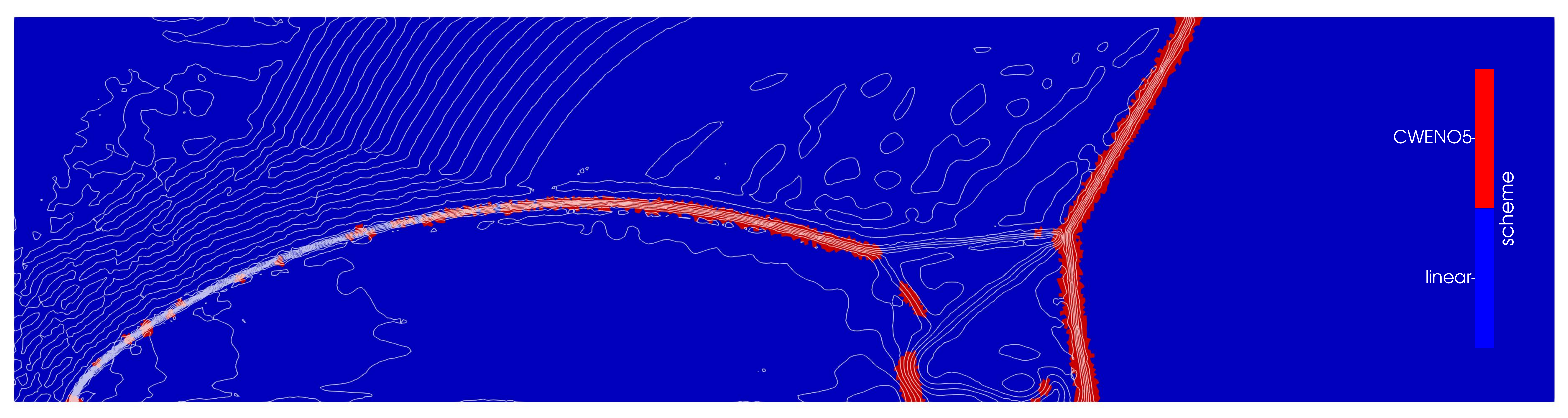}
        \vspace{-1.5ex}
        \caption{\label{fig:DMR_CWENO-MOOD_entropy_order} MOOD-CWENO hybrid with entropy trigger}
    \end{subfigure}
    \caption{\label{fig:DMR_orderPlot} The cells which have been flagged by MOOD triggers in the last time step of Double Mach Reflection simulation. The blue cells have not been flagged and used high-order linear reconstruction. In the upper two plots the green cells have been flagged once and used $2^{nd}$ order reconstruction and the red cells were flagged twise and used $1^{st}$ order scheme. In MOOD-CWENO hybrid plots red designates cells which have been flagged and full CWENO5 reconstruction was performed. For the corresponding simulation results see Figure \ref{fig:DMR_DensityPlots}}
\end{figure}

Comparing which cells are flagged the two decision algorithms produce substantially different patterns as seen in Figure \ref{fig:DMR_orderPlot}. The proposed entropy based trigger flags cells lying in clean lines along shocks and discontinuities, whereas cells flagged by the original detector are scattered. It is caused by the trigger flagging predominantly cells in front and behind, but not inside the shock, which happens because the trigger compares the potential solution to the local maximum and minimum, which are far apart inside the shock, so it is easy for the potential solution to pass the check. On the other hand, in the smooth flow regions the minimum and maximum are close together, which sometimes causes the trigger to flag additional, spurious cells.
On the contrary, the new entropy-production-based trigger flags only the cells inside the shock waves, because those are the only regions where substantial amounts of numerical entropy are produced. Unless a significant numerical error has occurred, very entropy should be produced in the smooth regions, so those cells would not be flagged, exactly as desired.

\begin{figure}[H]
    \vspace{-3ex}
    \centering
    \begin{subfigure}{0.49\textwidth}
        \centering
        \includegraphics[width=\textwidth, trim = {4cm 3cm 0 2cm}, clip]{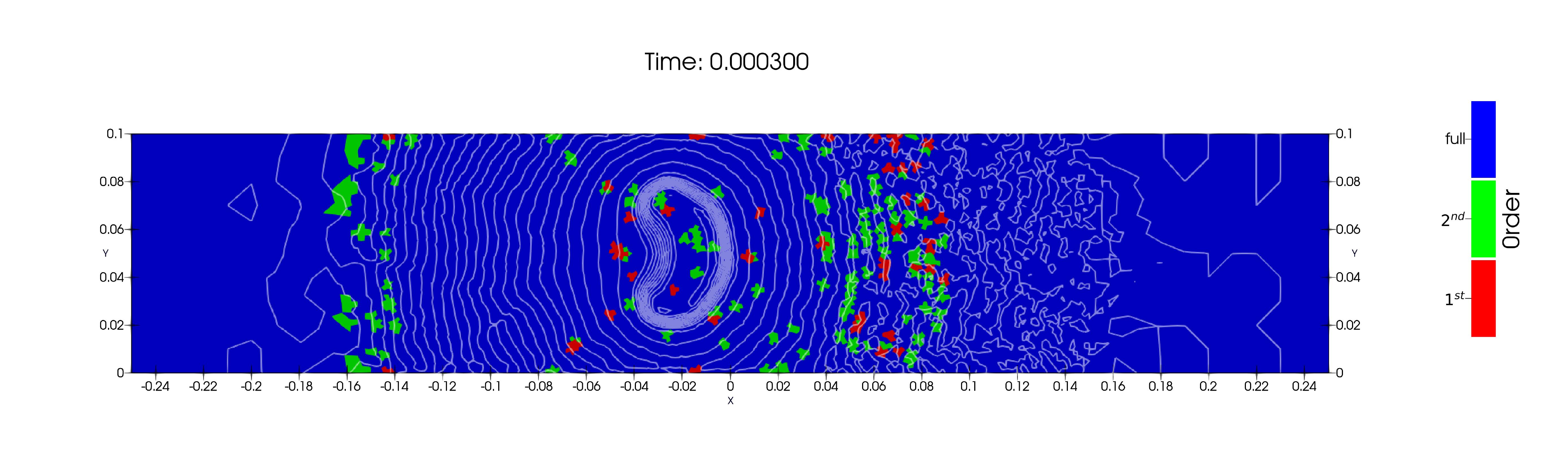}
        \caption{\label{fig:Bubble_ClassicMOOD_order} Original MOOD}
    \end{subfigure}
    \begin{subfigure}{0.49\textwidth}
        \centering
        \includegraphics[width=\textwidth, trim = {4cm 3cm 0 2cm}, clip]{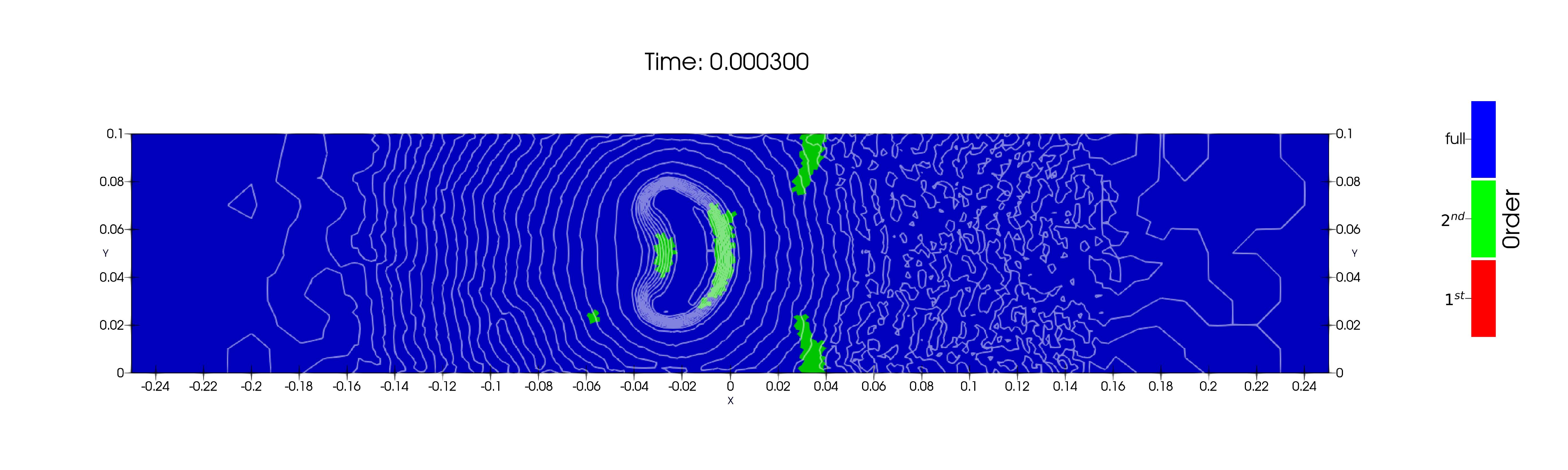}
        \caption{\label{fig:Bubble_EntropyMOOD_order} Entropy MOOD}
    \end{subfigure}
    \begin{subfigure}{0.49\textwidth}
        \centering
        \includegraphics[width=\textwidth, trim = {4cm 3cm 0 1cm}, clip]{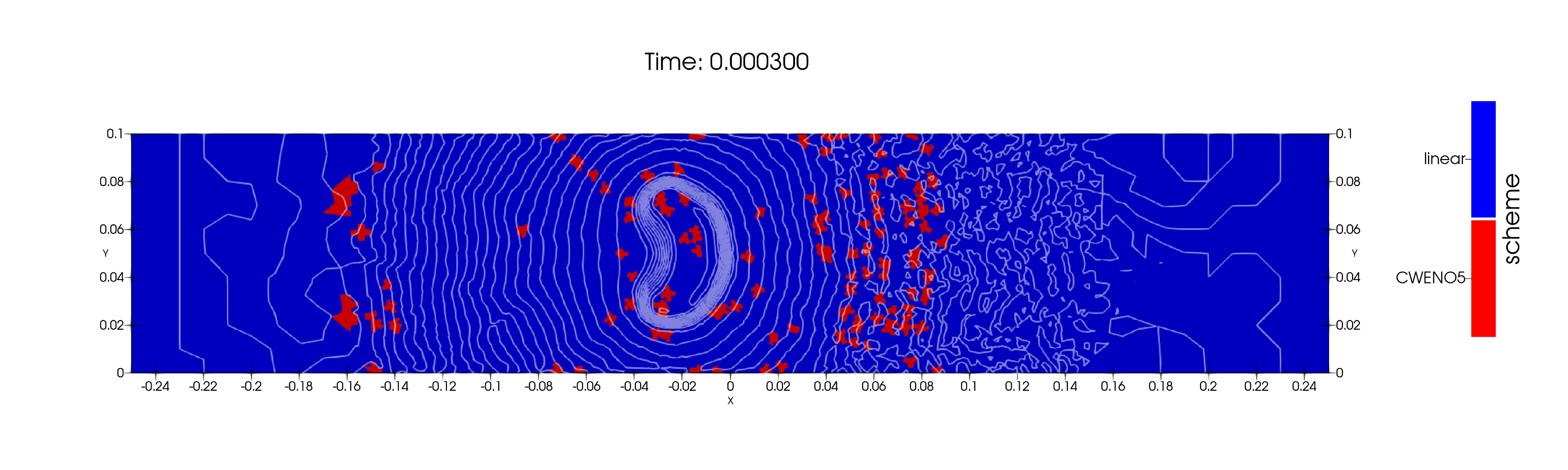}
        \caption{\label{fig:Bubble_CWENO-MOOD_classic_order} MOOD-CWENO hybrid with original MOOD trigger}
    \end{subfigure}
    \begin{subfigure}{0.49\textwidth}
        \centering
        \includegraphics[width=\textwidth, trim = {4cm 3cm 0 1cm}, clip]{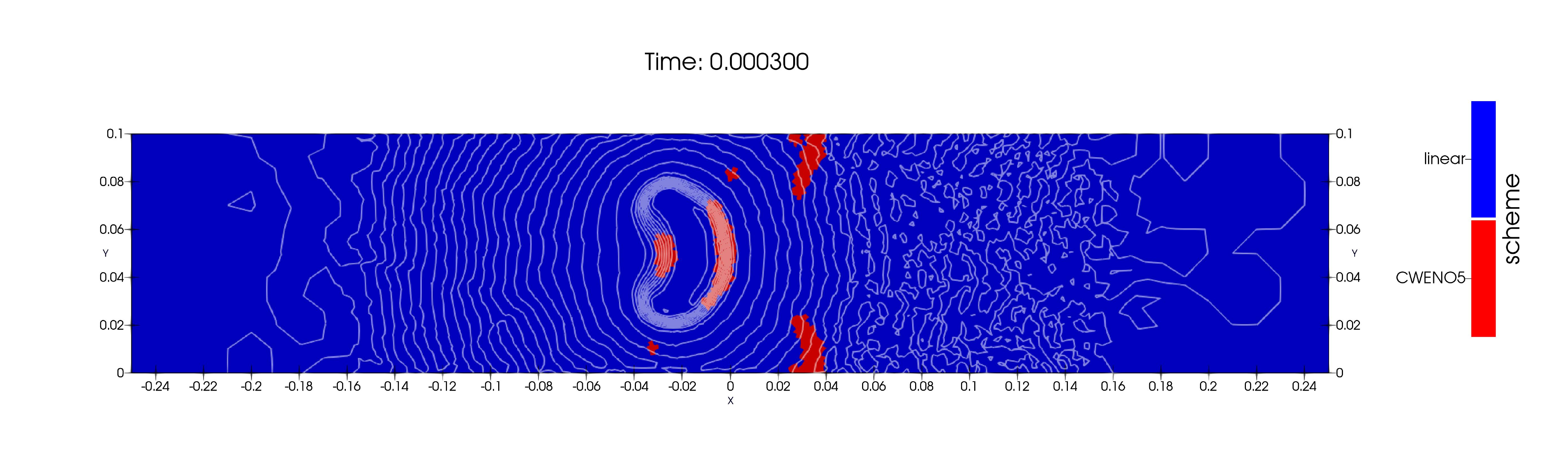}
        \caption{\label{fig:Bubble_CWENO-MOOD_entropy_order} MOOD-CWENO hybrid with entropy trigger}
    \end{subfigure}
    \caption{\label{fig:Bubble_orderPlot} The polynomial reconstruction methods which have been used in each cell during helium bubble simulation. Blue cells have not been flagged by MOOD trigger at all. For context, uniformly spaced density contours are overlaid on the plots.}
\end{figure}
\begin{figure}[H]
    \vspace{-3ex}
    \centering
    \begin{subfigure}{0.49\textwidth}
        \centering
        \includegraphics[width=\textwidth, trim = {2cm 3cm 0 2cm}, clip]{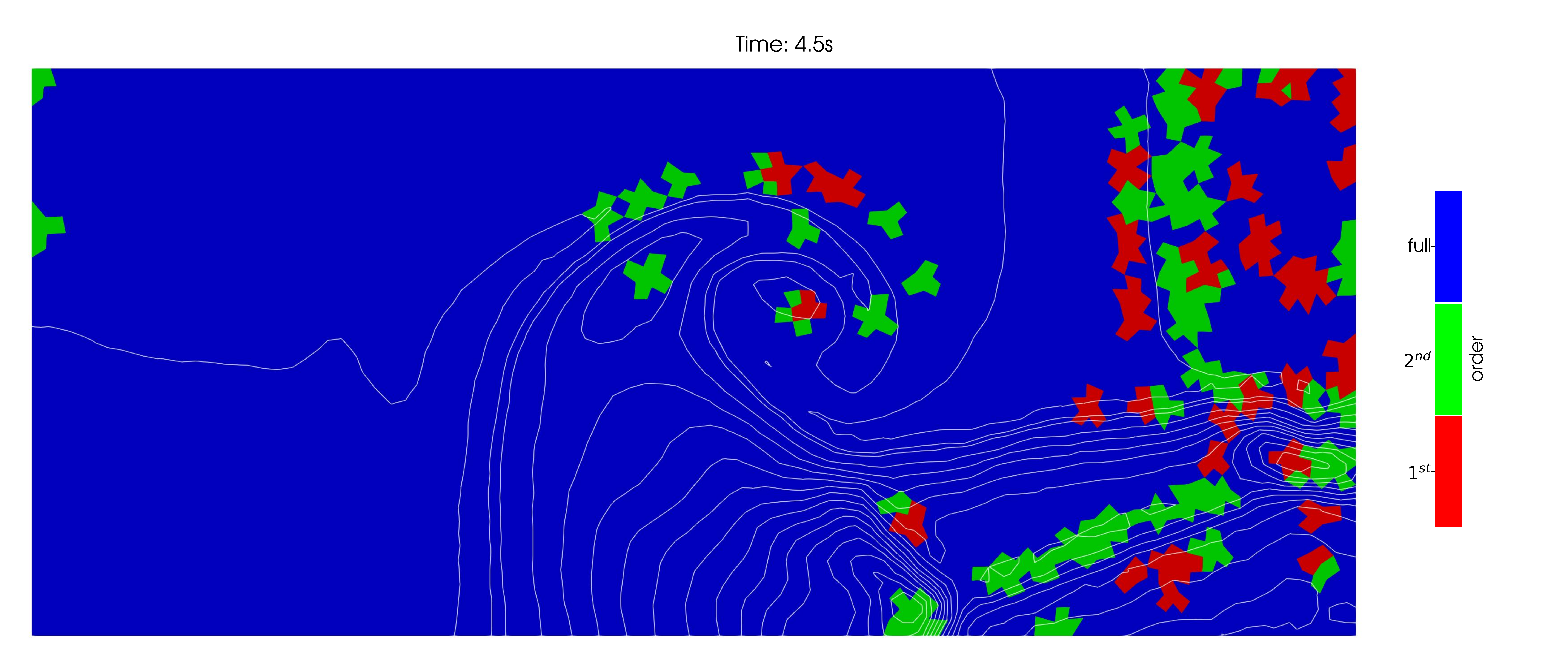}
        \caption{\label{fig:TripleMaterial_ClassicMOOD_order} Original MOOD}
    \end{subfigure}
    \begin{subfigure}{0.49\textwidth}
        \centering
        \includegraphics[width=\textwidth, trim = {0 3cm 2cm 2cm}, clip]{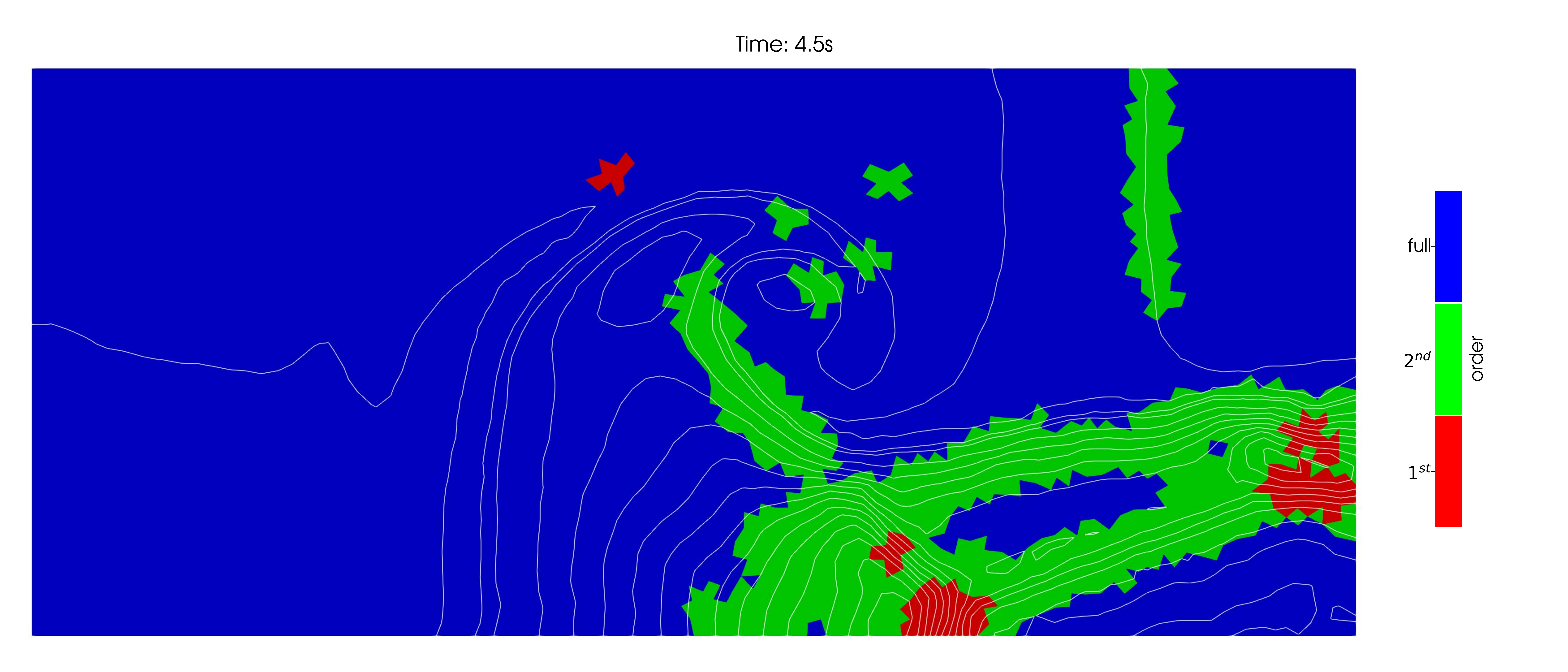}
        \caption{\label{fig:TripleMaterial_EntropyMOOD_order} Entropy MOOD}
    \end{subfigure}
    \begin{subfigure}{0.49\textwidth}
        \centering
        \includegraphics[width=\textwidth, trim = {2cm 3cm 0 1cm}, clip]{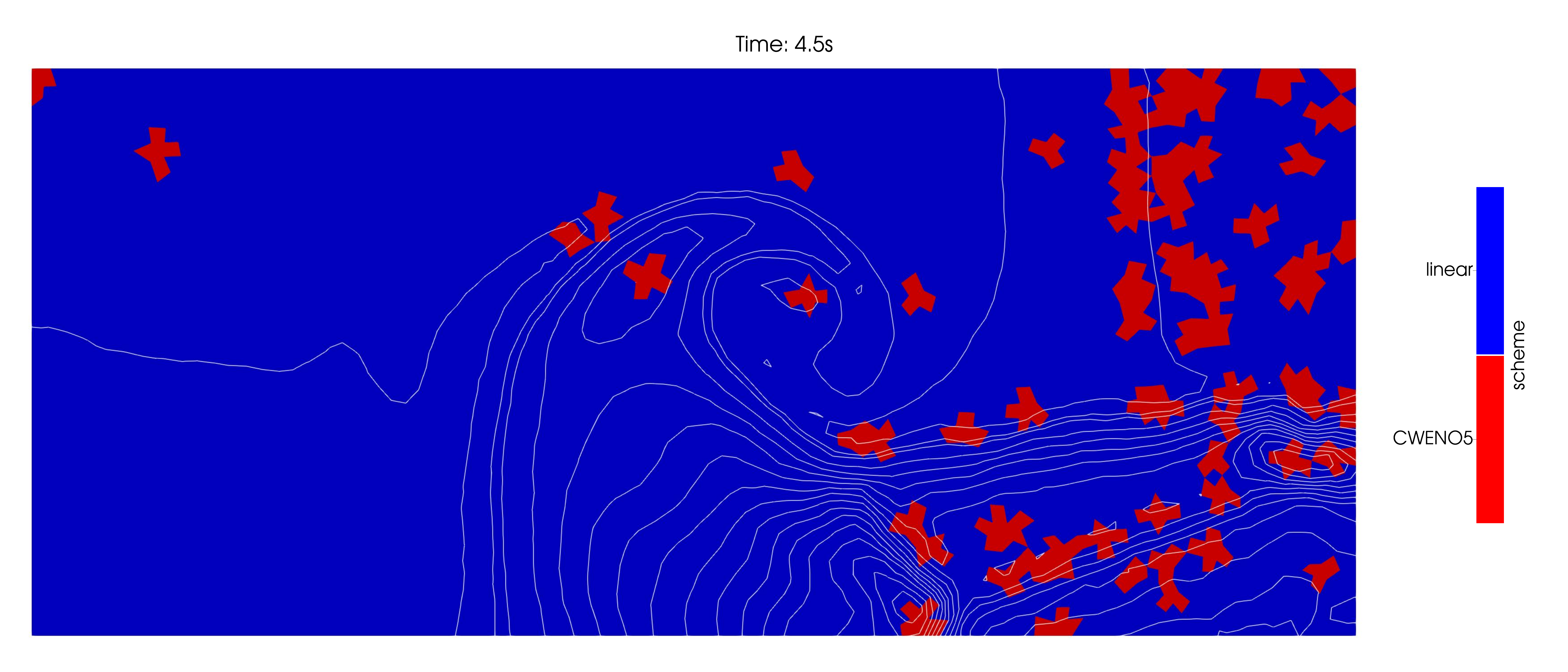}
        \caption{\label{fig:TripleMaterial_CWENO-MOOD_classic_order} MOOD-CWENO hybrid with original MOOD trigger}
    \end{subfigure}
    \begin{subfigure}{0.49\textwidth}
        \centering
        \includegraphics[width=\textwidth, trim = {0 3cm 2cm 1cm}, clip]{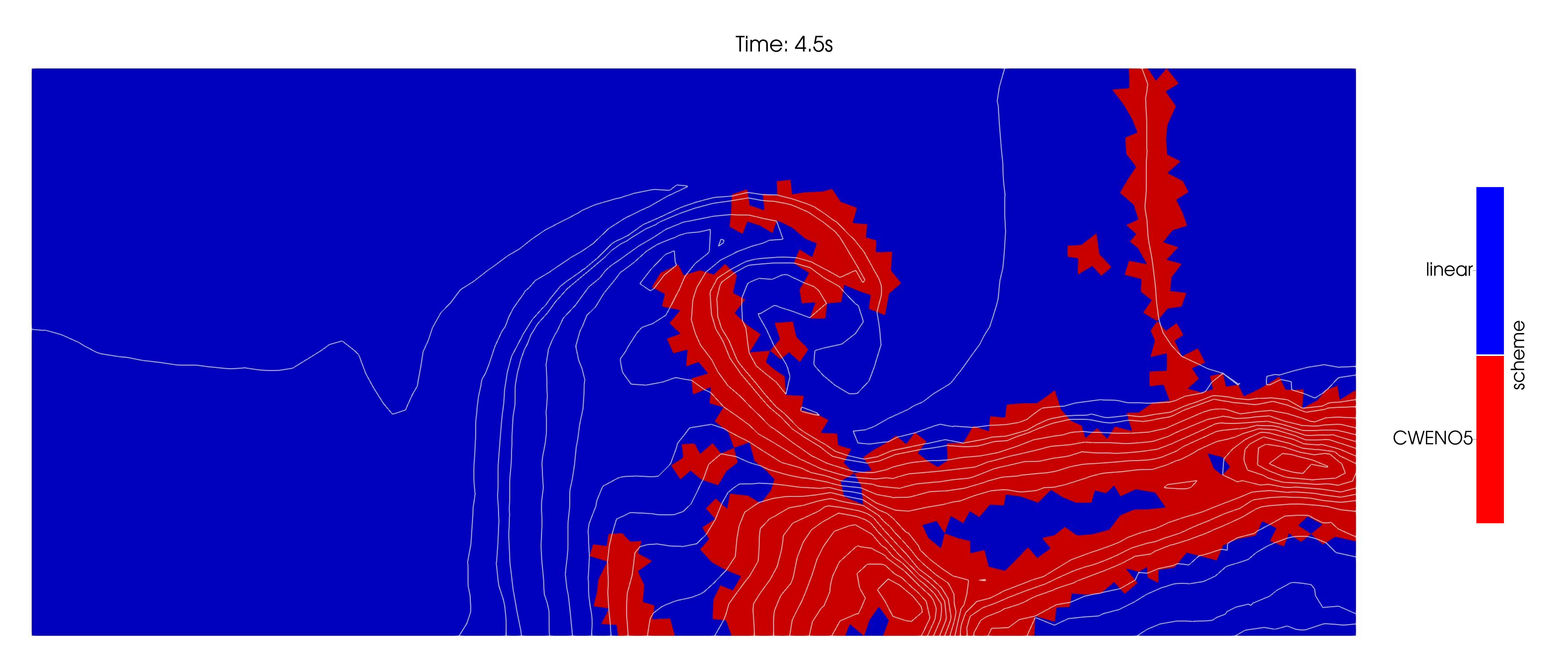}
        \caption{\label{fig:TripleMaterial_CWENO-MOOD_entropy_order} MOOD-CWENO hybrid with entropy trigger}
    \end{subfigure}
    \caption{\label{fig:TripleMaterial_orderPlot} The polynomial reconstruction methods which have been used in each cell during Triple Material Point simulation. Blue cells have not been flagged by MOOD trigger at all. For context, uniformly spaced density contours are overlaid on the plots. The corresponding simulation results are found in the Figure \ref{fig:TripleMaterial_DensityPlots}}
\end{figure}

The exact same difference in patterns is also seen in the multi-phase flow simulations - Figures \ref{fig:Bubble_orderPlot} and \ref{fig:TripleMaterial_orderPlot}.
The entropy production trigger flags cells correctly in both phases, proving that the simple, linear blend of two phases' entropies scales correctly and does not cause the trigger to flag calls spuriously in either phase.

\subsection{MOOD-CWENO hybrid accuracy}\label{Conclusion-CWENO-MOOD_hybrid_accuracy}

The MOOD-CWENO hybrid aims to provide accuracy improvement over the traditional MOOD by eliminating the $2^{nd}$ and $1^{st}$ order methods used in flagged cells and replacing them with high-order, robust scheme. In the tests presented so far the effect is not very pronounced, because the number of flagged cells, hence the number of cells where the difference is applicable, is kept to minimum ($<5\%$ of cells). It can be noticed in the sharper material interface in the Triple Point test problem (compare Figures \ref{fig:TripleMaterial_CWENO-MOOD_classic} and \ref{fig:TripleMaterial_CWENO-MOOD_entropy} to \ref{fig:TripleMaterial_ClassicMOOD} \& \ref{fig:TripleMaterial_EntropyMOOD}). And the negative interaction of the high order scheme with the scattered of cells flagged by the classic trigger is likely the cause of wrong shock position in the Double Mach Reflection test (see Figure \ref{fig:DMR_CWENO-MOOD_classic}).

\begin{figure}[!h]
    \centering
    \begin{subfigure}{0.475\textwidth}
        \centering
        \includegraphics[width=\textwidth, trim = {0 3cm 0 6cm}, clip]{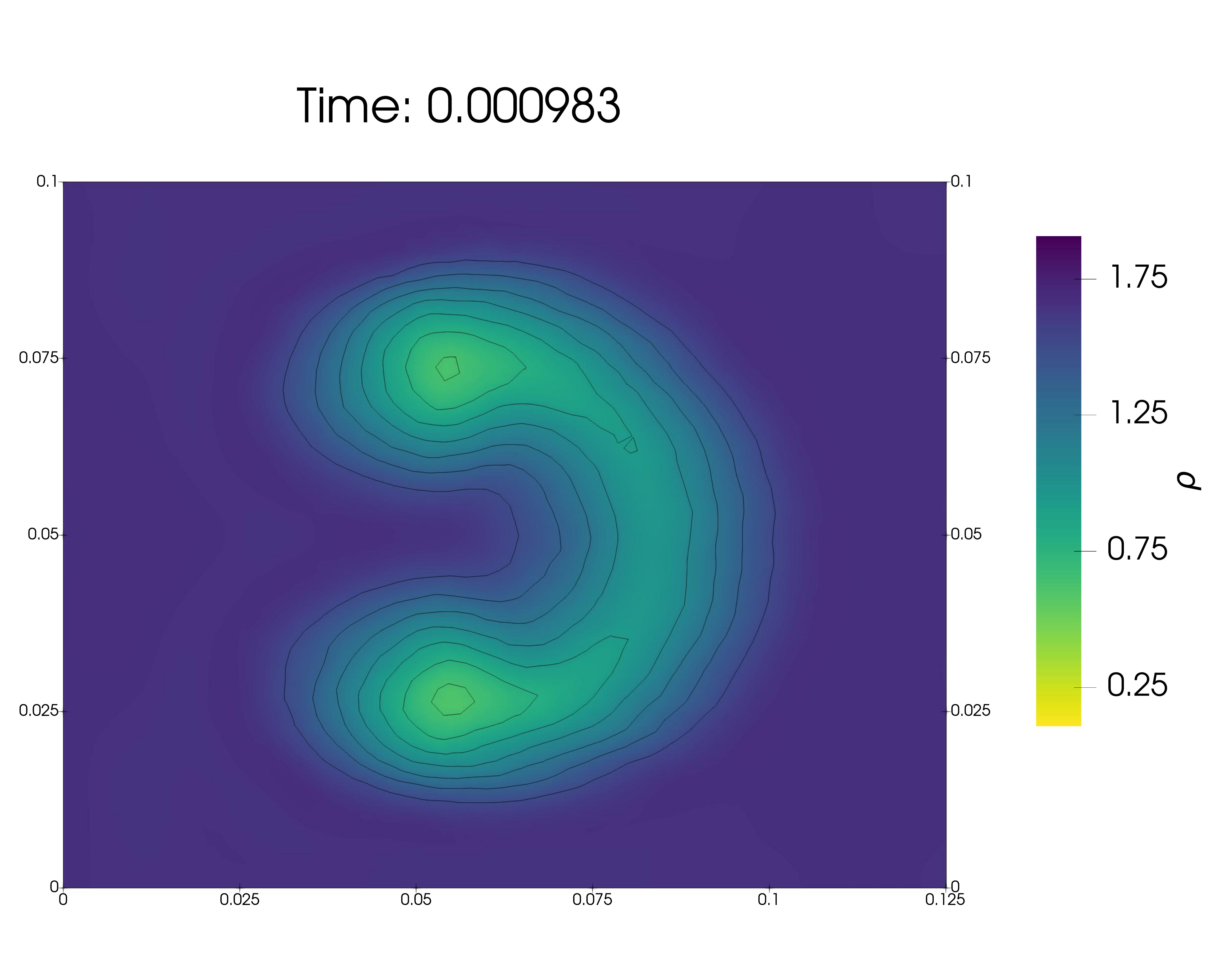}
        \includegraphics[width=\textwidth, trim = {4cm 3cm 0 2cm}, clip]{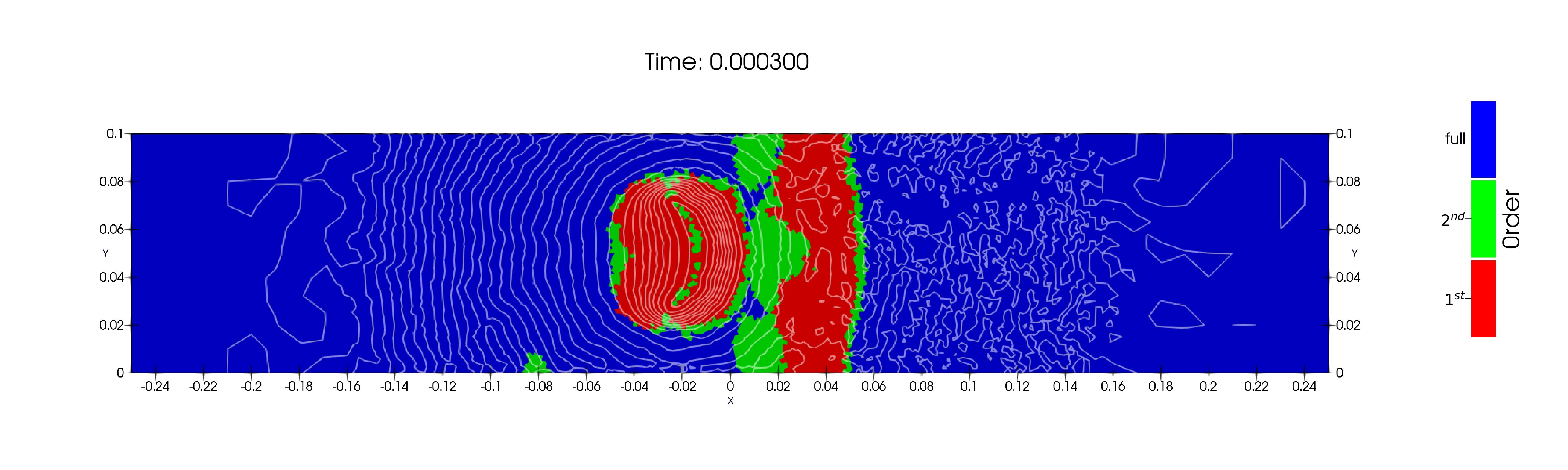}
        \caption{\label{fig:Bubble_MOODvsHybrid-MOOD} Entropy MOOD}
    \end{subfigure}
    \begin{subfigure}{0.475\textwidth}
        \centering
        \includegraphics[width=\textwidth, trim = {0 3cm 0 6cm}, clip]{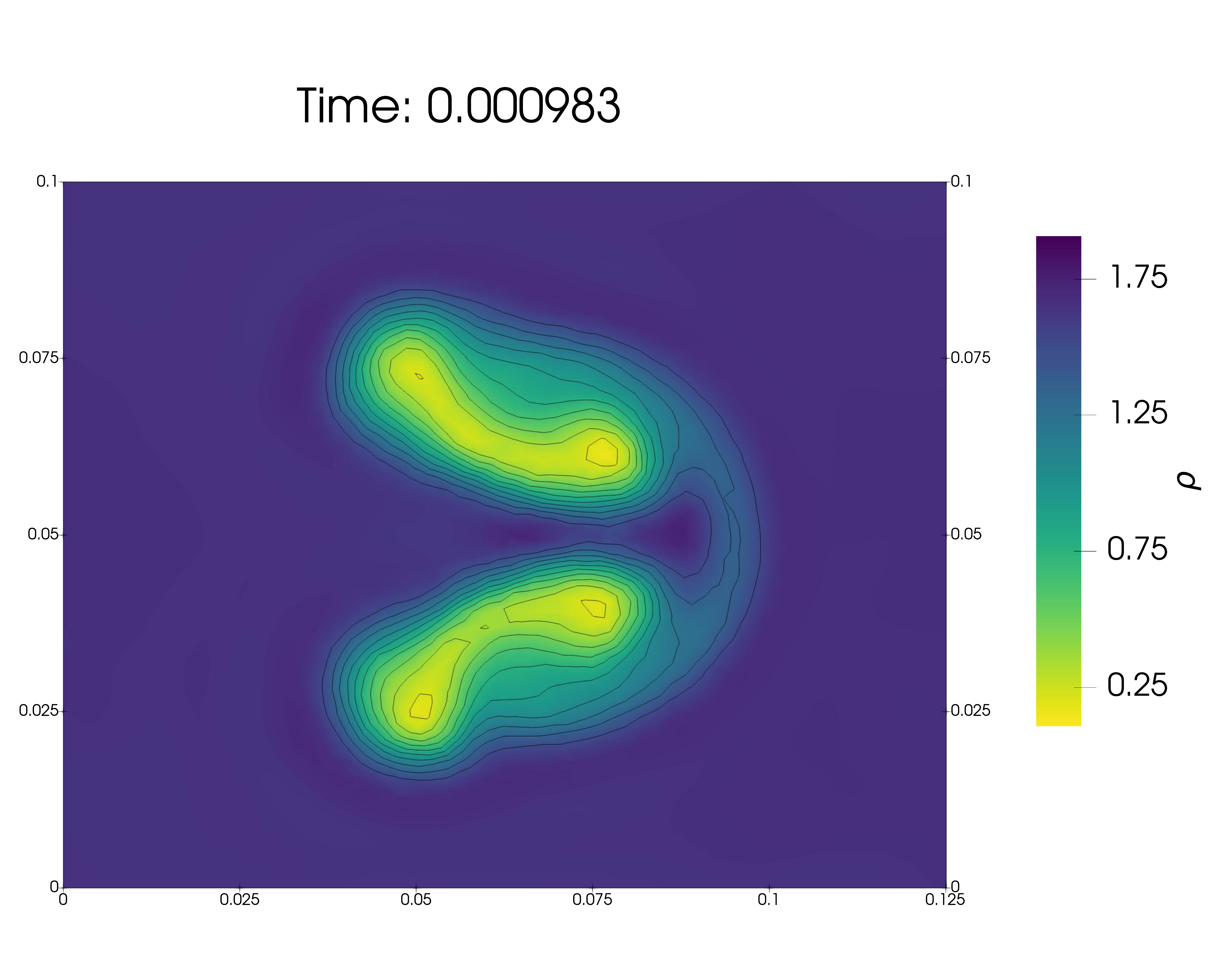}
        \includegraphics[width=\textwidth, trim = {4cm 3cm 0 2cm}, clip]{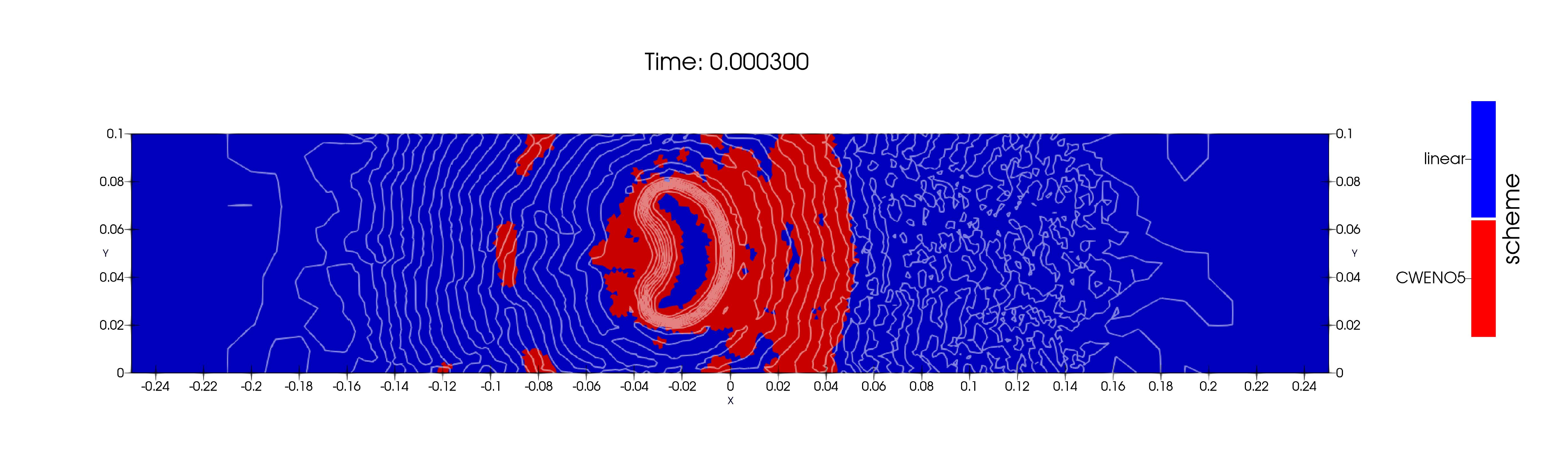}
        \caption{\label{fig:Bubble_MOODvsHybrid-Hybrid} MOOD-CWENO hybrid with entropy trigger}
    \end{subfigure}
    \caption{\label{fig:Bubble_MOODvsHybrid} Comparison of the results achieved with MOOD and MOOD-CWENO methods with intentionally increased number of flagged cells. Top: density - colour, volume fraction - contours at the final time $t=983\mu s$. Bottom: the flagged cells at an intermediate time $t=300\mu s$.}
\end{figure}

To emphasise the difference between the two methods, a contrived test was performed. The Helium Bubble test (\ref{HeliumBubble}) was repeated with the entropy production threshold $\beta$ reduced ten times, which results in excessively many cells being flagged. The results are presented in the Figure \ref{fig:Bubble_MOODvsHybrid}. The accuracy improvement due to using high-order scheme everywhere is immediately visible.
This improvement may not be applicable to simple test cases with perfectly adjusted MOOD trigger, but it is important nonetheless because it can save simulations where the trigger has not been configured correctly and it enables simulating very complex scenarios where adjusting the trigger perfectly is prohibitively hard or it is even impossible, because it must inevitably flag large number of cells.

\subsection{Computational cost}\label{Computational_cost}

\subsubsection{Timing}\label{Timing}

The computational cost of the methods were analysed by comparing the runtimes of simulations. The total time of simulation was measured including reading in the configuration, initialization, the main simulation loop and writing the output. The code and compilation optimizations were modified as little as possible between methods, in particular the input and output, the mesh handling and the parallelisation mechanisms were kept constant and the same as in the main branch UCNS3D, so they do not affect the comparison. 

To represent different use-cases the timing tests were performed on two different machines: a personal workstation and a large national HPC cluster. 
The workstation was equipped with 16 core Intel\textsuperscript{\textregistered} Core\textsuperscript{\texttrademark} i7-10700 and 64GB of RAM and running Linux Ubuntu 24.04. The simulations were run parallelised between all 16 cores.
The HPC cluster was the largest one in the UK - the Archer2 at the University of Edinburgh. The simulations were run on a 128-core partition with an aim of utilising one full node.
Each test combination was run four times and the runtimes were averaged to limit the impact of uncontrollable fluctuations associated with the load on the machines and performance of the scheduler.
\begin{table}[H]
  \begin{center}
  \caption{\label{tab:DMR_Archer2_classic_vs_entropy_timing} Timing of the Double Mach Reflection test case performed on the Archer2 supercomputer}
  \vspace{-1ex}
  \begin{tabular}{|l | c | c |} 
    \hline
    & classic MOOD & entropy MOOD\\
    \hline
    coarse & $10.92 \pm0.06$s & $9.72 \pm 0.02$s \\
    \hline
    medium & $66 \pm 5$s & $63.1 \pm 1.4$s \\ 
    \hline
    fine & $680 \pm 20$s & $ 520 \pm 30$s\\
    \hline
    super fine & $5150 \pm 130$s & $4280 \pm 110$s \\
    \hline
  \end{tabular}
  \end{center}
  \vspace{-2ex}
\end{table}

\begin{table}[H]
  \begin{center}
  \caption{\label{tab:ShuOsher_Workstation_MOOD_vs_hybrid_timing} Timing of the Shu-Osher test case performed on a personal workstation}
  \vspace{-1ex}
  \begin{tabular}{|l | c | c | c|} 
    \hline
    & CWENO & entropy MOOD & MOOD-CWENO hybrid with entropy trigger\\
    \hline
    coarse: & $20.6 \pm 0.6$s & $14.9 \pm 0.9$s & $16.3 \pm 0.2$s \\
    \hline
    medium: & $232 \pm 14$s & $149 \pm 3$s & $169 \pm 5$s \\
    \hline
    fine: & $1505 \pm 45$s & $1243 \pm 7$s & $1288 \pm 5$s \\ 
    \hline
  \end{tabular}  
  \end{center}
  \vspace{-2ex}
\end{table}
The results of the test on the personal workstation are presented in the Table \ref{tab:ShuOsher_Workstation_MOOD_vs_hybrid_timing} and the Archer2 results are shown in the Table \ref{tab:DMR_Archer2_classic_vs_entropy_timing}.

\subsubsection{Triggers comparison}

As seen in Tables \ref{tab:DMR_Archer2_classic_vs_entropy_timing} and \ref{tab:ShuOsher_Workstation_MOOD_vs_hybrid_timing} MOOD with entropy production trigger is up to $10\%$ faster than equivalent simulation using classic trigger.

\begin{figure}[!h]
    \centering
    \begin{subfigure}{0.49\linewidth}
        \centering
        \includegraphics[width=\linewidth]{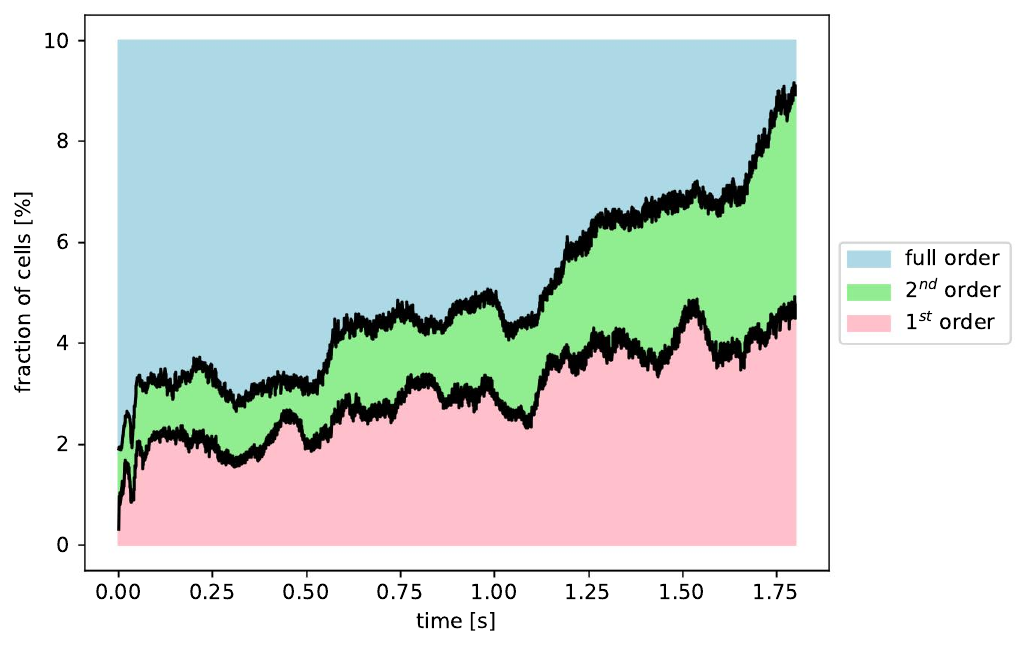}
        \vspace{-3.5ex}
        \caption{\label{fig:ShuOsher_ClassicMOOD_OrderFraction} Original MOOD}
    \end{subfigure}
    \begin{subfigure}{0.49\linewidth}
        \centering
        \includegraphics[width=\linewidth]{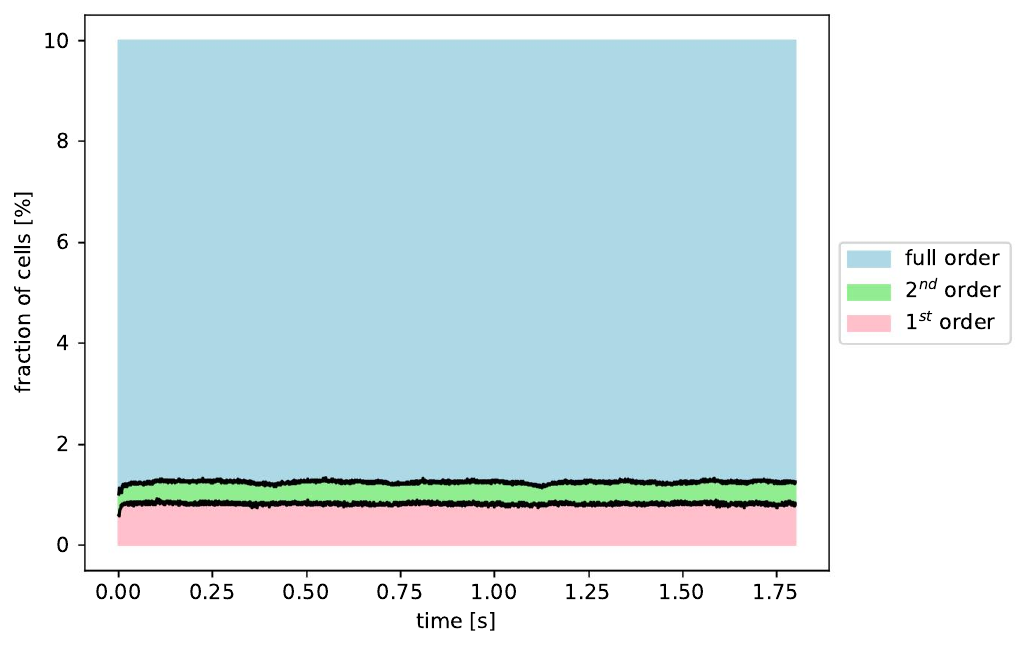}
        \vspace{-3.5ex}
        \caption{\label{fig:ShuOsher_EntropyMOOD_OrderFraction} MOOD entropy production trigger}
    \end{subfigure}
    \begin{subfigure}{0.49\linewidth}
        \centering
        \includegraphics[width=\linewidth]{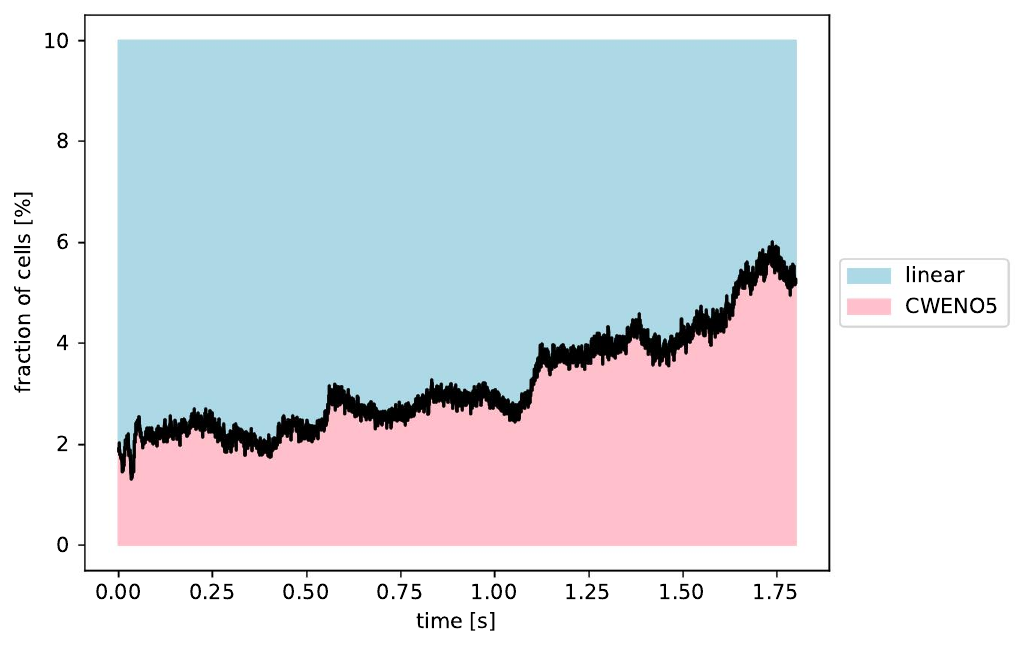}
        \vspace{-3.5ex}
        \caption{\label{fig:ShuOsher_CWENO-MOODclassic_OrderFraction} MOOD-CWENO classic trigger}
    \end{subfigure}
    \begin{subfigure}{0.49\linewidth}
        \centering
        \includegraphics[width=\linewidth]{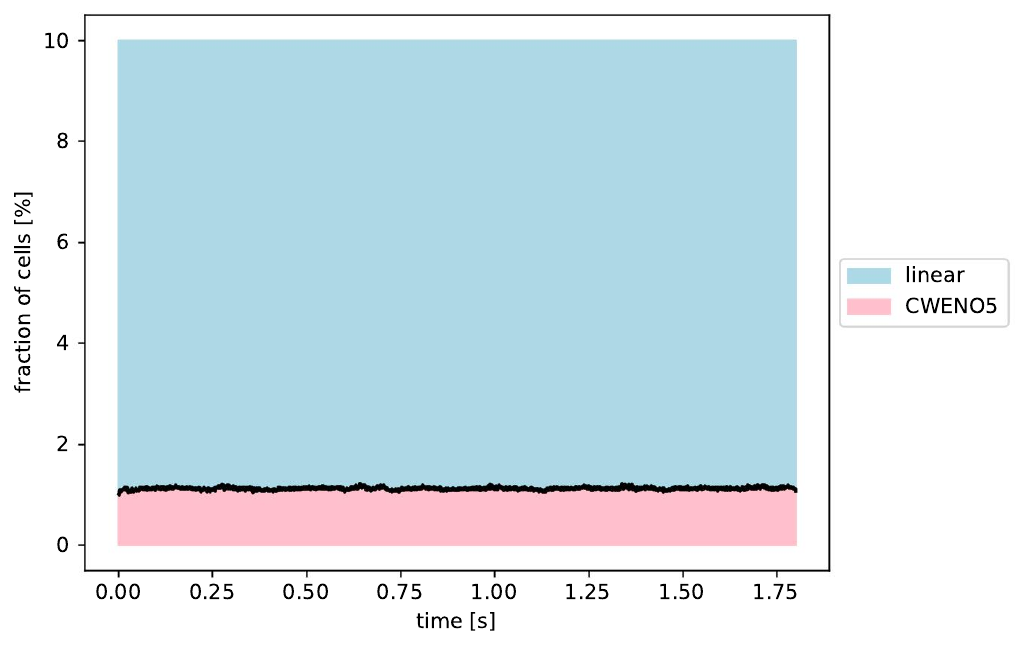}
        \vspace{-3.5ex}
        \caption{\label{fig:ShuOsher_CWENO-MOODentropy_OrderFraction} \centering MOOD-CWENO entropy trigger}
    \end{subfigure}
    \caption{\label{fig:ShuOsher_orderFraction} Fraction of cells flagged by the MOOD's decision algorithms during simulation of the Shu-Osher test problem on the fine mesh. For MOOD the bottom line shows the fraction of cells which have been flagged in both passes and the top one only in the first pass. Hence the height of the green band is the fraction of cells that have been updated using $2^{nd}$ order scheme and the red band is the fraction of cells that have been updated using $1^{st}$ order scheme.}
\end{figure}

\begin{figure}[!h]
    \centering
    \includegraphics[width=0.66\textwidth]{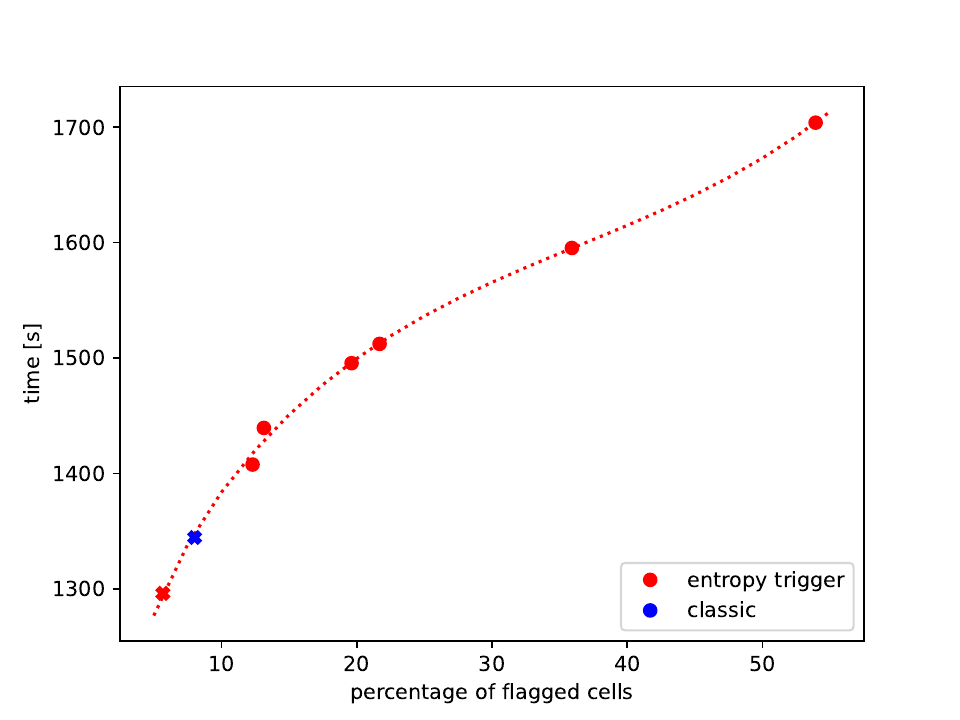}
    \caption{\label{fig:TimeVsCellFraction} Effects of fraction of flagged cells,represented here by the percentage in the last time step, on the simulation time. Data from medium mesh in Figure \ref{fig:resolution_scaling} compared with equivalent simulation using classic MOOD trigger.}
\end{figure}

The overall computations cost of MOOD methods depends both on the complexity of the decision algorithm and the number of recomputed cells. As seen in the Figure \ref{fig:ShuOsher_orderFraction}, the entropy production trigger has an advantage when it comes to the number of flagged cells, because it can be dialled in to flag only the necessary minimum of cells inside the shocks.
\FloatBarrier
Comparing the complexity of the decision algorithms is more involved, because different algorithms inherently result in different fractions of flagged cells, so to separate the computational cost of the trigger the data from Figure \ref{fig:resolution_scaling} is reanalysed in the Figure \ref{fig:TimeVsCellFraction} and compared with the 
classic MOOD trigger.
The simulation using relaxed classic trigger has runtime exactly in line with entropy trigger simulations with similar fraction of flagged cells, so it can be concluded that both decision algorithms have practically identical computational cost.

\subsubsection{MOOD-CWENO hybrid}

As seen in Table \ref{tab:ShuOsher_Workstation_MOOD_vs_hybrid_timing} using MOOD-CWENO hybrid gives up to $20\%$ speed up over classic CWENO5 and optimising communication algorithms specifically for this algorithm (all tests were performed using the same communication subroutines originally developed in UCNS3D for traditional high-order nonlinear schemes) could potentially further increase the difference.

Pure MOOD which uses lower order schemes in flagged cells is faster then the proposed MOOD-CWENO hybrid. This is inevitable because the full CWENO5 high-order nonlinear reconstruction is a much more involved algorithm, than a slop-limited second order reconstruction or the first-order strategy of using cell-centre data. However, in many scenarios, it can compensate for the slower performance by improved accuracy and robustness.

\FloatBarrier
\section{Conclusions}\label{Conclusions}

In this work, we introduced \textit{eMOOD}, an entropy-guided extension of the MOOD framework for high-order finite-volume methods on unstructured meshes. The proposed troubled-cell indicator is based on numerical entropy production and provides a compact, physically motivated alternative to the classical DMP-based MOOD criterion. By directly targeting shocks, contact discontinuities, and under-resolved flow structures, the method significantly improves the localisation of troubled regions while avoiding excessive activation in smooth areas.

Across the benchmark problems considered in this work, the entropy-based detector reduced the number of flagged cells by up to 75\%, leading to improved computational efficiency while preserving the robustness and non-oscillatory properties of the original MOOD framework. Owing to its compact stencil and fully local formulation, the proposed approach also improves memory-access locality and reduces communication overhead, making it well suited for massively parallel simulations on modern CPU and GPU architectures of heterogeneous systems.

In addition, we proposed a hybrid MOOD-CWENO strategy that accelerates WENO-type reconstructions through an \textit{a posteriori} adaptive switching mechanism between inexpensive linear reconstruction and full nonlinear CWENO blending. Compared to the reference CWENO implementation, the proposed hybrid approach achieved speed-ups of up to 25\% while maintaining comparable accuracy and robustness. Unlike conventional MOOD methods, the hybrid formulation preserves high-order reconstruction even within flagged regions, improving its ability to resolve fine-scale flow structures and reducing sensitivity to imperfections in the troubled-cell detector.

The proposed methods were implemented within the open-source high-order parallel UCNS3D CFD framework and validated on a diverse range of challenging compressible-flow benchmarks involving strong shock interactions, shock--entropy-wave coupling, vortex generation, and multi-material interface dynamics. Overall, the results demonstrate that entropy-guided MOOD strategies provide a robust and computationally efficient framework for next-generation high-order CFD simulations on unstructured meshes.

Future work will focus on extending the proposed methodology to direct and indirect Arbitrary Lagrangian-Eulerian (ALE) formulations on high-order unstructured meshes in two and three dimensions for compressible-flow applications, with particular emphasis on moving-mesh robustness, scalability, and complex multi-material interactions.

\section*{Acknowledgements}
The authors acknowledge the computing time on ARCHER2 through the UK Turbulence Consortium [EPSRC grant EP/X035484/1]. The authors acknowledge the support from EPSRC-CASE award [Grant No. 2924903] and AWE plc.

\section*{Data Availability}
The datasets of the test problems in this article will be made available in an open-access repository once the manuscript is accepted.

\bibliography{bibtex}

@ARTICLE{Shu198932,
    author = {Shu, C.-W. and Osher, S.},
    title  = {Efficient implementation of essentially non-oscillatory shock-capturing schemes, {I}{I}},
    journal= {Journal of Computational Physics},
    year   = {1989},
    volume = {83},
    number = {1},
    pages  = {32-78},
    doi    = {10.1016/0021-9991(89)90222-2}
}

@ARTICLE{Woodward1984115,
    author = {Woodward, P. and Colella, P.},
    title  = {The numerical simulation of two-dimensional fluid flow with strong shocks},
    journal= {Journal of Computational Physics},
    year   = {1984},
    volume = {54},
    number = {1},
    pages  = {115-173},
    doi    = {10.1016/0021-9991(84)90142-6},
}

@article{HeliumBubblePaperExperimental,
    author = {J-F Haas and B. Sturtevant},
    title  = {Interaction of weak shock waves with cylindrical and spherical gas inhomogeneities. Journal of Fluid Mechanics},
    journal= {Journal of Fluid Mechanics},
    year   = {1987},
    doi    = {10.1017/S0022112087002003},
    pages  = {41-76},
    volume = {181}
}

@article{HeliumBubblePaperSimulation,
    author = {J. M. Picone and J. P.Boris},
    title  = {Vorticity generation by shock propagation through bubbles in a gas},
    journal= {Journal of Fluid Mechanics},
    year   = {1988},
    doi    = {10.1017/S0022112088000904},
    pages  = {23-51},
    volume = {189}
}

@article{TriplePointPaper,
    title  = {High-order curvilinear finite elements for axisymmetric Lagrangian hydrodynamics},
    journal= {Computers \& Fluids},
    volume = {83},
    pages  = {58-69},
    year   = {2013},
    doi    = {10.1016/j.compfluid.2012.06.004},
    author = {Veselin A. Dobrev and Truman E. Ellis and Tzanio V. Kolev and Robert N. Rieben},
}

@article{SharkEntropyPaper,
    author = {Matteo Semplice and Alessandra Zappa},
    title  = {Numerical entropy production in finite volume {$P_0P_M$} {ADER} schemes},
    journal= {ArXiv},
    year   = {2025},
    doi    = {10.48550/arXiv.2501.05872}
}

@article{HLL_OriginalPaper,
    author = {Harten, Amiram and Lax, Peter D. and van Leer, Bram},
    title  = {On Upstream Differencing and Godunov-Type Schemes for Hyperbolic Conservation Laws},
    journal= {SIAM Review},
    volume = {25},
    number = {1},
    pages  = {35-61},
    year   = {1983},
    doi    = {10.1137/1025002}
}

@article{HLLC_OriginalPaper,
    author = {Toro, EF. and Chakraborty, A.},
    title  = {The development of a Riemann solver for the steady supersonic Euler equations.},
    journal= {The Aeronautical Journal},
    year   = {1994},
    doi    = {10.1017/S0001924000026890},
    pages  = {325-339},
    volume = {98},
    issue  = {979}
}

@ARTICLE{TsoutsanisMood,
    author={Farmakis, P.S. and Tsoutsanis, P. and Nogueira, X.},
    title={{W}{E}{N}{O} schemes on unstructured meshes using a relaxed a posteriori {M}{O}{O}{D} limiting approach},
    journal={Computer Methods in Applied Mechanics and Engineering},
    year={2020},
    volume={363},
    doi={10.1016/j.cma.2020.112921},
    art_number={112921},
}

@article{MOOD_OriginalPaper1,
   author = {S. Clain and S. Diot and R. Loubère},
   doi = {10.1016/j.jcp.2011.02.026},
   issn = {00219991},
   issue = {10},
   journal = {Journal of Computational Physics},
   month = {5},
   pages = {4028-4050},
   title = {{A high-order finite volume method for systems of conservation laws—Multi-dimensional Optimal Order Detection (MOOD)}},
   volume = {230},
   year = {2011}
}

@article{MOOD_OriginalPaper2,
    title  = {{An Overview on the Multidimensional Optimal Order Detection Method}},
    author = {Stéphane Clain and Jorge Figueiredo and Raphael Loubère and Steven Diot},
    journal= {SYMCOMP 2015 Faro},
    year   = {2015},
    doi    = {10.1093/gji/ggy472}
}

@article{Diot2013,
   title = {{The Multidimensional Optimal Order Detection method in the three-dimensional case: very high-order finite volume method for hyperbolic systems}},
   author = {Diot, S. and Loub{\`{e}}re, R. and Clain, S.},
   doi = {10.1002/fld.3804},
   journal = {International Journal for Numerical Methods in Fluids},
   number = {4},
   pages = {362--392},
   volume = {73},
   year = {2013}
}

@article{Dumbser2014,
    title = {{A posteriori subcell limiting of the discontinuous Galerkin finite element method for hyperbolic conservation laws}},
    author = {Dumbser, Michael and Zanotti, Olindo and Loub{\`{e}}re, Rapha{\"{e}}l and Diot, Steven},
    doi = {10.1016/j.jcp.2014.08.009},
    issn = {10902716},
    journal = {Journal of Computational Physics},
    pages = {47--75},
    publisher = {Elsevier Inc.},
    volume = {278},
    year = {2014}
}

@article{Vilar2019,
    title = {{A posteriori correction of high-order discontinuous {G}alerkin scheme through subcell finite volume formulation and flux reconstruction}},
    author = {Vilar, Fran{\c{c}}ois},
    doi = {10.1016/j.jcp.2018.10.050},
    issn = {10902716},
    journal = {Journal of Computational Physics},
    pages = {245-279},
    publisher = {Elsevier Inc.},
    volume = {387},
    year = {2019}
}

@article{Fidalgo2018,
    title = {{An a posteriori, efficient, high-spectral resolution hybrid finite-difference method for compressible flows}},
    author = {Fern{\'{a}}ndez-Fidalgo, Javier and Nogueira, Xes{\'{u}}s and Ram{\'{i}}rez, Luis and Colominas, Ignasi},
    doi = {10.1016/j.cma.2018.02.013},
    issn = {00457825},
    journal = {Computer Methods in Applied Mechanics and Engineering},
    pages = {91-127},
    volume = {335},
    year = {2018}
}

@article{Tsoutsanis_stencil_selection,
    author  = {Panagiotis Tsoutsanis},
    title   = {{Stencil selection algorithms for WENO schemes on unstructured meshes}},
    journal = {Journal of Computational Physics: X},
    volume  = {4},
    pages   = {100037},
    year    = {2019},
    doi     = {10.1016/j.jcpx.2019.100037}
}

@article{WENO_OriginalPaper,
    title   = {{Weighted Essentially Non-oscillatory Schemes}},
    author  = {Xu-Dong Liu and Stanley Osher and Tony Chan},
    journal = {Journal of Computational Physics},
    volume  = {115},
    number  = {1},
    pages   = {200-212},
    year    = {1994},
    doi     = {10.1006/jcph.1994.1187},
    issn    = {0021-9991}
}

@article{MUSCL_OriginalPaper,
    title   = {{Towards the ultimate conservative difference scheme. V. A second-order sequel to Godunov's method}},
    author  = {Bram {van Leer}},
    journal = {Journal of Computational Physics},
    volume  = {32},
    number  = {1},
    pages   = {101-136},
    year    = {1979},
    issn    = {0021-9991},
    doi     = {10.1016/0021-9991(79)90145-1}
}

@article{Dumbser_2017,
    title={{Central Weighted ENO Schemes for Hyperbolic Conservation Laws on Fixed and Moving Unstructured Meshes}},
    author={Dumbser, Michael and Boscheri, Walter and Semplice, Matteo and Russo, Giovanni},
    volume={39},
    doi={10.1137/17m1111036},
    number={6},
    journal={SIAM Journal on Scientific Computing},
    publisher={Society for Industrial & Applied Mathematics (SIAM)},
    year={2017},
    month = {1},
    pages={A2564–A2591}
}

@phdthesis{GodunovPhD,
    author = {Godunov, S.},
    title = {{Different Methods for Shock Waves}},
    school = {Moscow State University},
    year = {1954}
}

@article{GodunovMethod_OriginalPaper,
    author = {Sergei Godunov},
    title = {{A Difference Scheme for Numerical Solution of Discontinuous Solution of Hydrodynamic Equations}},
    journal = {Mat. Sb.},
    year = {1959},
    volume = {47(89)},
    pages = {271-306}
}

@article{Harten1983,
    title = {High resolution schemes for hyperbolic conservation laws},
    author = {Ami Harten},
    journal = {Journal of Computational Physics},
    volume = {49},
    number = {3},
    pages = {357-393},
    year = {1983},
    issn = {0021-9991},
    doi = {https://doi.org/10.1016/0021-9991(83)90136-5}
}

@article{Tsoutsanis_2018_bounds_limiter,
    author  = {Panagiotis Tsoutsanis},
    title   = {{Extended bounds limiter for high-order finite-volume schemes on unstructured meshes}},
    journal = {Journal of Computational Physics},
    volume  = {362},
    pages   = {69--94},
    year    = {2018},
    doi     = {10.1016/j.jcp.2018.02.045}
}

@article{Tsoutsanis_Dumbser_CWENO,
    author  = {Panagiotis Tsoutsanis and Michael Dumbser},
    title   = {{Arbitrary high order central non-oscillatory schemes on mixed-element unstructured meshes}},
    journal = {Computers and Fluids},
    volume  = {225},
    pages   = {104993},
    year    = {2021},
    doi     = {10.1016/j.compfluid.2021.104993}
}

@article{Tsoutsanis_Adebayo_Skote_2021,
    author  = {Panagiotis Tsoutsanis and Ebenezer M. Adebayo and Alberto Carriba Merino and Antonio Perez Arjona and Martin Skote},
    title   = {{CWENO finite-volume interface capturing schemes for multicomponent flows using unstructured meshes}},
    journal = {Journal of Scientific Computing},
    volume  = {89},
    number  = {3},
    year    = {2021},
    doi     = {10.1007/s10915-021-01671-z}
}

@article{Maltsev_Skote_Tsoutsanis_dg_fv,
    author  = {Vladimir Maltsev and Martin Skote and Panagiotis Tsoutsanis},
    title   = {{High-order hybrid DG-FV framework for compressible multi-fluid problems on unstructured meshes}},
    journal = {Journal of Computational Physics},
    volume  = {502},
    pages   = {112804},
    year    = {2024}
}

@book{ToroTextbook,
    author = {Eleuterio Toro},
    title = {{Riemann Solvers and Numerical Methods for Fluid Dynamics}},
    publisher = {Springer},
    year = {2009}
}

@article{SSP_Runge-Kutta,
   title = {{On high order strong stability preserving runge-kutta and multi step time discretizations}},
   author = {Gottlieb, Sigal},
   date = {2005/11/01},
   doi = {10.1007/BF02728985},
   isbn = {1573-7691},
   journal = {Journal of Scientific Computing},
   number = {1},
   pages = {105--128},
   volume = {25},
   year = {2005},
}

@article{Ioriatti2019,
   title = {{A posteriori sub-cell finite volume limiting of staggered semi-implicit Discontinuous Galerkin schemes for the shallow water equations}},
   author = {Matteo Ioriatti and Michael Dumbser},
   doi = {10.1016/j.apnum.2018.08.018},
   issn = {01689274},
   journal = {Applied Numerical Mathematics},
   month = {1},
   pages = {443-480},
   publisher = {Elsevier B.V.},
   volume = {135},
   year = {2019}
}

@article{Qing2023,
   title = {{A cell-centered Discontinuous Galerkin multi-material Arbitrary Lagrangian-Eulerian method in axisymmetric geometry}},
   author = {Fang Qing and Zupeng Jia and Liqi Liu},
   doi = {10.1016/j.jcp.2022.111745},
   issn = {10902716},
   journal = {Journal of Computational Physics},
   year = {2023},
   month = {1},
   publisher = {Academic Press Inc.},
   volume = {473}
}

@article{Gaburro2021,
   title = {{A Posteriori Subcell Finite Volume Limiter for General $P_NP_M$ Schemes: Applications from Gasdynamics to Relativistic Magnetohydrodynamics}},
   author = {Elena Gaburro and Michael Dumbser},
   doi = {10.1007/s10915-020-01405-8},
   issn = {15737691},
   issue = {3},
   journal = {Journal of Scientific Computing},
   month = {3},
   publisher = {Springer},
   volume = {86},
   year = {2021}
}

@article{Wang2025,
   title = {{An extremum properties (EP)-based discontinuous sensor and hybrid weighted essentially non-oscillatory scheme on tetrahedral meshes}},
   author = {Zhenming Wang and Jun Zhu and Yan Tan and Linlin Tian and Ning Zhao},
   doi = {10.1016/j.jcp.2025.113979},
   issn = {10902716},
   journal = {Journal of Computational Physics},
   month = {7},
   publisher = {Academic Press Inc.},
   volume = {533},
   year = {2025}
}

@article{Tong2026,
   title = {A hybrid finite-volume reconstruction framework for efficient high-order shock-capturing on unstructured meshes},
   author = {Yiren Tong and Panagiotis Tsoutsanis},
   doi = {10.1016/j.compfluid.2026.106988},
   issn = {00457930},
   journal = {Computers \& Fluids},
   month = {3},
   publisher = {Elsevier Ltd},
   volume = {308},
   year = {2026}
}

@article{ExaHyPE,
   title = {{ExaHyPE}: An engine for parallel dynamically adaptive simulations of wave problems},
   author = {Anne Reinarz and Dominic E. Charrier and Michael Bader and Luke Bovard and Michael Dumbser and Kenneth Duru and Francesco Fambri and Alice Agnes Gabriel and Jean Matthieu Gallard and Sven Köppel and Lukas Krenz and Leonhard Rannabauer and Luciano Rezzolla and Philipp Samfass and Maurizio Tavelli and Tobias Weinzierl},
   doi = {10.1016/J.CPC.2020.107251},
   issn = {0010-4655},
   journal = {Computer Physics Communications},
   month = {9},
   pages = {107251},
   publisher = {North-Holland},
   volume = {254},
   year = {2020}
}

@inproceedings{Popov2021,
   title = {{Efficient implementation of space-time adaptive ADER-DG finite element method with LST-DG predictor and a posteriori sub-cell WENO finite-volume limiting for simulation of non-stationary compressible multicomponent reactive flows}},
   author = {Ivan S. Popov},
   doi = {10.1088/1742-6596/1740/1/012059},
   issn = {17426596},
   issue = {1},
   booktitle = {Journal of Physics: Conference Series},
   month = {1},
   publisher = {IOP Publishing Ltd},
   volume = {1740},
   year = {2021}
}

@phdthesis{Ioriatti_PhD,
    author = {Matteo Ioriatti},
    title = {{Semi-implicit schemes for compressible fluids in elastic pipes and a-posteriori sub-cell finite volume limiting techniques for semi-implicit Discontinuous Galerkin schemes for hyperbolic conservation laws on staggered meshes}},
    school = {University of Trento},
    year = {2018},
}

@article{Cirrottola2021,
   title = {{Adaptive deformation of 3D unstructured meshes with curved body fitted boundaries with application to unsteady compressible flows}},
   author = {Luca Cirrottola and Mario Ricchiuto and Algiane Froehly and Barbara Re and Alberto Guardone and Giuseppe Quaranta},
   doi = {10.1016/j.jcp.2021.110177},
   issn = {10902716},
   journal = {Journal of Computational Physics},
   month = {5},
   publisher = {Academic Press Inc.},
   volume = {433},
   year = {2021}
}

@article{BORGES20083191,
   title = {{An improved weighted essentially non-oscillatory scheme for hyperbolic conservation laws}},
   author = {Rafael Borges and Monique Carmona and Bruno Costa and Wai Sun Don},
   journal = {Journal of Computational Physics},
   volume = {227},
   number = {6},
   pages = {3191-3211},
   year = {2008},
   issn = {0021-9991},
   doi = {10.1016/j.jcp.2007.11.038}
}

@article{Castro2011,
   title = {{High order weighted essentially non-oscillatory WENO-Z schemes for hyperbolic conservation laws}},
   author = {Marcos Castro and Bruno Costa and Wai Sun Don},
   doi = {10.1016/J.JCP.2010.11.028},
   issn = {0021-9991},
   issue = {5},
   journal = {Journal of Computational Physics},
   month = {3},
   pages = {1766-1792},
   publisher = {Academic Press},
   volume = {230},
   year = {2011}
}

@article{Ioriatti2020,
   title = {{A Staggered Semi-implicit Discontinuous Galerkin Scheme with a Posteriori Subcell Finite Volume Limiter for the Euler Equations of Gasdynamics}},
   author = {Matteo Ioriatti and Michael Dumbser and Raphaël Loubère},
   doi = {10.1007/s10915-020-01209-w},
   issn = {15737691},
   issue = {2},
   journal = {Journal of Scientific Computing},
   month = {5},
   publisher = {Springer},
   volume = {83},
   year = {2020}
}

@article{MOOD_OriginalPaper3,
   title = {{Improved detection criteria for the Multi-dimensional Optimal Order Detection (MOOD) on unstructured meshes with very high-order polynomials}},
   author = {S. Diot and S. Clain and R. Loubère},
   doi = {10.1016/j.compfluid.2012.05.004},
   issn = {00457930},
   journal = {Computers \& Fluids},
   month = {7},
   pages = {43-63},
   publisher = {Elsevier Ltd},
   volume = {64},
   year = {2012}
}

@article{Nogueira2016,
   title  = {{High-accurate SPH method with Multidimensional Optimal Order Detection limiting}},
   author = {Xesús Nogueira and Luis Ramírez and Stéphane Clain and Raphaël Loubère and Luis Cueto-Felgueroso and Ignasi Colominas},
   doi    = {10.1016/j.cma.2016.06.032},
   issn   = {00457825},
   journal= {Computer Methods in Applied Mechanics and Engineering},
   year   = {2016},
   month  = {10},
   pages  = {134-155},
   volume = {310}
}

@article{TENO_Review2023,
   title = {{Review of the High-Order TENO Schemes for Compressible Gas Dynamics and Turbulence}},
   author = {Lin Fu},
   doi = {10.1007/S11831-022-09877-7},
   issn = {18861784},
   issue = {4},
   journal= {Archives of Computational Methods in Engineering},
   month = {5},
   pages = {2493-2526},
   volume = {30},
   year = {2023}
}

@article{Farmakis2020,
   title = {{WENO schemes on unstructured meshes using a relaxed a posteriori MOOD limiting approach}},
   author = {Pericles S. Farmakis and Panagiotis Tsoutsanis and Xesús Nogueira},
   doi = {10.1016/J.CMA.2020.112921},
   issn = {0045-7825},
   journal = {Computer Methods in Applied Mechanics and Engineering},
   month = {5},
   pages = {112921},
   publisher = {North-Holland},
   volume = {363},
   year = {2020}
}

@misc{Nguyen2026,
   title = {A Conservative Runge–Kutta Discontinuous Galerkin ConRKDG Method for Inviscid Compressible Flows in One-Dimensional Computational Fluid Dynamics Simulations},
   author = {Thien Binh Nguyen and Nguyen Minh Hieu Pham},
   doi = {10.20944/preprints202603.1611.v1},
   month = {3},
   year = {2026}
}

@article{Pandey2022,
   title = {{High-resolution WENO schemes using local variation-based smoothness indicator}},
   author = {Prashant Kumar Pandey and Farzad Ismail and Ritesh Kumar Dubey},
   doi = {10.1007/s40314-022-01916-0},
   issn = {18070302},
   issue = {5},
   journal = {Computational and Applied Mathematics},
   month = {7},
   publisher = {Springer Science and Business Media Deutschland GmbH},
   volume = {41},
   year = {2022}
}

@article{Fernndez-Fidalgo2018,
   title = {An a posteriori, efficient, high-spectral resolution hybrid finite-difference method for compressible flows},
   author = {Javier Fernández-Fidalgo and Xesús Nogueira and Luis Ramírez and Ignasi Colominas},
   doi = {10.1016/j.cma.2018.02.013},
   issn = {00457825},
   journal = {Computer Methods in Applied Mechanics and Engineering},
   month = {6},
   pages = {91-127},
   publisher = {Elsevier B.V.},
   volume = {335},
   year = {2018}
}

@article{Henrick2005,
   title = {{Mapped weighted essentially non-oscillatory schemes: Achieving optimal order near critical points}},
   author = {Andrew K. Henrick and Tariq D. Aslam and Joseph M. Powers},
   doi = {10.1016/J.JCP.2005.01.023},
   issn = {0021-9991},
   issue = {2},
   journal = {Journal of Computational Physics},
   year = {2005},
   month = {8},
   pages = {542-567},
   publisher = {Academic Press},
   volume = {207}
}

@article{Fernndez-Fidalgo2019,
   title = {An a posteriori very efficient hybrid method for compressible flows},
   author = {Javier Fernández-Fidalgo and Xesús Nogueira and Luis Ramírez and Ignasi Colominas},
   doi = {10.1007/978-3-030-11887-7_13},
   issn = {21959870},
   journal = {Springer Tracts in Mechanical Engineering},
   pages = {137-148},
   publisher = {Springer International Publishing},
   year = {2019}
}

@article{Bogey2004,
   title = {A family of low dispersive and low dissipative explicit schemes for flow and noise computations},
   author = {Christophe Bogey and Christophe Bailly},
   doi = {10.1016/J.JCP.2003.09.003},
   issn = {0021-9991},
   issue = {1},
   journal = {Journal of Computational Physics},
   month = {2},
   pages = {194-214},
   publisher = {Academic Press},
   volume = {194},
   year = {2004}
}

@article{Kozhanova2025,
   title = {{A hybrid a posteriori MOOD limited lattice Boltzmann method to solve compressible fluid flows – LBMOOD}},
   author = {Ksenia Kozhanova and Song Zhao and Raphaël Loubère and Pierre Boivin},
   doi = {10.1016/J.JCP.2024.113570},
   issn = {0021-9991},
   journal = {Journal of Computational Physics},
   month = {1},
   pages = {113570},
   publisher = {Academic Press},
   volume = {521},
   year = {2025}
}

@article{Eirs2021,
   title = {{SPH-ALE} scheme for weakly compressible viscous flow with a posteriori stabilization},
   author = {Antonio Eirís and Luis Ramírez and Javier Fernández-Fidalgo and Iván Couceiro and Xesús Nogueira},
   doi = {10.3390/w13030245},
   issn = {20734441},
   issue = {3},
   journal = {Water (Switzerland)},
   month = {2},
   publisher = {DPI AG},
   volume = {13},
   year = {2021}
}

@article{Jiang2023,
   title = {A simple a posteriori indicator for discontinuous Galerkin method on unstructured grids},
   author = {Zhen Hua Jiang and Chao Yan and Jian Yu},
   doi = {10.1007/s10409-022-22296-x},
   issn = {16143116},
   issue = {5},
   journal = {Acta Mechanica Sinica/Lixue Xuebao},
   month = {5},
   publisher = {Springer Verlag},
   volume = {39},
   year = {2023}
}

@article{Huang2023,
   title = {{Multidimensional WENO-AO Reconstructions Using a Simplified Smoothness Indicator and Applications to Conservation Laws}},
   author = {Chieh Sen Huang and Todd Arbogast and Chenyu Tian},
   doi = {10.1007/s10915-023-02319-x},
   issn = {15737691},
   issue = {1},
   journal = {Journal of Scientific Computing},
   year = {2023},
   month = {10},
   publisher = {Springer},
   volume = {97}
}

@article{WENO_Review2025,
   title = {{Recent Advancements in Fluid Flow Simulation Using the WENO Scheme: A Comprehensive Review}},
   author = {Reza Bozorgpour and Hossein Mahmoodi Darian},
   doi = {10.1007/s44198-025-00269-6},
   issn = {17760852},
   issue = {1},
   journal = {Journal of Nonlinear Mathematical Physics},
   month = {12},
   publisher = {Springer Nature},
   volume = {32},
   year = {2025}
}

@article{Zhang2007,
   title = {A new smoothness indicator for the WENO schemes and its effect on the convergence to steady state solutions},
   author = {Shuhai Zhang and Chi Wang Shu},
   doi = {10.1007/s10915-006-9111-y},
   issn = {08857474},
   issue = {1-2},
   journal = {Journal of Scientific Computing},
   month = {6},
   pages = {273-305},
   volume = {31},
   year = {2007}
}

@inproceedings{Tsoutsanis_ECCOMAS2022,
    author = {P Tsoutsanis and P S Farmakis},
    title = {A Posteriori {MOOD} Limiting Approach for Multicomponent Flows On Unstructured Meshes},
    booktitle = {ECCOMAS Congress},
    year = {2022},
    month = {06}
}

@article{Blanchard201683,
   title = {High order accurate conservative remapping scheme on polygonal meshes using a posteriori {MOOD} limiting},
   author = {Ghislain Blanchard and Raphaël Loub\'ere},
   journal = {Computers \& Fluids},
   volume = {136},
   pages = {83-103},
   year = {2016},
   issn = {0045-7930},
   doi = {10.1016/j.compfluid.2016.06.002}
}

@article{Thornber2008,
   title = {On entropy generation and dissipation of kinetic energy in high-resolution shock-capturing schemes},
   author = {B. Thornber and D. Drikakis and R. J.R. Williams and D. Youngs},
   doi = {10.1016/j.jcp.2008.01.035},
   issn = {10902716},
   issue = {10},
   journal = {Journal of Computational Physics},
   month = {5},
   pages = {4853-4872},
   publisher = {Academic Press Inc.},
   volume = {227},
   year = {2008}
}

@article{Boscheri2015,
   title = {{Direct Arbitrary-Lagrangian-Eulerian ADER-MOOD finite volume schemes for multidimensional hyperbolic conservation laws}},
   author = {Walter Boscheri and Raphaël Loubère and Michael Dumbser},
   doi = {10.1016/j.jcp.2015.03.015},
   issn = {10902716},
   journal = {Journal of Computational Physics},
   month = {7},
   pages = {56-87},
   publisher = {Academic Press Inc.},
   volume = {292},
   year = {2015}
}

@article{Semplice2018,
   title = {{Adaptive-Mesh-Refinement for hyperbolic systems of conservation laws based on a posteriori stabilized high order polynomial reconstructions}},
   author = {Matteo Semplice and Raphaël Loubère},
   doi = {10.1016/j.jcp.2017.10.031},
   issn = {00219991},
   journal = {Journal of Computational Physics},
   year = {2018},
   month = {2},
   pages = {86-110},
   volume = {354}
}

@article{USVevek2021,
   title = {{An efficient algorithm for computing smoothness indicators for WENO schemes}},
   author = {U S Vevek},
   doi = {10.1016/j.jcp.2021.110619},
   issn = {00219991},
   journal = {Journal of Computational Physics},
   month = {12},
   pages = {110619},
   volume = {446},
   year = {2021}
}

@article{Titarev2010,
   title = {{WENO Schemes for Mixed-Element Unstructured Meshes}},
   author = {Vladimir Titarev and Panagiotis Tsoutsanis and Dimitry Drikakis},
   doi = {10.4208/cicp.040909.080110a},
   issn = {1991-7120},
   issue = {3},
   journal = {Communications in Computational Physics},
   month = {1},
   pages = {585-609},
   volume = {8},
   year = {2010}
}

@article{Teng2011,
   title = {A speed-up strategy for finite volume {WENO} schemes for hyperbolic conservation laws},
   author = {Fei Teng and Li Yuan and Tao Tang},
   doi = {10.1007/s10915-010-9407-9},
   issn = {08857474},
   issue = {3},
   journal = {Journal of Scientific Computing},
   month = {3},
   pages = {359-378},
   volume = {46},
   year = {2011}
}

@article{Baeza2019,
   title = {{On the Efficient Computation of Smoothness Indicators for a Class of WENO Reconstructions}},
   author = {Antonio Baeza and Raimund Bürger and Pep Mulet and David Zorío},
   doi = {10.1007/s10915-019-00974-7},
   issn = {0885-7474},
   issue = {2},
   journal = {Journal of Scientific Computing},
   month = {8},
   pages = {1240-1263},
   volume = {80},
   year = {2019},
}

@article{Braeunig2018,
   title = {{A posteriori limiting for 2D Lagrange plus Remap schemes solving the hydrodynamics system of equations}},
   author = {Jean-Philippe Braeunig and Raphaël Loubère and Renaud Motte and Mathieu Peybernes and Raphaël Poncet},
   doi = {10.1016/j.compfluid.2017.08.020},
   issn = {00457930},
   journal = {Computers \& Fluids},
   month = {6},
   pages = {249-262},
   volume = {169},
   year = {2018}
}

@article{Bourgeois2022,
   title = {{GP-MOOD: A positivity-preserving high-order finite volume method for hyperbolic conservation laws}},
   author = {R. Bourgeois and D. Lee},
   doi = {10.1016/j.jcp.2022.111603},
   journal = {Journal of Computational Physics},
   volume = {471},
   year = {2022}
}

@article{Lee2020,
   title = {{GP-MOOD: a positivity-preserving high-order finite volume method for hyperbolic conservation laws}},
   author = {D. Lee and R. Bourgeois},
   doi = {10.1017/S1743921322001363},
   journal = {Proceedings of the International Astronomical Union},
   pages = {373-379},
   volume = {16},
   year = {2020}
}

@article{Boscheri2022,
   title = {{A 3D cell-centered ADER MOOD Finite Volume method for solving updated Lagrangian hyperelasticity on unstructured grids}},
   author = {Walter Boscheri and Raphaël Loubère and Pierre-Henri Maire},
   doi = {10.1016/j.jcp.2021.110779},
   issn = {00219991},
   journal = {Journal of Computational Physics},
   month = {1},
   pages = {110779},
   volume = {449},
   year = {2022}
}

@article{Tsoutsanis2023,
   title = {{A relaxed a posteriori MOOD algorithm for multicomponent compressible flows using high-order finite-volume methods on unstructured meshes}},
   author = {Panagiotis Tsoutsanis and Machavolu Sai Santosh Pavan Kumar and Pericles S. Farmakis},
   doi = {10.1016/J.AMC.2022.127544},
   issn = {0096-3003},
   journal = {Applied Mathematics and Computation},
   month = {1},
   pages = {127544},
   publisher = {Elsevier},
   volume = {437},
   year = {2023}
}

@article{Tsoutsanis2011,
   title = {{WENO} schemes on arbitrary mixed-element unstructured meshes in three space dimensions},
   author = {P. Tsoutsanis and V. A. Titarev and D. Drikakis},
   doi = {10.1016/J.JCP.2010.11.023},
   issn = {0021-9991},
   issue = {4},
   journal = {Journal of Computational Physics},
   month = {2},
   pages = {1585-1601},
   publisher = {Academic Press},
   volume = {230},
   year = {2011}
}

@article{Boscheri2017,
   title = {{Arbitrary-Lagrangian–Eulerian Discontinuous Galerkin schemes with a posteriori subcell finite volume limiting on moving unstructured meshes}},
   author = {Walter Boscheri and Michael Dumbser},
   doi = {10.1016/J.JCP.2017.06.022},
   issn = {0021-9991},
   journal = {Journal of Computational Physics},
   month = {10},
   pages = {449-479},
   publisher = {Academic Press},
   volume = {346},
   year = {2017}
}

@article{Michel-Dansac2022,
   title = {{TVD-MOOD} schemes based on implicit-explicit time integration},
   author = {Victor Michel-Dansac and Andrea Thomann},
   doi = {10.1016/J.AMC.2022.127397},
   issn = {0096-3003},
   journal = {Applied Mathematics and Computation},
   month = {11},
   pages = {127397},
   publisher = {Elsevier},
   volume = {433},
   year = {2022},
}

@article{Bacigaluppi2023,
   title = {{“A posteriori” limited high order and robust schemes for transient simulations of fluid flows in gas dynamics}},
   author = {Paola Bacigaluppi and Rémi Abgrall and Svetlana Tokareva},
   doi = {10.1016/J.JCP.2022.111898},
   issn = {0021-9991},
   journal = {Journal of Computational Physics},
   month = {3},
   pages = {111898},
   publisher = {Academic Press},
   volume = {476},
   year = {2023}
}

@article{Li2024,
   title = {{A family of TENOA-THINC-MOOD schemes based on diffuse-interface method for compressible multiphase flows}},
   author = {Qichao Li and Lin Fu},
   doi = {10.1016/J.JCP.2024.113375},
   issn = {0021-9991},
   journal = {Journal of Computational Physics},
   year = {2024},
   month = {12},
   pages = {113375},
   publisher = {Academic Press},
   volume = {519}
}

@article{VelascoRomero2023,
   title = {{Spectral difference method with a posteriori limiting: application to the Euler equations in one and two space dimensions}},
   author = {David A. Velasco Romero and Maria Han-Veiga and Romain Teyssier},
   doi = {10.1093/mnras/stad226},
   issn = {13652966},
   issue = {3},
   journal = {Monthly Notices of the Royal Astronomical Society},
   year = {2023},
   month = {4},
   pages = {3591-3608},
   publisher = {Oxford University Press},
   volume = {520}
}

@article{TENO_OriginalPaper,
   title = {{A family of high-order targeted ENO schemes for compressible-fluid simulations}},
   author = {Lin Fu and Xiangyu Y. Hu and Nikolaus A. Adams},
   doi = {10.1016/j.jcp.2015.10.037},
   issn = {10902716},
   journal = {Journal of Computational Physics},
   month = {1},
   pages = {333-359},
   publisher = {Academic Press Inc.},
   volume = {305},
   year = {2016}
}

@article{ENO_OriginalPaper,
    title = {Uniformly High Order Accurate Essentially Non-oscillatory Schemes, {III}},
    author = {Ami Harten and Bjorn Angquist and Stanley Osher and Sukumar Chakravarthv},
    journal = {Journal of Computational Physics},
    year = {1987},
    issue = {71},
    pages = {231-303}
}

@article{CWENO_OriginalPaper,
    title   = {{Central {WENO} Schemes for Hyperbolic Systems of Conservation Laws}},
    author  = {Doron, Levy and Gabriella, Puppo and Giovanni Russo},
    journal = {Mathematical Modelling and Numerical Analysis},
    year    = {1999},
    volume  = {33},
    pages   = {547-571},
    doi     = {10.1051/m2an:1999152}
}

@article{WENO-JS_OriginalPaper,
   title   = {{Efficient Implementation of Weighted ENO Schemes}},
   author  = {Guang-Shan Jiang and Chi-Wang Shu},
   journal = {Journal of Computational Physics},
   volume = {126},
   number = {1},
   pages = {202-228},
   year = {1996},
   issn = {0021-9991},
   doi = {10.1006/jcph.1996.0130}
}

@article{UCNS3D_Paper,
   title   = {{UCNS3D: An open-source high-order finite-volume unstructured CFD solver}},
   author  = {A. F. Antoniadis and D. Drikakis and P. S. Farmakis and L. Fu and I. Kokkinakis and X. Nogueira and P. A. S. F. Silva and M. Skote and V. Titarev and P. Tsoutsanis},
   journal = {Computer Physics Communications},
   year    = {2022},
   volume  = {279},
   doi     = {10.1016/j.cpc.2022.108453}
}

@article{HWENO_OriginalPaper1,
   title = {{Hermite WENO schemes and their application as limiters for Runge–Kutta discontinuous Galerkin method: one-dimensional case}},
   author = {Jianxian Qiu and Chi-Wang Shu},
   journal = {Journal of Computational Physics},
   volume = {193},
   number = {1},
   pages = {115-135},
   year = {2004},
   issn = {0021-9991},
   doi = {10.1016/j.jcp.2003.07.026}
}

@article{HWENO_OriginalPaper2,
   title = {{Multi-resolution HWENO schemes for hyperbolic conservation laws}},
   author = {Jiayin Li and Chi-Wang Shu and Jianxian Qiu},
   journal = {Journal of Computational Physics},
   volume = {446},
   pages = {110653},
   year = {2021},
   issn = {0021-9991},
   doi = {10.1016/j.jcp.2021.110653}
}

@article{HUANG2018498,
   title = {A simple smoothness indicator for the {WENO} scheme with adaptive order},
   author = {Cong Huang and Li Li Chen},
   journal = {Journal of Computational Physics},
   volume = {352},
   pages = {498-515},
   year = {2018},
   issn = {0021-9991},
   doi = {10.1016/j.jcp.2017.10.005}
}

@article{KUMAR20181059,
   title = {Simple smoothness indicator and multi-level adaptive order {WENO} scheme for hyperbolic conservation laws},
   author = {Rakesh Kumar and Praveen Chandrashekar},
   journal = {Journal of Computational Physics},
   volume = {375},
   pages = {1059-1090},
   year = {2018},
   issn = {0021-9991},
   doi = {10.1016/j.jcp.2018.09.027}
}

@article{Cravero2016,
   author = {I. Cravero and M. Semplice},
   title = {On the Accuracy of WENO and CWENO Reconstructions of Third Order on Nonuniform Meshes},
   doi = {10.1007/s10915-015-0123-3},
   journal = {Journal of Scientific Computing},
   issue = {3},
   month = {6},
   pages = {1219-1246},
   volume = {67},
   year = {2016}
}
\bibliographystyle{elsarticle-num}

\end{document}